# Gödel coding on fibrations and geminal categories

*Yuto Ikeda*
Graduate School of Mathematical Sciences, the University of Tokyo
yuto@ikeda.ac

## Abstract

In their 2023 dissertation, Ramesh has proposed two categorical notions, *introspective theories* and *geminal categories,* to unify various structures sharing the form of *Löb's theorem:* $\Box A \vdash A$ implies $\vdash A$. These include classical logical systems, Joyal's arithmetic universes, Kripke semantics for the provability logic GL, and categorical models for guarded recursion. The remarkable feature of Ramesh's approach is that it does not impose Löb's theorem as a definition; instead, it directly formalizes the "self-internalizing" nature, such as the formalization of a theory within itself as seen in the proof of the incompleteness theorems. Indeed, Ramesh demonstrates that Löb's theorem can be derived from these structures through highly non-trivial arguments.

In this thesis, we provide a mathematical and conceptual reorganization of the theory of geminal categories in a self-contained manner. As a consequence, we establish a novel categorical counterpart of the *Gödel–Löb axiom* $\Box(\Box A \to A) \to \Box A$ for any geminal category. We also provide a significant simplification of the proof of Löb's theorem for geminal categories, along with a slight generalization of the result.

To achieve such an organization, we introduce the notion of *code structures* on fibrations, which serves as a natural abstraction of Gödel coding, and establish the *fixed point theorem* for them, generalizing the classical diagonal lemma. This result acts as a universal theorem that unifies various classical results possessing *intensional* features. We then introduce *pre-geminal categories,* a slight weakening of geminal categories, as a natural structure that induces code structures on codomain fibrations. *Löb's theorem for pre-geminal categories* is then derived from our fixed point theorem, which simplifies and highlights the core idea of Ramesh's original proof. The notion of *geminal categories* and the *Gödel–Löb axiom* for them are realized by internalizing our arguments for pre-geminal categories, providing a new perspective on geminal categories via code structures. We also compare our framework with the structures appearing in modal calculi for metaprogramming.

Our formulation does not rely on the notion of introspective theories or any informal internal reasoning, making the theory of geminal categories more accessible to a broader audience. We propose that the framework of geminal categories offers a promising perspective for abstracting and unifying the interaction between meta- and object-levels arising in both logic and computer science.

## Acknowledgments

I am deeply grateful to my supervisor, Ryu Hasegawa, for his invaluable academic guidance and constant support. I would also like to thank my secondary supervisor, Toshitake Kohno, for helpful advice during our meetings.

I am sincerely grateful to Sridhar Ramesh for our insightful discussions and generous cooperation. These personal communications provided me with valuable insights and warm encouragement that motivated my research. Ramesh's pioneering dissertation also provided the foundational framework for this study and continues to be a profound source of inspiration.

I wish to express my special thanks to Hisashi Aratake for the kind encouragement since the start of my studies and for a suggestion regarding Joyal's arithmetic universe, which served as a turning point that guided me toward the topics explored in this thesis. I am also deeply grateful to Yuito Murase for introducing me to the extensive literature on modal calculi. I am also indebted to Kazuyuki Asada, Ryuya Hora, Keisuke Hoshino, Koshiro Ichikawa, Yuto Kawase, Taishi Kurahashi, Yutaka Maita, Akinori Maniwa, Satoshi Nakata, Hayato Nasu, Takako Nemoto, Christoph Schweigert, Haruki Toyota, Takeshi Tsukada and Taichi Uemura for their helpful discussions and encouragement, as well as for pointing out relevant literature.

I am also grateful to the (former) graduate students, Yuta Yamamoto, Yuto Arai, Ryuya Hora, Haruya Minoura, Sangwoo Kim, Yuki Kiyomori, Koki Kurahashi, Ryota Kuroki, Hiromasa Kondo, Haruki Toyota, and Minoru Sekiyama. Their insightful talks and comments in the graduate seminar were indispensable to the development of my ideas.

Finally, I am deeply grateful to my family for their unwavering support in my daily life. Their constant presence and encouragement have been a vital source of strength throughout my studies.

## Funding

This research was supported by Forefront Physics and Mathematics Program to Drive Transformation (FoPM), a World-leading Innovative Graduate Study (WINGS) Program, the University of Tokyo.

## Note on this version

This document is a slightly modified version of the master's thesis submitted to the University of Tokyo. The original submitted thesis is available at: https://ikeda.ac/pubs/ikeda2026/

# Contents

# 1 Introduction

## 1.1 Background

Despite its fundamental importance, Gödel's second incompleteness theorem remains mysterious. Unlike the first incompleteness theorem, the second theorem's statement depends on the specific construction of Gödel numbering and a provability predicate. This sensitivity to the choice of specific encoding, often called the problem of *intensionality* [6], obscures the essential structure of the theorem and its development.

Categorical logic may offer a promising framework to address this issue. Beginning with Lawvere's seminal work [25], categorical logic is founded on the idea of reformulating logical theories as categories rather than formal symbols, a shift that allows for extracting the signature-independent structure of the theories. In the same spirit, one may seek the essence of incompleteness by extracting the categorical structures common to "self-internalizing" systems where the second incompleteness and Löb's theorem hold. Were such a categorical structure established, it would serve as an essential structure governing the interaction between meta- and object-levels, thereby isolating the notorious encoding-dependent arguments.

There are also more concrete motivations for categorical analysis arising from computer science. In this field, computational principles analogous to the provability logic GL have been observed empirically. For instance, Löb induction in guarded recursion [2, 32] is known to correspond to the *strong Löb axiom*, $(\Box A \to A) \to A$. Even more direct is the connection suggested by Kavvos [21], where the *Gödel–Löb axiom* $\Box(\Box A \to A) \to \Box A$ corresponds to *intensional recursion* motivated by metaprogramming. These observations raise a question: why does provability logic, originating in pure logic, reappear so naturally in computational contexts? To provide an answer, we need a unified theory that extracts essential structures shared by both domains.

There are few categorical studies on the second incompleteness theorem, though some do exist. The most pioneering is Joyal's categorical interpretation of the theorem using *arithmetic universes*. Regrettably, Joyal's original work remains unpublished, although an attempt at its reconstruction can be found in van Dijk and Gietelink Oldenziel [8]. Joyal's key insight is to capture the "self-internalization" of a theory — where the theory is formalized within itself — through the concept of *internal categories*. Subsequent research has largely focused on the concept of an arithmetic universe itself, leading to its axiomatization as a list-arithmetic pretopos [27, 29, 31]. However, very few studies have focused on the original motivation of applying this framework to the incompleteness theorem. To the best of our knowledge, this direction has been pursued only by van Dijk and Gietelink Oldenziel [8] and by Ramesh [38], whom we shall discuss later. We should also mention that Joyal himself has revisited these ideas [17], though no full paper was published.

On the other hand, in computer science, several categorical structures have been proposed to analyze relevant computational systems. Most notably, Kavvos [20] intro-duces *Gödel–Löb categories* as models for the modal calculus DGL, which precisely corresponds to the provability logic GL. In the study of guarded recursion, *guarded fixpoint categories* [30] are known as a structure corresponding to the strong Löb axiom.

More recently, Ramesh [38] has proposed a general categorical framework in their 2023 Ph.D. dissertation, "Introspective theories and geminal categories". This framework provides a unified account of Gödel–Löb phenomena across several different structures: the syntactic categories of logical theories, Joyal's arithmetic universes, the Kripke semantics of provability logic, and the topos of trees, which serves as a categorical model for guarded recursion. In contrast to existing categorical models of modal calculi, the remarkable feature of Ramesh's framework is that it does not impose Löb's theorem or the Gödel–Löb axiom as a part of the definitions. Instead, it directly formalizes self-internalizing structures through the use of internal categories, in line with Joyal's pioneering insight. In fact, Ramesh demonstrates through a highly non-trivial mathematical argument that Löb's theorem can be derived from these structures. We believe that Ramesh's theory offers a promising perspective on unifying meta- and object-levels interactions, as observed in the study of provability in proof theory and metaprogramming in computer science.

Ramesh's framework rests on two central concepts: *introspective theories* and *geminal categories*. A prototypical example of a geminal category is the syntactic category of a formal theory where Löb's theorem or the Gödel–Löb axiom holds. As mentioned above, geminal categories possess the property of self-internalization, having their own internalizations within themselves (see Theorem 6.8). On the other hand, an *introspective theory* is a categorical formalization of theories that describe such self-internalizing structures (more precisely, it is a *finite limit theory*; see Section 2.5). These two notions are intimately related: the theory of geminal categories is an (initial) introspective theory, and conversely, any introspective theory naturally gives rise to a geminal category.

Despite its significance, the theory has not yet gained wide recognition. To facilitate its broader adoption and further development, we believe that it is crucial to streamline the complexity of the framework and reorganize its core concepts. In particular, the original proof of the main result, *Löb's theorem for introspective theories* [38, Theorem 4.19], is remarkably intricate. Since Ramesh's development is centered on introspective theories, the theory of geminal categories is built upon them. Consequently, even *Löb's theorem for geminal categories* [38, Observation 5.24], which is a direct generalization of the standard Löb's theorem, is derived only as a corollary of the main result.

## 1.2 Our contributions

In this thesis, we focus on geminal categories and develop a self-contained theory that is mathematically and conceptually reorganized. As a consequence, we establish a categorical counterpart of the *Gödel–Löb axiom* which applies to any geminal category. Our formulation also provides a streamlined (and slightly generalized) proof of *Löb's theorem for (pre-)geminal categories*, which in turn yields a simplified proof of Ramesh's main theorem. Notably, our development relies neither on the notion of introspective theories nor on informal internal reasoning, such as "performing this argument internally." This approach makes the theory of geminal categories more accessible to a broader audience.

Conceptually, we provide a new motivation for the concept of geminal categories and the derivation of Löb's theorem for them by placing the notion of Gödel coding at the center. For this purpose, we introduce the notion of *code structures* on fibrations, which serves as a natural abstraction of Gödel coding from the perspective of fibrational semantics. In fact, by generalizing the classical diagonal lemma, we prove the *fixed point theorem for fibrations with codes.* This theorem is shown to unify various classical results studied under the name of *intensional recursion* by Kavvos [18, 19, 21].

We then introduce *pre-geminal categories,* a slight weakening of Ramesh's geminal categories, as a structure that canonically induces code structures on codomain fibrations. This formulation offers a new perspective: most of the definition of geminal categories can be viewed as natural requirements for inducing Gödel coding on the "universe," i.e., the codomain fibration. Our proof of *Löb's theorem for pre-geminal categories* is based on a *twofold* application of our fixed point theorem, highlighting the core idea of Ramesh's proof — what we call the *bootstrapping argument.*

Our strategy for establishing the *Gödel–Löb axiom for geminal categories* is simple: we internalize the derivation of Löb's theorem. The notion of geminal categories is thus re-introduced as a natural extension of a pre-geminal one, specifically one that induces *internal pre-geminal categories.* Finally, we clarify the modal structure arising from geminal categories and compare it with existing models in modal calculi. This suggests a potential link to *metaprogramming,* which shares our focus on meta-object interaction.

### 1.3 Organization of the thesis

The remainder of this thesis is organized as follows. Section 2 provides preliminaries and fixes our notation. In Section 3, we introduce *code structures* on fibrations and establish the *fixed point theorem for fibrations with codes.* Section 4 explores categorical structures that induce these code structures, leading to the notion of *pre-geminal categories.* In Section 5, after an overview of the *bootstrapping argument,* we prove *Löb's theorem for pre-geminal categories* by applying the preceding results. Finally, Section 6 establishes the *Gödel–Löb axiom for geminal categories* through internalization and compares the resulting modal structure with existing models of modal calculi. Section 7 concludes the thesis and discusses future work.

## 2 Preliminaries

This section provides the necessary preliminaries and fixes the notation used throughout this thesis. While most of the content is standard, we introduce several conventions and notations specific to this work. The reader may refer back to this section as needed.

## 2.1 Categories

### 2.1.1 General

Let $\mathcal{B}$ be a category. The class of objects of $\mathcal{B}$ is denoted by $|\mathcal{B}|$. As a shorthand, we often write $A \in \mathcal{B}$ to mean $A \in |\mathcal{B}|$. For $A, B \in \mathcal{B}$, the class of morphisms from $A$ to $B$ is denoted by $\mathcal{B}(A, B)$.

The arrow category of $\mathcal{B}$ is denoted by $\mathcal{B}^{\rightarrow}$. It is the functor category $[\mathbf{2}, \mathcal{B}]$ from $\mathbf{2} = (\cdot \rightarrow \cdot)$ to $\mathcal{B}$. Note that $|\mathcal{B}^{\rightarrow}|$ corresponds to the class of morphisms of $\mathcal{B}$.

A terminal object in a category $\mathcal{B}$ is denoted as $1_{\mathcal{B}}$ or simply $1$. The unique morphism $A \rightarrow 1$ is denoted as $!_A$ or simply $!$. Regarding morphisms from a terminal object, we follow the following convention:

**Definition 2.1 (Global section functor)** For a locally small category $\mathcal{B}$ with a terminal object $1$, a hom-functor $\mathcal{B}(1, -) : \mathcal{B} \rightarrow \mathbf{Set}$ is called a *global section functor* of $\mathcal{B}$. It is denoted by $\Gamma_{\mathcal{B}}$ or simply $\Gamma$.

**Definition 2.2 (Global element)** Let $\mathcal{B}$ be a category with a terminal object $1$. For an object $A \in \mathcal{B}$, a morphism $a : 1 \rightarrow A$ is called a *global element* of $A$. We use the notation $a \in_1 A$ for such an element. In other words, $a \in_1 A$ means $a \in \Gamma A$.

For a global element $a \in_1 A$ and a morphism $f : A \rightarrow B$, we write $f(a) \in_1 B$ for the composite $f \circ a$.

We write $\pi_1 : A \times B \rightarrow A$ and $\pi_2 : A \times B \rightarrow B$ for the projections from a binary product. For $f : C \rightarrow A$ and $g : C \rightarrow B$, we write $\langle f, g \rangle : C \rightarrow A \times B$ for the unique morphism satisfying $\pi_1 \circ \langle f, g \rangle = f$ and $\pi_2 \circ \langle f, g \rangle = g$. For $a \in_1 A, b \in_1 B$ and $f : A \times B \rightarrow C$, we write $f(a, b) \in_1 C$ for the composite $f \circ \langle a, b \rangle$.

We also use $\pi_1$ and $\pi_2$ for the projections from a pullback. In diagrams, a pullback square is indicated by a corner symbol, as shown below:

$$\begin{array}{ccc} A \times_C B & \xrightarrow{\pi_2} & B \\ {\scriptstyle \pi_1}\downarrow \lrcorner & & \downarrow{\scriptstyle g} \\ A & \xrightarrow{f} & C. \end{array}$$

Additionally, for $h : D \rightarrow A$ and $k : D \rightarrow B$ such that $f \circ h = g \circ k$, we write $\langle h, k \rangle : D \rightarrow A \times_C B$ for the unique morphism satisfying $\pi_1 \circ \langle h, k \rangle = h$ and $\pi_2 \circ \langle h, k \rangle = k$.

In Cartesian closed categories, we write $B^A$ for exponentials and $\mathrm{ev} : B^A \times A \rightarrow B$ for the evaluation morphism. For any $f : A \times B \rightarrow C$, we write $\lambda(f) : A \rightarrow C^B$ for its *Currying*, the unique morphism satisfying $\mathrm{ev} \circ (\lambda(f) \times \mathrm{id}_B) = f$.

### 2.1.2 Property-like structures

Structures on a category are often unique up to isomorphism, but not in the strict sense. Such structures are frequently referred to as *property-like structures*. Typical examples include terminal objects, binary products, and pullbacks. While these remarks might seem pedantic, they become relevant in certain places, such as in Remark 3.3.

We say a *category with* $X$ if at least one property-like structure $X$ exists on the category, though it is not specifically chosen. In contrast, if $X$ is specified as a fixed structure, we call it *a category with chosen* $X$. Following standard categorical practice, even when $X$ is not explicitly specified, we may implicitly choose and work with it, provided the argument does not depend on the particular choice.

A *functor preserving* $X$ refers to a functor that preserves $X$ up to isomorphism. In contrast, for functors between categories with chosen $X$, we use the term a *functor strictly preserving* $X$ if it preserves $X$ in the strict sense.

Throughout this thesis, *finite products* always mean "terminal objects and binary products," and *finite limits* always mean "terminal objects and pullbacks." Thus, a *category with chosen finite limits* is a category equipped with a specific choice of a terminal object and pullbacks.

For internal categories (to be introduced in Section 2.4), we always treat property-like structures as being chosen. Accordingly, we shall only consider *internal categories with chosen finite limits* and will avoid the term "internal categories with finite limits" without a specific choice. The same convention applies to internal functors.

## 2.2 Fibrations

We recall basic facts about fibrations. We refer to Jacobs [14] for the details.

### 2.2.1 Basic definitions

Consider a functor $p : \mathcal{E} \to \mathcal{B}$. An object $X \in \mathcal{E}$ is said to be *above* an object $A \in \mathcal{B}$ if $pX = A$. Similarly, a morphism $s \in |\mathcal{E}^{\to}|$ is *above* a morphism $f \in |\mathcal{B}^{\to}|$ if $ps = f$.

**Definition 2.3 (Cartesian morphism)** Let $p : \mathcal{E} \to \mathcal{B}$ be a functor. A morphism $s : X \to Y$ in $\mathcal{E}$ above $f : A \to B$ is said to be *Cartesian* over $f$ if, for any morphism $g : C \to A$ in $\mathcal{B}$ and any $u : Z \to Y$ above $f \circ g$, there exists a unique $t : Z \to X$ above $g$ such that $u = s \circ t$.

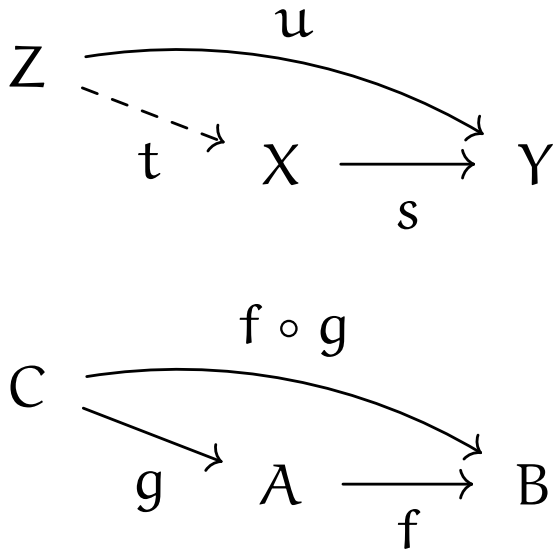

**Definition 2.4 (Fibration)** A functor $p : \mathcal{E} \to \mathcal{B}$ is called a *fibration* (or a *Grothendieck fibration*) if, for any morphism $f : A \to B$ in $\mathcal{B}$ and any object $X \in \mathcal{E}$ above $B$, there exists a Cartesian morphism $\overline{f}(X) : f^*(X) \to X$ over $f$. This morphism $\overline{f}(X)$ is called a *Cartesian lifting* of $f$ to $X$.

We refer to $\mathcal{B}$ as the *base category* and $\mathcal{E}$ as the *total category* of the fibration $p$.

$$\begin{array}{ccc} f^*(X) & \overset{\exists\, \overline{f}(X)}{\dashrightarrow} & X \\ & & \\ A & \xrightarrow{f} & B \end{array}$$

For the base category of a fibration, we use letters $A, B, C$ for its objects and $f, g, h$ for its morphisms. For the total category, we use $X, Y, Z$ for objects and $s, t, u$ for morphisms.

Cartesian liftings are unique only up to isomorphism, which leads to the following definition.

**Definition 2.5 (Cleavage)** A *cleavage* on a fibration $p$ is a choice of a Cartesian lifting $\overline{f}(X) : f^*(X) \to X$ to each pair $(f, X)$, where $f : A \to B$ in $\mathcal{B}$ and $X \in \mathcal{E}$ is above $B$. A fibration equipped with a cleavage is called a *cloven fibration.*

If a cleavage is regarded as a property-like structure on a functor, then a fibration is simply a functor with a cleavage, while a cloven fibration is a functor with a *chosen* cleavage.

**Definition 2.6 (Fiber category)** For a fibration $p : \mathcal{E} \to \mathcal{B}$ and $A \in \mathcal{B}$, the *fiber category* (or simply the *fiber*) at $A$ is the subcategory of $\mathcal{E}$ consisting of objects above $A$ and morphisms above $\mathrm{id}_A$. The fiber at $A$ is denoted by $\mathcal{E}_A$.

**Definition 2.7 (Reindexing functor)** For a cloven fibration $p : \mathcal{E} \to \mathcal{B}$ and a morphism $f : A \to B$ in $\mathcal{B}$, the assignment $X \mapsto f^*(X)$ from $|\mathcal{E}_B|$ to $|\mathcal{E}_A|$ extends to a functor $f^* : \mathcal{E}_B \to \mathcal{E}_A$. This functor is called the *reindexing functor* along $f$. We often write $X \cdot f$ as a shorthand for $f^*(X)$.

For a non-cloven fibration, reindexing functors are determined only up to natural isomorphism. We often write $f^*$ even in the non-cloven case, implicitly assuming a choice has been made.

**Proposition 2.8** *For a fibration $p : \mathcal{E} \to \mathcal{B}$, the following hold:*
(1) *For each $A \in \mathcal{B}$, there is a natural isomorphism $\mathrm{id}_{\mathcal{E}_A} \cong (\mathrm{id}_A)^*$.*
(2) *For each $f : A \to B$, $g : B \to C$ in $\mathcal{B}$, there is a natural isomorphism $f^* \circ g^* \cong (g \circ f)^*$.*

One of the most fundamental results is the equivalence between fibrations and indexed categories (see, e.g., Jacobs [14, Section 1.10]).

**Theorem 2.9** *For a category $\mathcal{B}$, there is a one-to-one correspondence, up to isomorphism, between:*

- *Cloven fibrations $p : \mathcal{E} \to \mathcal{B}$ over $\mathcal{B}$,*
- *Indexed categories over $\mathcal{B}$, i.e., pseudofunctors $\Phi : \mathcal{B}^{\mathrm{op}} \to$* **Cat.**

*Explicitly, for a cloven fibration* $p : \mathcal{E} \to \mathcal{B}$*, the corresponding indexed category* $\Phi : \mathcal{B}^{op} \to$ **Cat** *is given by*

$$\Phi(A) = \mathcal{E}_A, \qquad \Phi(f) = f^*.$$

*This correspondence extends to an equivalence between the category of cloven fibrations and the category of indexed categories.*

An advantage of fibrations over indexed categories lies in the concept of non-cloven fibrations (see Bénabou [5], Jacobs [14]). However, in this thesis, we primarily focus on cloven fibrations, as explicit reindexing functors are required to define code structures (Definition 3.1). Consequently, the arguments throughout this thesis could be developed using indexed categories with virtually no difference.

### 2.2.2 Weak equality between morphisms in the base category

We introduce a weak notion of equality between morphisms in the base category. This notion and its notation are specific to this paper and are not standard.

**Definition 2.10 (Weak equality)** Let $p : \mathcal{E} \to \mathcal{B}$ be a fibration on a category $\mathcal{B}$ with finite products. For any parallel morphisms $f, f' : A \to B$ in $\mathcal{B}$, we write $f \approx f'$ if the reindexing functors along $\langle id, f\rangle$ and $\langle id, f'\rangle : A \to A \times B$ are naturally isomorphic.

Obviously, the relation $f \approx f'$ is an equivalence relation. Specifically, the usual (or *strict, external*) equality $f = f'$ implies $f \approx f'$.

**Proposition 2.11** *Let* $f, f' : A \to B$ *be parallel morphisms in* $\mathcal{B}$ *such that* $f \approx f'$.

(1) *The reindexing functors* $f^*$ *and* $(f')^*$ *are naturally isomorphic.*

(2) *For any morphism* $g : A \to C$*, the reindexing functors along* $\langle g, f\rangle$ *and* $\langle g, f'\rangle : A \to C \times B$ *are naturally isomorphic.*

(3) *For any morphism* $h : B \to C$*, it holds that* $h \circ f \approx h \circ f'$.

*Proof.* (1) We have $f = \pi_2 \circ \langle id, f\rangle$ and $f' = \pi_2 \circ \langle id, f'\rangle$. Since $\langle id, f\rangle^* \cong \langle id, f'\rangle^*$ by assumption, it follows that $f^* \cong \langle id, f\rangle^* \circ \pi_2^* \cong \langle id, f'\rangle^* \circ \pi_2^* \cong (f')^*$.

(2) Follows from the same argument applied to $\langle g, f\rangle = (g \times id) \circ \langle id, f\rangle$ and $\langle g, f'\rangle = (g \times id) \circ \langle id, f'\rangle$.

(3) Follows from $\langle id, h \circ f\rangle = (id \times h) \circ \langle id, f\rangle$ and $\langle id, h \circ f'\rangle = (id \times h) \circ \langle id, f'\rangle$. □

In a fibration with certain properties, known as *equational fibrations*, one can define a notion of *internal equality* between morphisms [14, Section 3.4]. If $f$ and $f'$ are internally equal in this sense, then $f \approx f'$ holds. In fact, our motivation for introducing this weak equality is to subsume both internal equality in equational fibrations and strict (external) equality in non-equational cases. This notion may be replaced with strict equality throughout this thesis, except for examples constructed in Section 3.2.

### 2.2.3 Basic constructions

We recall several basic constructions of fibrations: the representable, codomain, subobject, and family fibrations. Among these, the first two play an important role in this thesis.

> **Definition 2.12 (Representable fibration)** For a category $\mathcal{B}$ and an object $A \in \mathcal{B}$, the functor $\mathrm{dom}_A : \mathcal{B}/A \to \mathcal{B}$ from the slice category $\mathcal{B}/A$ is defined as follows:
>
> **On objects** For an object $x : X \to A$ in $\mathcal{B}/A$, we define $\mathrm{dom}_A(x) = \mathrm{dom}(x) = X$.
>
> **On morphisms** For a morphism $s : x \to y$ in $\mathcal{B}/A$, we define $\mathrm{dom}_A(s) = s$.
>
> This functor forms a fibration, called a *representable fibration* over $\mathcal{B}$.

This terminology reflects the fact that $\mathrm{dom}_A$ corresponds to the representable functor $\mathcal{B}(-, A)$ via the equivalence in Theorem 2.9. The fiber at $B \in \mathcal{B}$ is the hom-set $\mathcal{B}(B, A)$, regarded as a discrete category. The cleavage is unique, and reindexing along $f : B \to C$ is given by precomposition: $f^* = (-) \circ f : \mathcal{B}(C, A) \to \mathcal{B}(B, A)$.

> **Definition 2.13 (Codomain fibration)** Let $\mathcal{B}$ be a category with pullbacks. The functor $\mathrm{cod}_{\mathcal{B}} : \mathcal{B}^{\to} \to \mathcal{B}$ is defined as follows:
>
> **On objects** For an object $x : X \to A$ in $\mathcal{B}^{\to}$, we define $\mathrm{cod}_{\mathcal{B}}(x) = \mathrm{cod}(x) = A$.
>
> **On morphisms** For a morphism $(f, s) : x \to y$ in $\mathcal{B}^{\to}$, we define $\mathrm{cod}_{\mathcal{B}}(f, s) = f$, where $f : A \to B$ and $s : X \to Y$ are such that $f \circ x = y \circ s$.
>
> $$\begin{array}{ccc} X & \xrightarrow{s} & Y \\ {\scriptstyle x}\downarrow & & \downarrow{\scriptstyle y} \\ A & \xrightarrow{f} & B \end{array}$$
>
> The functor $\mathrm{cod}_{\mathcal{B}}$ forms a fibration, called the *codomain fibration* over $\mathcal{B}$.

The fiber at $A \in \mathcal{B}$ is the slice category $\mathcal{B}/A$. A cleavage on this fibration is equivalent to choices of pullbacks in $\mathcal{B}$. If $\mathcal{B}$ has chosen pullbacks, we regard $\mathrm{cod}_{\mathcal{B}}$ as a cloven fibration. The reindexing $f^*$ is given by the pullback along $f$, as illustrated below:

$$\begin{array}{ccc} X & \longrightarrow & Y \\ {\scriptstyle y \cdot f}\downarrow & \lrcorner & \downarrow{\scriptstyle y} \\ A & \xrightarrow{f} & B. \end{array}$$

The codomain fibration is restricted to the subobject fibration:

> **Definition 2.14 (Subobject fibration)** Let $\mathcal{B}$ be a category with pullbacks. The category $\mathrm{Sub}(\mathcal{B})$ is defined as follows:
>
> **Objects** Pairs $(A, [m])$ where $A \in \mathcal{B}$ and $[m]$ is a *subobject* of $A$, i.e., an equivalence class of monomorphisms into $A$.
>
> **Morphisms** $(A, [m]) \to (B, [n])$ are morphisms $f : A \to B$ such that there exists a morphism $s$ such that the following commutes:

$$\begin{array}{ccc} X & \overset{s}{\dashrightarrow} & Y \\ {\scriptstyle m}\downarrow & & \downarrow{\scriptstyle n} \\ A & \xrightarrow{f} & B. \end{array}$$

The functor $p : \mathrm{Sub}(\mathcal{B}) \to \mathcal{B}$, defined by $p(A, [m]) = A$ and $p(f) = f$, forms a fibration called the *subobject fibration* over $\mathcal{B}$.

Finally, any category $\mathcal{C}$ gives rise to a fibration over **Set**:

**Definition 2.15 (Family fibration)** Let $\mathcal{C}$ be a category. The category $\mathrm{Fam}(\mathcal{C})$ is defined as follows:

**Objects** Pairs $\big(A, (X_a)_{a \in A}\big)$ where $A$ is a set and $(X_a)_{a \in A}$ is a family of objects of $\mathcal{C}$ indexed by $A$.

**Morphisms** $(A, (X_a)) \to (B, (Y_b))$ are pairs $\big(f, (s_a)_{a \in A}\big)$ where $f : A \to B$ is a function and $\big(s_a : X_a \to Y_{f(a)}\big)_{a \in A}$ is a family of morphisms in $\mathcal{C}$ indexed by $A$.

The functor $p : \mathrm{Fam}(\mathcal{C}) \to \mathbf{Set}$, defined by $p(A, (X_a)) = A, p(f, (s_a)) = f$, forms a fibration called the *family fibration* of $\mathcal{C}$.

Note that a cleavage on this fibration is canonically determined.

### 2.2.4 Syntactic fibrations

We introduce another class of examples called *syntactic fibrations* (referred to as *fibrations of contexts* in Jacobs [14]). This provides a powerful framework for organizing various logical systems into categorical structures. Our description follows Jacobs [14, Section 3.1]; we refer the reader there for details.

Let $\Sigma$ be a *many-typed signature*. It consists of a set $|\Sigma|$ of *basic types* and a set of *function symbols*. Each function symbol $f$ is assigned a finite list of *input types* $\sigma_1, ..., \sigma_n$ and exactly one *output type* $\sigma_{n+1}$, denoted by $f : \sigma_1, ..., \sigma_n \longrightarrow \sigma_{n+1}$. For instance, a standard signature for first-order arithmetic is defined by a single basic type $N$ and the following function symbols:

$$0 : () \longrightarrow N, \qquad S : N \longrightarrow N, \qquad + : N, N \longrightarrow N, \qquad \times : N, N \longrightarrow N.$$

Fix an infinite sequence of variables $v_1, v_2, ....$ A *type context* $\Gamma$ is an assignment of basic types to variables $v_1, ..., v_n$ for some non-negative integer $n$, denoted by $\Gamma = (v_1 : A_1, ..., v_n : A_n)$.

The notion of *terms* in a type context is defined in a standard inductive manner. Each term is assigned a unique type, and all variables appearing in a term must be present in its context. For instance, the typing judgment

$$x : N, y : N, z : N \vdash (x \times y) + (x \times z) : N$$

denotes that the expression $(x \times y) + (x \times z)$, a shorthand for $+(\times(x, y), \times(x, z))$, is a term of type $N$ in the context $(x : N, y : N, z : N)$.

*Simultaneous substitution* for terms is defined in the standard way. For terms $M, N_1, ..., N_n$ and variables $v_1, ..., v_n$, we write $M[N_1, ..., N_n/v_1, ..., v_n]$ for the term obtained by the simultaneous substitution of $N_i$ for $v_i$ in $M$.

Terms over a given signature $\Sigma$ can be organized into the *classifying category* $\mathrm{Cl}(\Sigma)$. This category serves as the base category of syntactic fibrations.

**Definition 2.16 (Classifying category)** For a many-typed signature $\Sigma$, its *classifying category* $\mathrm{Cl}(\Sigma)$ is defined as follows:

**Objects** Type contexts $\Gamma$ over $\Sigma$.

**Morphisms** For contexts $\Gamma$ and $\Delta = (v_1 : \sigma_1, ..., v_n : \sigma_n)$, morphisms $\Gamma \to \Delta$ are finite lists of terms $(M_1, ..., M_n)$ such that $\Gamma \vdash M_i : \sigma_i$ for each $i$.

**Identities** For $\Gamma = (v_1 : \sigma_1, ..., v_n : \sigma_n)$, the identity $\mathrm{id}_\Gamma$ is the list of variables $(v_1, ..., v_n)$.

**Composition** For $f = (M_1, ..., M_m) : \Gamma \to \Delta$ and $g = (N_1, ..., N_n) : \Delta \to \Theta$, their composition $g \circ f$ is defined by

$$(N_1[M_1, ..., M_m/v_1, ..., v_m], ..., N_n[M_1, ..., M_m/v_1, ..., v_m]) : \Gamma \to \Theta.$$

The concatenation of contexts provides finite products in $\mathrm{Cl}(\Sigma)$. Specifically, the empty context $()$ is a terminal object.

Next, we consider logical systems built over a given signature. In addition to the signature, a logical system may involve additional symbols (such as predicate symbols) and axioms. Such data, together with its underlying signature, is called a *specification* and denoted by $(\Sigma, \Pi)$. Specifications can represent various logical systems, such as equational logic, first-order logic, and higher-order logic. We refer the reader to Jacobs [14] for their formal definitions.

Given a specification, *propositions* are constructed within type contexts. As with terms, any free variables occurring in a proposition must be present in its context. For example, in a first-order specification extending the signature for arithmetic described above with a binary predicate symbol $<$, we have the following proposition:

$$x : N \vdash (x > 0) \to (\exists y : N.\, y + 1 =_N x) \quad \mathsf{Prop}.$$

Notice that both quantification and equality are typed.

A finite list of propositions within a common type context $\Gamma$ is called a *proposition context* in $\Gamma$. Given a specification, for a proposition context $\Phi$ and a proposition $\varphi$ in $\Gamma$, one can define the derivation of judgments of the form:

$$\Gamma \mid \Phi \vdash \varphi.$$

Intuitively, this judgment expresses that the proposition $\varphi$ follows from the assumptions in $\Phi$ for any assignment of the variables in $\Gamma$. For instance, in our arithmetic example, we may have a derivable judgment:

$$x : N, y : N, z : N \mid x + z =_N y + z \vdash x =_N y.$$

Now we can describe the *syntactic fibration* for a given specification.

**Definition 2.17 (Syntactic fibration)** Given a specification $(\Sigma, \Pi)$, its *syntactic fibration* $p : \mathcal{L}(\Sigma, \Pi) \to \mathrm{Cl}(\Sigma)$ is defined as follows. The base category is the classifying category $\mathrm{Cl}(\Sigma)$. The total category $\mathcal{L}(\Sigma, \Pi)$ is defined by:

**Objects** Pairs $(\Gamma \mid \Phi)$ of a type context $\Gamma$ and a proposition context $\Phi$ in $\Gamma$.

**Morphisms** $(\Gamma \mid \Phi) \to (\Delta \mid \Psi)$ are morphisms $(M_1, ..., M_n) : \Gamma \to \Delta$ in $\mathrm{Cl}(\Sigma)$ such that, for every proposition $\psi$ in $\Psi$, the following is derivable:

$$\Gamma \mid \Phi \vdash \psi[M_1, ..., M_n/v_1, ..., v_n].$$

**Identities and composition** Inherited from $\mathrm{Cl}(\Sigma)$.

Finally, the functor $p : \mathcal{L}(\Sigma, \Pi) \to \mathrm{Cl}(\Sigma)$ is the obvious projection $(\Gamma \mid \Phi) \mapsto \Gamma$.

The functor $p : \mathcal{L}(\Sigma, \Pi) \to \mathrm{Cl}(\Sigma)$ forms a cloven fibration (in fact, a split one), with the reindexing operation given by simultaneous substitution:

$$(\Delta \mid \Psi) \cdot (M_1, ..., M_n) = (\Gamma \mid \Psi[M_1, ..., M_n/v_1, ..., v_n]).$$

We further assume that the logical system possesses internal (propositional) equality $=_\sigma$ for each type $\sigma$. This equality is required to satisfy the following rule:

$$\frac{\Gamma \mid \Phi \vdash \varphi[M/v] \qquad \Gamma \mid () \vdash M =_\sigma N}{\Gamma \mid \Phi \vdash \varphi[N/v].}$$

From this rule, the following property follows:

**Proposition 2.18** *Consider a syntactic fibration* $p : \mathcal{L}(\Sigma, \Pi) \to \mathrm{Cl}(\Sigma)$. *Let* $f = (M_1, ..., M_n)$ *and* $g = (N_1, ..., N_n) : \Gamma \to \Delta$ *be parallel arrows in* $\mathrm{Cl}(\Sigma)$. *If each judgment*

$$\Gamma \mid () \vdash M_i =_{\sigma_i} N_i \quad (i = 1, ..., n)$$

*is derivable, then* $f \approx g$ *holds in this fibration.*

## 2.3 Enriched categories

We review basic definitions about enriched categories, which will be primarily used in Section 4.2. Throughout this thesis, we consider enrichment only over Cartesian monoidal structures, even if not explicitly stated.

### 2.3.1 Basic definitions

**Definition 2.19 (Enriched category)** Let $\mathcal{B}$ be a category with finite products. A $\mathcal{B}$-*enriched category* $\mathcal{C}$ consists of the following data:

- A set $|\mathcal{C}|$ of *objects*.
- For each pair of objects $X, Y \in |\mathcal{C}|$, a *hom-object* $\mathcal{C}(X, Y) \in \mathcal{B}$.
- For each object $X \in |\mathcal{C}|$, an *identity element* $e_X \in_1 \mathcal{C}(X, X)$, i.e., a morphism $e_X : 1 \to \mathcal{C}(X, X)$.
- For each triple $X, Y, Z \in |\mathcal{C}|$, a *composition morphism* $\circ_{X,Y,Z} : \mathcal{C}(X, Y) \times \mathcal{C}(Y, Z) \to \mathcal{C}(X, Z)$ in $\mathcal{B}$.

These data must satisfy the following commutative diagrams:

$$
\begin{array}{ccc}
1 \times \mathcal{C}(X,Y) & \xrightarrow{e_X \times \mathrm{id}} & \mathcal{C}(X,X) \times \mathcal{C}(X,Y) \\
& \searrow^{\cong} & \downarrow{\circ_{X,X,Y}} \\
& & \mathcal{C}(X,Y),
\end{array}
\qquad
\begin{array}{ccc}
\mathcal{C}(X,Y) \times 1 & \xrightarrow{\mathrm{id} \times e_Y} & \mathcal{C}(X,Y) \times \mathcal{C}(Y,Y) \\
& \searrow^{\cong} & \downarrow{\circ_{X,Y,Y}} \\
& & \mathcal{C}(X,Y),
\end{array}
$$

$$
\begin{array}{ccc}
(\mathcal{C}(X,Y) \times \mathcal{C}(Y,Z)) \times \mathcal{C}(Z,W) & \xrightarrow{\cong} & \mathcal{C}(X,Y) \times (\mathcal{C}(Y,Z) \times \mathcal{C}(Z,W)) \\
\downarrow{\circ_{X,Y,Z} \times \mathrm{id}} & & \downarrow{\mathrm{id} \times \circ_{Y,Z,W}} \\
\mathcal{C}(X,Z) \times \mathcal{C}(Z,W) & & \mathcal{C}(X,Y) \times \mathcal{C}(Y,W) \\
& {\scriptstyle \circ_{X,Z,W}} \searrow \quad \swarrow {\scriptstyle \circ_{X,Y,W}} & \\
& \mathcal{C}(X,W). &
\end{array}
$$

Notice that we use the diagrammatic order for composition morphisms. As for ordinary categories, we often abbreviate $X \in |\mathcal{C}|$ as $X \in \mathcal{C}$. The morphisms $e_X$ and $\circ_{X,Y,Z}$ are often written simply as $e$ and $\circ$, or as $e_{\mathcal{C}}$ and $\circ_{\mathcal{C}}$ to clarify the category. We also note that a **Set**-enriched category is simply a locally small category.

**Definition 2.20 (Enriched functor)** Let $\mathcal{C}, \mathcal{D}$ be $\mathcal{B}$-enriched categories. A *$\mathcal{B}$-enriched functor* $H : \mathcal{C} \to \mathcal{D}$ consists of the following data:

- A function $|H| : |\mathcal{C}| \to |\mathcal{D}|$. We typically write $HX$ for $|H|X$.
- For each $X, Y \in \mathcal{C}$, a morphism $H_{X,Y} : \mathcal{C}(X,Y) \to \mathcal{D}(HX, HY)$ in $\mathcal{B}$.

These data must satisfy the following commutative diagrams:

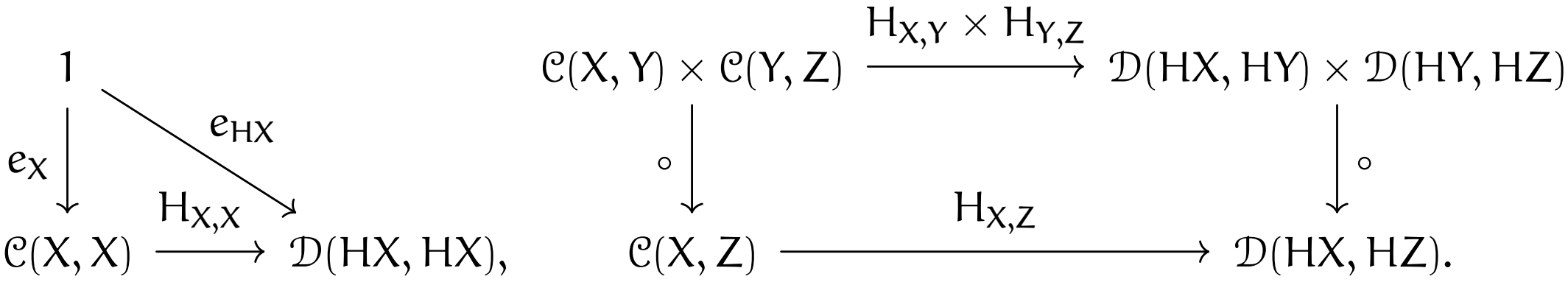

### 2.3.2 Change of base and underlying categories

Given a $\mathcal{B}$-enriched category $\mathcal{C}$ and a functor $F : \mathcal{B} \to \mathcal{B}'$ preserving finite products, one can obtain a $\mathcal{B}'$-enriched category by "pushing forward" the structure of $\mathcal{C}$ along $F$. This construction is called the *change of base*.

**Definition 2.21 (Change of base)** Let $\mathcal{C}$ be a $\mathcal{B}$-enriched category and $F : \mathcal{B} \to \mathcal{B}'$ be a functor preserving finite products. The $\mathcal{B}'$-enriched category $F_*\mathcal{C}$, called the *change of base* of $\mathcal{C}$ along $F$, is defined as follows:

**Objects** $|F_*\mathcal{C}| = |\mathcal{C}|$.

**Hom-objects** For $X, Y \in F_*\mathcal{C}$, the hom-object is $(F_*\mathcal{C})(X,Y) = F(\mathcal{C}(X,Y))$.

Identities and composition of $F_*\mathcal{C}$ are induced by those of $\mathcal{C}$ via $F$:

**Identity elements** $e_{F_*\mathcal{C}}$ is given by

$$1_{\mathcal{B}'} \xrightarrow{\cong} F(1_{\mathcal{B}}) \xrightarrow{F(e_{\mathcal{C}})} F(\mathcal{C}(X, X)).$$

**Composition morphisms** $\circ_{F_*\mathcal{C}}$ is given by the following morphism:

$$F(\mathcal{C}(X, Y)) \times F(\mathcal{C}(Y, Z)) \xrightarrow{\cong} F(\mathcal{C}(X, Y) \times \mathcal{C}(Y, Z)) \xrightarrow{F(\circ_{\mathcal{C}})} F(\mathcal{C}(X, Z)).$$

**Definition 2.22 (Change of base for enriched functors)** Let $H : \mathcal{C} \to \mathcal{D}$ be a $\mathcal{B}$-enriched functor between $\mathcal{B}$-enriched categories. For a functor $F : \mathcal{B} \to \mathcal{B}'$ preserving finite products, the $\mathcal{B}'$-enriched functor $F_*H : F_*\mathcal{C} \to F_*\mathcal{D}$ is defined as follows:

**On objects** $|F_*H| = |H|$.

**On hom-objects** For $X, Y \in \mathcal{C}$, the morphism $(F_*H)_{X,Y}$ is defined to be $F(H_{X,Y}) : F(\mathcal{C}(X, Y)) \to F(\mathcal{D}(X, Y))$.

If $\mathcal{B}$ is locally small, there exists a global section functor $\Gamma : \mathcal{B} \to \mathbf{Set}$, which preserves finite products. Applying the change of base construction along this functor, we can obtain an ordinary category from a $\mathcal{B}$-enriched category $\mathcal{C}$.

**Definition 2.23 (Underlying category)** Let $\mathcal{C}$ be a $\mathcal{B}$-enriched category and assume that $\mathcal{B}$ is locally small. The *underlying category* of $\mathcal{C}$, denoted by $\Gamma_*\mathcal{C}$, is the change of base of $\mathcal{C}$ along the global section functor $\Gamma : \mathcal{B} \to \mathbf{Set}$. Explicitly, the category $\Gamma_*\mathcal{C}$ consists of:

**Objects** $X \in |\mathcal{C}|$.

**Morphisms** $X \to Y$ are global elements $f \in_1 \mathcal{C}(X, Y)$.

Whenever we consider the underlying category of a $\mathcal{B}$-enriched category, we implicitly assume $\mathcal{B}$ to be locally small.

**Definition 2.24 (Underlying functor)** Let $H : \mathcal{C} \to \mathcal{D}$ be a $\mathcal{B}$-enriched functor and assume that $\mathcal{B}$ is locally small. The *underlying functor* of $H$, denoted by $\Gamma_*H$, is the change of base of $H$ along $\Gamma : \mathcal{B} \to \mathbf{Set}$. Explicitly, the functor $\Gamma_*H : \Gamma_*\mathcal{C} \to \Gamma_*\mathcal{D}$ is described as follows:

**On objects** For $X \in \mathcal{C}$, $(\Gamma_*H)X = HX$.

**On morphisms** For $f \in_1 \mathcal{C}(X, Y)$, we define $(\Gamma_*H)(f) = H_{X,Y}(f)$, i.e., the composite

$$1 \xrightarrow{f} \mathcal{C}(X, Y) \xrightarrow{H_{X,Y}} \mathcal{D}(HX, HY).$$

Although the underlying category $\Gamma_*\mathcal{C}$ loses some information from the original enriched category $\mathcal{C}$, it is useful for establishing a connection between enriched categories and ordinary categories. For instance, an enriched analog of a hom-functor can be formulated as an ordinary functor from the underlying category $\Gamma_*\mathcal{C}$ to $\mathcal{B}$.

**Definition 2.25 (Hom-functor)** Let $\mathcal{C}$ be a $\mathcal{B}$-enriched category and $X \in \mathcal{C}$. The *hom-functor* $\mathcal{C}(X, -)$ is the ordinary functor $\Gamma_*\mathcal{C} \to \mathcal{B}$ defined as follows:

**On objects** For each object $Y \in \mathcal{C}$, the functor maps $Y$ to the hom-object $\mathcal{C}(X, Y) \in \mathcal{B}$.

**On morphisms** For a morphism $f \in_1 \mathcal{C}(Y, Z)$ in $\Gamma_*\mathcal{C}$, the morphism $\mathcal{C}(X, f) : \mathcal{C}(X, Y) \to \mathcal{C}(X, Z)$ is defined by the following composite in $\mathcal{B}$:

$$\mathcal{C}(X, Y) \xrightarrow{\langle \mathrm{id}, f \circ ! \rangle} \mathcal{C}(X, Y) \times \mathcal{C}(Y, Z) \xrightarrow{\circ_{\mathcal{C}}} \mathcal{C}(X, Y).$$

In fact, if $\mathcal{B}$ is Cartesian closed, one can define hom-functors as $\mathcal{B}$-enriched functors from $\mathcal{C}$ to $\mathcal{B}$ (with its self-enrichment), which provides a more suitable notion for enriched categories. Yet, for our purposes, the definition as an ordinary functor from $\Gamma_*\mathcal{C}$ given above suffices. For a detailed account, see, e.g., Kelly [22].

The following definition and its internal analog (Definition 2.44) are used throughout this thesis. The superscript "e" in $F^{\#e}$ stands for enriched categories.

**Definition 2.26 ($F^{\#e}$)** Let $F : \mathcal{B} \to \mathcal{B}'$ be a functor preserving finite products. For a $\mathcal{B}$-enriched category $\mathcal{C}$, we define a functor $F^{\#e}_{\mathcal{C}} : \Gamma_*\mathcal{C} \to \Gamma_*F_*\mathcal{C}$ (or simply $F^{\#e}$) as follows:

**On objects** $F^{\#e}(X) = X$ for $X \in \mathcal{C}$. Note that $|\Gamma_*\mathcal{C}| = |\Gamma_*F_*\mathcal{C}| = |\mathcal{C}|$.

**On morphisms** For $f \in_1 \mathcal{C}(X, Y)$, the morphism $F^{\#e}(f)$ is defined as the global element of $(F_*\mathcal{C})(X, Y) = F(\mathcal{C}(X, Y))$ given by the following composite in $\mathcal{B}'$:

$$1_{\mathcal{B}'} \xrightarrow{\cong} F(1_{\mathcal{B}}) \xrightarrow{F(f)} F(\mathcal{C}(X, Y)).$$

### 2.3.3 Finite products

We define terminal objects and binary products in an enriched category as follows.

**Definition 2.27 (Finite products)** Let $\mathcal{C}$ be a $\mathcal{B}$-enriched category.

(1) A *terminal object* in $\mathcal{C}$ is an object $1_{\mathcal{C}} \in \mathcal{C}$ such that the morphism $! : \mathcal{C}(X, 1_{\mathcal{C}}) \to 1_{\mathcal{B}}$ is an isomorphism for any $X \in \mathcal{C}$.

(2) For objects $X, Y \in \mathcal{C}$, their *binary product* consists of an object $X \times Y \in \mathcal{C}$ and global elements $p_1 \in_1 \mathcal{C}(X \times Y, X)$ and $p_2 \in_1 \mathcal{C}(X \times Y, Y)$ such that the morphism

$$\langle \mathcal{C}(Z, p_1), \mathcal{C}(Z, p_2) \rangle : \mathcal{C}(Z, X \times Y) \to \mathcal{C}(Z, X) \times \mathcal{C}(Z, Y)$$

is an isomorphism for any $Z \in \mathcal{C}$.

If $\mathcal{C}$ has a terminal object and binary products, we say that $\mathcal{C}$ has *finite products*. For the general notion of limits in enriched categories, we refer the reader to Kelly [22]. The definition provided here is sufficient for our purpose.

As in ordinary categories, any hom-functor of enriched categories preserves finite products:

**Proposition 2.28** *Let $\mathcal{C}$ be a $\mathcal{B}$-enriched category with a terminal object (resp. binary products).*

(1) *The underlying category $\Gamma_*\mathcal{C}$ has a terminal object (resp. binary products) inherited from $\mathcal{C}$.*

(2) *For any $X \in \mathcal{C}$, the hom-functor $\mathcal{C}(X, -) : \Gamma_*\mathcal{C} \to \mathcal{B}$ preserves the terminal object (resp. binary products).*

*Proof.* These are immediate from the previous definition. □

The following proposition is an enriched analog of the identity $\langle f, g\rangle \circ h = \langle f \circ h, g \circ h\rangle$ in ordinary categories.

**Proposition 2.29** *Let $\mathcal{C}$ be a $\mathcal{B}$-enriched category with binary products. For any $W, X, Y, Z \in \mathcal{C}$, the following diagram commutes:*

$$\begin{array}{ccc}
\mathcal{C}(W,X) \times \mathcal{C}(X, Y \times Z) & \xrightarrow{\cong} & \mathcal{C}(W,X) \times (\mathcal{C}(X,Y) \times \mathcal{C}(X,Z)) \\
\Big\downarrow \circ_{\mathcal{C}} & & \Big\downarrow \langle\langle \pi_1, \pi_1 \circ \pi_2\rangle, \langle \pi_1, \pi_2 \circ \pi_2\rangle\rangle \\
 & & (\mathcal{C}(W,X) \times \mathcal{C}(X,Y)) \times (\mathcal{C}(W,X) \times \mathcal{C}(X,Z)) \\
 & & \Big\downarrow \circ_{\mathcal{C}} \times \circ_{\mathcal{C}} \\
\mathcal{C}(W, Y \times Z) & \xrightarrow{\cong} & \mathcal{C}(W,Y) \times \mathcal{C}(W,Z).
\end{array}$$

*Proof.* The first components of the two paths coincide since the following diagram commutes:

$$\begin{array}{ccc}
\mathcal{C}(W,X) \times \mathcal{C}(X, Y \times Z) & \xrightarrow{\circ_{W,X,Y\times Z}} & \mathcal{C}(W, Y \times Z) \\
\Big\downarrow \langle \mathrm{id}, p_1 \circ !\rangle & & \Big\downarrow \langle \mathrm{id}, p_1 \circ !\rangle \\
(\mathcal{C}(W,X) \times \mathcal{C}(X, Y \times Z)) \times \mathcal{C}(Y \times Z, Y) & \xrightarrow{\circ_{W,X,Y\times Z} \times \mathrm{id}} & \mathcal{C}(W, Y \times Z) \times \mathcal{C}(Y \times Z, Y) \\
\Big\downarrow \cong & & \\
\mathcal{C}(W,X) \times (\mathcal{C}(X, Y \times Z) \times \mathcal{C}(Y \times Z, Y)) & & \Big\downarrow \circ_{W,Y\times Z,Y} \\
\Big\downarrow \mathrm{id} \times \circ_{X,Y\times Z,Y} & & \\
\mathcal{C}(W,X) \times \mathcal{C}(X,Y) & \xrightarrow{\circ_{W,X,Y}} & \mathcal{C}(W,Y).
\end{array}$$

Similarly for the second components. □

The notions of hom-functors and terminal objects allow us to define global section functors for enriched categories. We write them as $\gamma$ to avoid the confusion with the global section functor $\Gamma$ for ordinary categories.

**Definition 2.30 (Global section functor)** Let $\mathcal{C}$ be a $\mathcal{B}$-enriched category with a terminal object $1_{\mathcal{C}}$. The hom-functor $\mathcal{C}(1_{\mathcal{C}}, -) : \Gamma_*\mathcal{C} \to \mathcal{B}$ is called the *global section functor* of $\mathcal{C}$ and is denoted by $\gamma_{\mathcal{C}}$ (or simply $\gamma$).

**Proposition 2.31 ($\Gamma \circ \gamma = \Gamma$)** *For a $\mathcal{B}$-enriched category $\mathcal{C}$ with a terminal object, the identity $\Gamma_{\mathcal{B}} \circ \gamma_{\mathcal{C}} = \Gamma_{\Gamma_*\mathcal{C}}$ holds.*

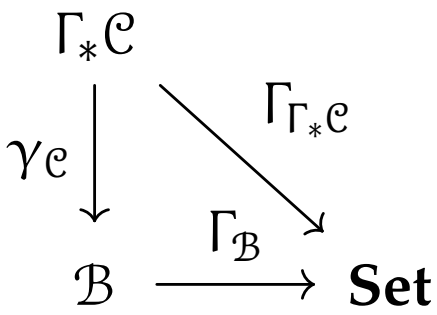

*Proof.* The functor $\Gamma_{\Gamma_*\mathcal{C}}$ sends an object $X \in |\Gamma_*\mathcal{C}| = |\mathcal{C}|$ to

$$(\Gamma_*\mathcal{C})(1_{\mathcal{C}}, X) = \mathcal{B}(1, \mathcal{C}(1_{\mathcal{C}}, X)) = \Gamma_{\mathcal{B}}(\gamma_{\mathcal{C}}(X)).$$

Similarly, for a morphism $f \in_1 \mathcal{C}(X, Y)$ in $\Gamma_*\mathcal{C}$, the functor $\Gamma_{\Gamma_*\mathcal{C}}$ sends $f$ to the map

$$(\Gamma_*\mathcal{C})(1_{\mathcal{C}}, X) \xrightarrow{f \circ (-)} (\Gamma_*\mathcal{C})(1_{\mathcal{C}}, Y).$$

By the definition of composition in $\Gamma_*\mathcal{C}$, this map coincides with

$$\mathcal{B}(1, \mathcal{C}(1_{\mathcal{C}}, X)) \xrightarrow{\gamma_{\mathcal{C}}(f) \circ (-)} \mathcal{B}(1, \mathcal{C}(1_{\mathcal{C}}, Y)).$$

This shows that the actions on morphisms also coincide, which completes the proof. □

## 2.4 Internal categories

The notion of an enriched category captures the idea of a category in which each hom-set is replaced by an object in another category. By contrast, the notion of an internal category is even stronger: not only the hom-objects but also the collections of all objects and morphisms are required to be objects of a category.

There is an intimate connection between internal categories and fibrations, but this topic will not be pursued here. We refer the reader to Jacobs [14, Chapter 7] for details.

### 2.4.1 Basic definitions

**Definition 2.32 (Internal category)** Let $\mathcal{B}$ be a category with pullbacks. A $\mathcal{B}$-*internal category* $\mathbb{C}$ consists of the following data:

- An object $|\mathbb{C}| \in \mathcal{B}$, called the *object of objects*.
- An object $|\mathbb{C}^{\to}| \in \mathcal{B}$, called the *object of morphisms*.
- A pair of morphisms $\mathrm{dom}, \mathrm{cod} : |\mathbb{C}^{\to}| \to |\mathbb{C}|$, called the *domain* and *codomain morphisms*.
- A morphism $e : |\mathbb{C}| \to |\mathbb{C}^{\to}|$, called the *identity-assigning morphism*.
- A morphism $\circ : |\mathbb{C}^{\to\to}| \to |\mathbb{C}^{\to}|$, called the *composition morphism*, where $|\mathbb{C}^{\to\to}|$ denotes the *object of composable pairs* defined by the following pullback:

$$
\begin{array}{ccc}
|\mathbb{C}^{\to\to}| & \xrightarrow{\pi_2} & |\mathbb{C}^{\to}| \\
{\scriptstyle \pi_1}\downarrow & \lrcorner & \downarrow{\scriptstyle \mathrm{dom}} \\
|\mathbb{C}^{\to}| & \xrightarrow{\mathrm{cod}} & |\mathbb{C}|.
\end{array}
$$

These data must satisfy the following commutative diagrams:

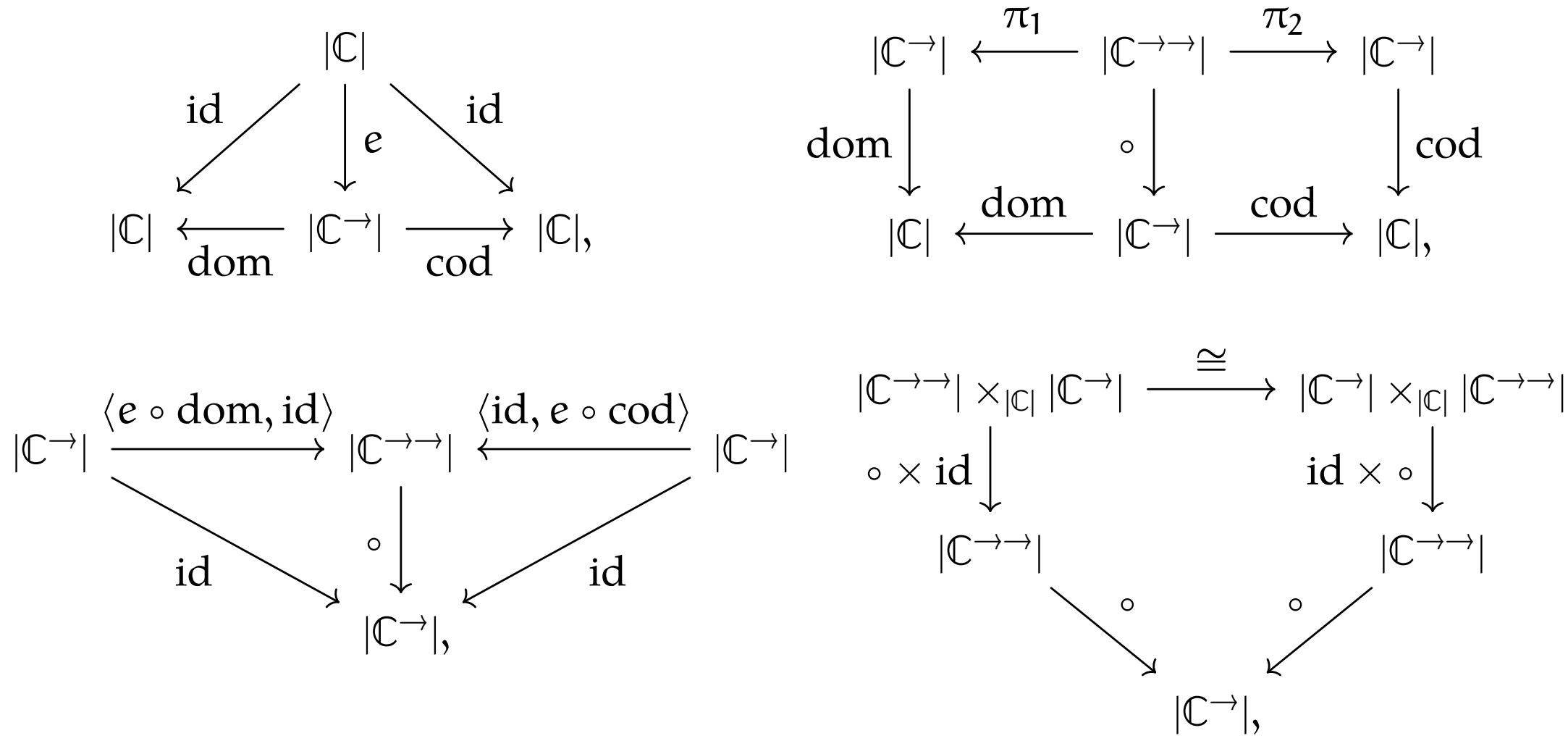

where the objects in the final diagram are defined by the following pullbacks:

$$
\begin{array}{ccc}
|\mathbb{C}^{\to\to}| \times_{|\mathbb{C}|} |\mathbb{C}^{\to}| & \xrightarrow{\pi_2} & |\mathbb{C}^{\to}| \\
{\scriptstyle \pi_1}\downarrow & \lrcorner & \downarrow{\scriptstyle \mathrm{dom}} \\
|\mathbb{C}^{\to\to}| & \xrightarrow{\mathrm{cod} \circ \pi_2} & |\mathbb{C}|,
\end{array}
\qquad
\begin{array}{ccc}
|\mathbb{C}^{\to}| \times_{|\mathbb{C}|} |\mathbb{C}^{\to\to}| & \xrightarrow{\pi_2} & |\mathbb{C}^{\to\to}| \\
{\scriptstyle \pi_1}\downarrow & \lrcorner & \downarrow{\scriptstyle \mathrm{dom} \circ \pi_1} \\
|\mathbb{C}^{\to}| & \xrightarrow{\mathrm{cod}} & |\mathbb{C}|.
\end{array}
$$

We often write $e_{\mathbb{C}}, \circ_{\mathbb{C}}$ for the morphisms $e, \circ$ to clarify the category. Note that a **Set**-internal category is simply a small category.

Since the definition of internal categories is involved, verifying its axioms for specific constructions is often straightforward but tedious. In this thesis, we provide only the basic structure for such constructions and omit the detailed verification.

**Definition 2.33 (Internal functor)** Let $\mathbb{C}$ and $\mathbb{D}$ be $\mathcal{B}$-internal categories. A *$\mathcal{B}$-internal functor* $\mathbb{H}$ from $\mathbb{C}$ to $\mathbb{D}$ consists of a pair of morphisms in $\mathcal{B}$,

$$|\mathbb{H}| : |\mathbb{C}| \to |\mathbb{D}|, \qquad |\mathbb{H}^{\to}| : |\mathbb{C}^{\to}| \to |\mathbb{D}^{\to}|,$$

such that they commute with the structural morphisms of $\mathbb{C}$ and $\mathbb{D}$ (i.e., the domain, codomain, identity-assigning, and composition morphisms).

### 2.4.2 Internal slice categories and hom-objects

When $\mathcal{B}$ has a terminal object 1, we write $X \in_1 \mathbb{C}$ for $X \in_1 |\mathbb{C}|$, called a *global object* of the internal category $\mathbb{C}$. For these global objects, one can construct slice categories and hom-objects analogously to ordinary categories.

**Definition 2.34 (Internal slice category)** Let $\mathbb{C}$ be a $\mathcal{B}$-internal category. For a global object $X \in_1 \mathbb{C}$, the $\mathcal{B}$-internal category $\mathbb{C}/X$ is constructed as follows:

**Object of objects** $|\mathbb{C}/X|$ is defined by the following pullback:

$$\begin{array}{ccc} |\mathbb{C}/X| & \xrightarrow{\iota_X} & |\mathbb{C}^{\to}| \\ {\scriptstyle !}\downarrow \lrcorner & & \downarrow{\scriptstyle \mathrm{cod}} \\ 1 & \xrightarrow{X} & |\mathbb{C}|. \end{array}$$

The projection $|\mathbb{C}/X| \to |\mathbb{C}^{\to}|$ is denoted by $\iota_X$ (or simply $\iota$).

**Object of morphisms** $\left|(\mathbb{C}/X)^{\to}\right|$ is defined by the following pullback:

$$\begin{array}{ccc} \left|(\mathbb{C}/X)^{\to}\right| & \xrightarrow{\pi_2} & |\mathbb{C}/X| \\ {\scriptstyle \pi_1}\downarrow \lrcorner & & \downarrow{\scriptstyle \mathrm{dom} \circ \iota_X\ (= \mathrm{dom}_X)} \\ |\mathbb{C}^{\to}| & \xrightarrow{\mathrm{cod}} & |\mathbb{C}|. \end{array}$$

**Domain morphism** The unique morphism $f$ making the following commute:

$$\begin{array}{ccc} \left|(\mathbb{C}/X)^{\to}\right| & \xrightarrow{\langle \pi_1, \iota_X \circ \pi_2 \rangle} & |\mathbb{C}^{\to\to}| \\ & {\scriptstyle f} \searrow & \downarrow{\scriptstyle \circ_{\mathbb{C}}} \\ & |\mathbb{C}/X| \xrightarrow{\iota_X} & |\mathbb{C}^{\to}| \\ {\scriptstyle !} & {\scriptstyle !}\downarrow \lrcorner & \downarrow{\scriptstyle \mathrm{cod}} \\ & 1 \xrightarrow{X} & |\mathbb{C}|. \end{array}$$

**Codomain morphism** The projection $\pi_2 : \left|(\mathbb{C}/X)^{\to}\right| \to |\mathbb{C}/X|$.

**Remark 2.35** The projection $\iota_X$ is regarded as an inclusion $\iota_X : |\mathbb{C}/X| \rightarrowtail |\mathbb{C}^{\to}|$. Indeed, since $X : 1 \to |\mathbb{C}|$ is a monomorphism, its pullback $\iota_X$ is also a monomorphism. While this fact is trivial, it plays a crucial role in the proof of Lemma 5.2. Additionally, the composite $\mathrm{dom} \circ \iota_X : |\mathbb{C}/X| \to |\mathbb{C}|$ is written as $\mathrm{dom}_X$.

We proceed to the construction of hom-objects.

**Definition 2.36 (Hom-object)** Let $\mathbb{C}$ be a $\mathcal{B}$-internal category. For $X, Y \in_1 \mathbb{C}$, the *hom-object* $\mathbb{C}(X, Y) \in \mathcal{B}$ is defined by the following pullback:

$$\begin{array}{ccc} \mathbb{C}(X, Y) & \xrightarrow{\iota_{X,Y}} & |\mathbb{C}^{\to}| \\ {\scriptstyle !}\downarrow \lrcorner & & \downarrow{\scriptstyle \langle \mathrm{dom}, \mathrm{cod} \rangle} \\ 1 & \xrightarrow{\langle X, Y \rangle} & |\mathbb{C}| \times |\mathbb{C}|. \end{array}$$

The projection $\mathbb{C}(X, Y) \to |\mathbb{C}^{\to}|$ is written as $\iota_{X,Y}$ (or simply $\iota$). There is also an inclusion $\iota'_{X,Y} : \mathbb{C}(X, Y) \to |\mathbb{C}/Y|$ defined by the commutativity of the following diagram:

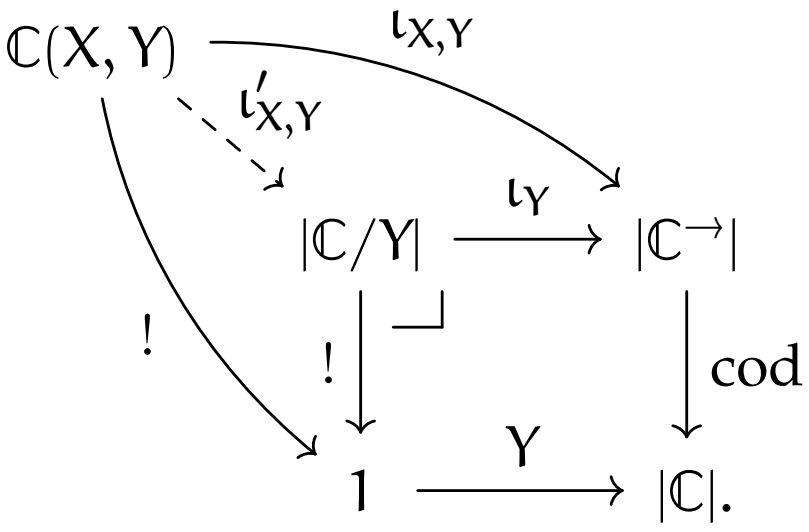

For any pair of global objects $X, Y \in_1 \mathbb{C}$, a global element $f \in_1 \mathbb{C}(X, Y)$ of the hom-object is called a *global morphism* of $\mathbb{C}$ from X to Y.

The composition morphism naturally induces morphisms between hom-objects:

**Definition 2.37** Let $\mathbb{C}$ be a $\mathcal{B}$-internal category. For global objects $X, Y, Z \in_1 \mathbb{C}$, the morphism $\circ_{X,Y,Z} : \mathbb{C}(X, Y) \times \mathbb{C}(Y, Z) \to \mathbb{C}(X, Z)$ is defined as the unique morphism such that the following commutes:

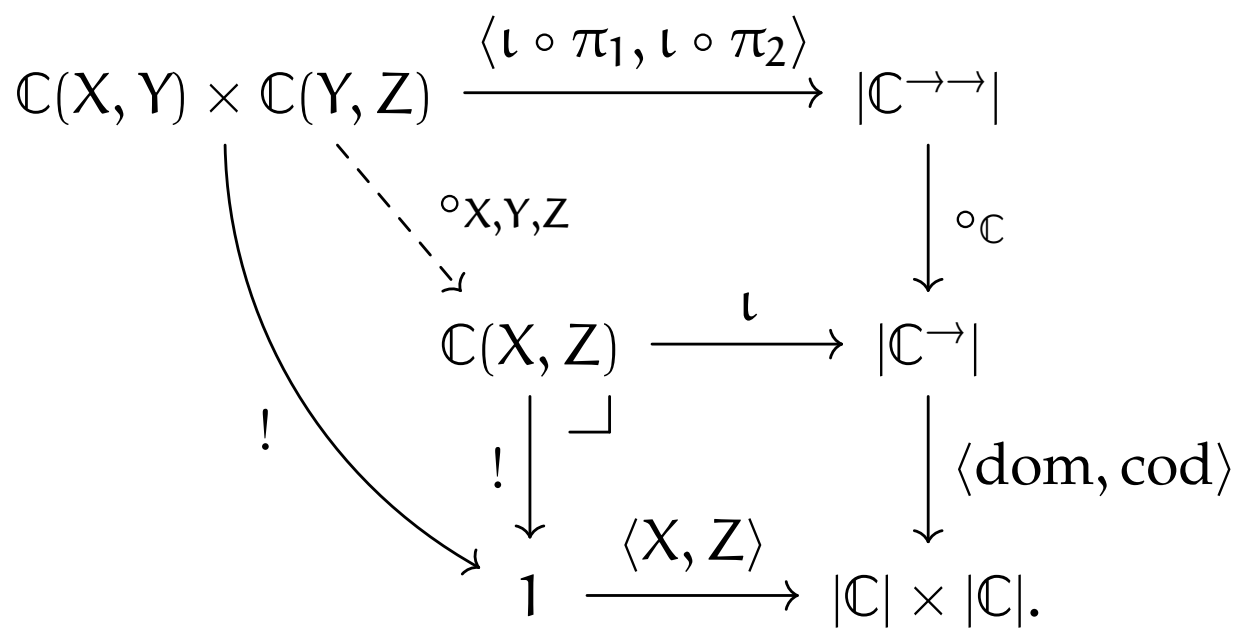

When there is no risk of confusion, we simply write $\circ$ or $\circ_{\mathbb{C}}$ for the morphism $\circ_{X,Y,Z}$.

Given a $\mathcal{B}$-internal functor $\mathbb{H} : \mathbb{C} \to \mathbb{D}$ and a global object $X \in_1 \mathbb{C}$, we write $\mathbb{H}X$ for the global object $|\mathbb{H}|(X) \in_1 \mathbb{D}$.

**Definition 2.38** Let $\mathbb{H} : \mathbb{C} \to \mathbb{D}$ be a $\mathcal{B}$-internal functor. For global objects $X, Y \in_1 \mathbb{C}$, the morphism $\mathbb{H}_{X,Y} : \mathbb{C}(X, Y) \to \mathbb{D}(\mathbb{H}X, \mathbb{H}Y)$ is defined as the unique morphism making the following diagram commute:

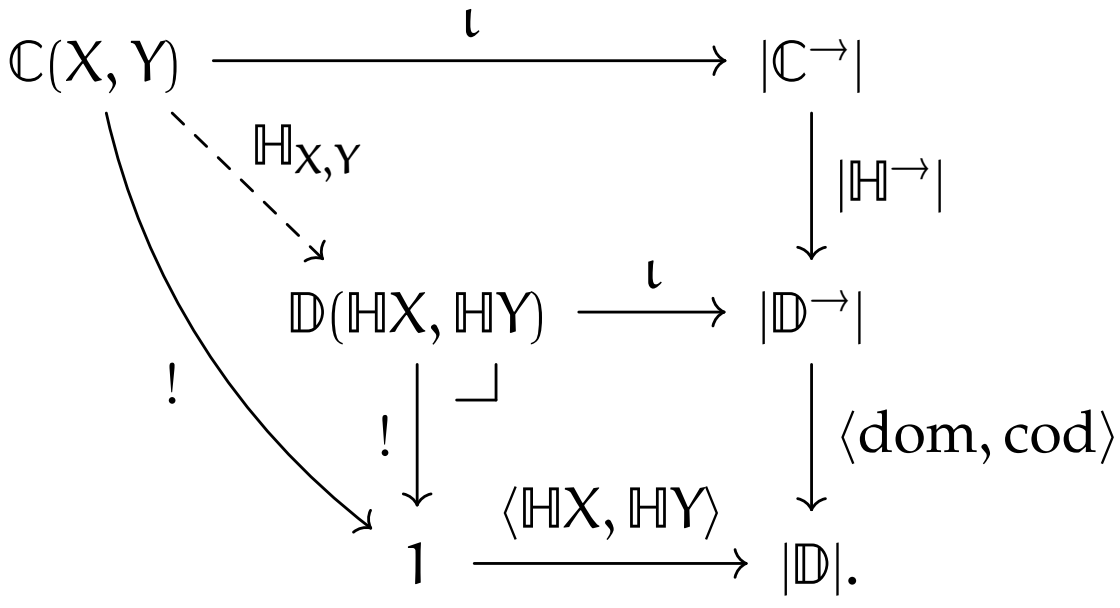

### 2.4.3 Change of base and underlying categories

There is an internal analog of the change of base for enriched categories (Definition 2.21), which we simply refer to as the *change of base.*

**Definition 2.39 (Change of base)** Let $\mathbb{C}$ be a $\mathcal{B}$-internal category and $F : \mathcal{B} \to \mathcal{B}'$ be a functor preserving pullbacks. Applying the functor $F$ to all objects and morphisms constituting $\mathbb{C}$ yields a $\mathcal{B}'$-internal category. This internal category is denoted as $F\mathbb{C}$ and called the *change of base* of $\mathbb{C}$ along $F$.

In other words, for a $\mathcal{B}$-internal category $\mathbb{C} = (|\mathbb{C}|, |\mathbb{C}^{\to}|, \mathrm{dom}, \mathrm{cod}, e, \circ)$, we define the $\mathcal{B}'$-internal category $F\mathbb{C}$ as $(F|\mathbb{C}|, F|\mathbb{C}^{\to}|, F(\mathrm{dom}), F(\mathrm{cod}), F(e), F(\circ))$, up to canonical isomorphisms between certain pullbacks.

**Definition 2.40 (Change of base for internal functors)** Let $\mathbb{H} : \mathbb{C} \to \mathbb{D}$ be a $\mathcal{B}$-internal functor. For a functor $F : \mathcal{B} \to \mathcal{B}'$ preserving pullbacks, the $\mathcal{B}'$-internal functor $F\mathbb{H} : F\mathbb{C} \to F\mathbb{D}$ is defined by $|F\mathbb{H}| = F|\mathbb{H}|$ and $\left|(F\mathbb{H})^{\to}\right| = F|\mathbb{H}^{\to}|$.

By applying the change of base construction along the global section functor, one can obtain ordinary categories from internal categories. Following the terminology for enriched categories, we call these *underlying categories.*

**Definition 2.41 (Underlying category)** Let $\mathcal{B}$ be a locally small category with finite limits. For a $\mathcal{B}$-internal category $\mathbb{C}$, its *underlying category* $\Gamma\mathbb{C}$ is the change of base of $\mathbb{C}$ along the global section functor $\Gamma : \mathcal{B} \to \mathbf{Set}$. Explicitly, the category $\Gamma\mathbb{C}$ is described as follows:

**Objects** Global objects $X \in_1 \mathbb{C}$.

**Morphisms** $X \to Y$ are global elements $f \in_1 |\mathbb{C}^{\to}|$ such that $\mathrm{dom}(f) = X$ and $\mathrm{cod}(f) = Y$. Such global elements are in one-to-one correspondence with global morphisms from $X$ to $Y$ via the following diagram, hence we identify them:

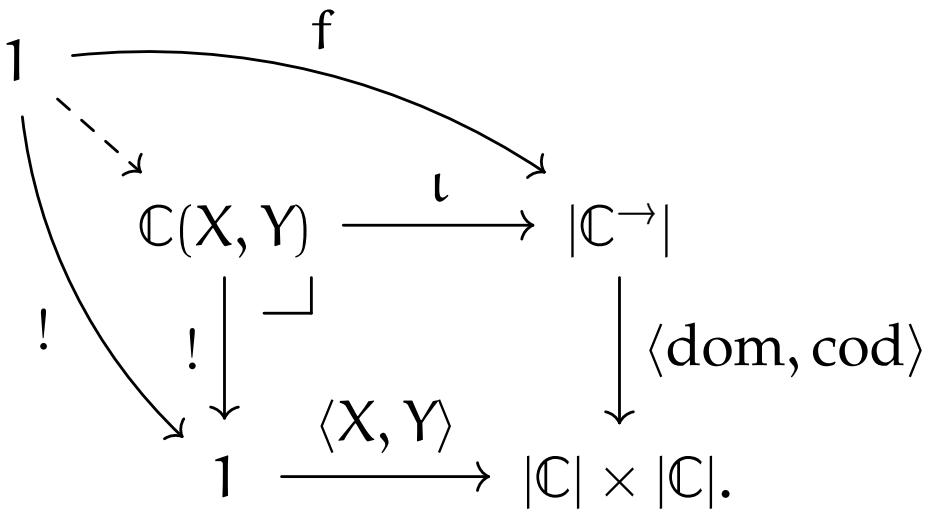

**Identities and composition** Induced by $e_{\mathbb{C}}$ and $\circ_{\mathbb{C}}$, respectively.

The underlying category of $\mathbb{C}$ is isomorphic to the fiber of the externalization $\int \mathbb{C} \to \mathcal{B}$ of $\mathbb{C}$ at the terminal object $1$ (see, e.g., Jacobs [14, Section 7.3]). Whenever we consider the underlying category of a $\mathcal{B}$-internal category, we implicitly assume that $\mathcal{B}$ is locally small and has finite limits.

**Definition 2.42 (Underlying functor)** Let $\mathcal{B}$ be a category with finite limits and $\mathbb{H} : \mathbb{C} \to \mathbb{D}$ be a $\mathcal{B}$-internal functor. The *underlying functor* of $\mathbb{H}$, denoted by $\Gamma\mathbb{H}$, is the change of base of $\mathbb{H}$ along $\Gamma : \mathcal{B} \to \mathbf{Set}$. Explicitly, the functor $\Gamma\mathbb{H} : \Gamma\mathbb{C} \to \Gamma\mathbb{D}$ is described as follows:

**On objects** For a global object $X \in_1 \mathbb{C}$, $(\Gamma\mathbb{H})X = \mathbb{H}X$.

**On morphisms** For a global morphism $f \in_1 \mathbb{C}(X, Y)$, we define $(\Gamma\mathbb{H})(f) \in_1 \mathbb{D}(\mathbb{H}X, \mathbb{H}Y)$ as $\mathbb{H}_{X,Y}(f)$, i.e., the following composite:

$$1 \xrightarrow{f} \mathbb{C}(X, Y) \xrightarrow{\mathbb{H}_{X,Y}} \mathbb{D}(\mathbb{H}X, \mathbb{H}Y).$$

**Definition 2.43 (Hom-functor)** Let $\mathbb{C}$ be a $\mathcal{B}$-internal category. For any global object $X \in_1 \mathbb{C}$, a *hom-functor* $\mathbb{C}(X, -) : \Gamma\mathbb{C} \to \mathcal{B}$ is defined as follows:

**On objects** For a global object $Y \in_1 \mathbb{C}$, the functor maps $Y$ to the hom-object $\mathbb{C}(X, Y)$.

**On morphisms** For a global morphism $f \in_1 \mathbb{C}(Y, Z)$, the morphism $\mathbb{C}(X, f) : \mathbb{C}(X, Y) \to \mathbb{C}(X, Z)$ is defined by the following composite:

$$\mathbb{C}(X, Y) \xrightarrow{\langle \mathrm{id}, f \circ ! \rangle} \mathbb{C}(X, Y) \times \mathbb{C}(Y, Z) \xrightarrow{\circ_{X,Y,Z}} \mathbb{C}(X, Z).$$

The following definition provides an internal analog of $F^{\#e}$ for enriched categories (Definition 2.26).

**Definition 2.44 ($F^{\#}$)** Let $F : \mathcal{B} \to \mathcal{B}'$ be a functor preserving finite limits. For a $\mathcal{B}$-internal category $\mathbb{C}$, we define the functor $F^{\#}_{\mathbb{C}} : \Gamma\mathbb{C} \to \Gamma F\mathbb{C}$ (or simply $F^{\#}$) as follows:

**On objects** For each global object $X \in_1 \mathbb{C}$, the global object $F^{\#}X \in_1 F\mathbb{C}$ is defined by the following composite:

$$1_{\mathcal{B}'} \xrightarrow{\cong} F(1_{\mathcal{B}}) \xrightarrow{F(X)} F|\mathbb{C}|.$$

**On morphisms** For each global morphism $f \in_1 \mathbb{C}(X, Y)$, the global morphism $F^{\#}f \in_1 (F\mathbb{C})\big(F^{\#}X, F^{\#}Y\big)$ is defined by the following composite:

$$1_{\mathcal{B}'} \xrightarrow{\cong} F(1_{\mathcal{B}}) \xrightarrow{F(f)} F(\mathbb{C}(X, Y)) \xrightarrow{\cong} (F\mathbb{C})\big(F^{\#}X, F^{\#}Y\big),$$

where the last isomorphism is induced by the fact that $F$ preserves pullbacks.

**Proposition 2.45** *Let $F : \mathcal{B} \to \mathcal{B}'$ be a functor preserving finite limits and $\mathbb{H} : \mathbb{C} \to \mathbb{D}$ be a $\mathcal{B}$-internal functor. Then, we have $\Gamma F\mathbb{H} \circ F^{\#}_{\mathbb{C}} = F^{\#}_{\mathbb{D}} \circ \Gamma\mathbb{H}$.*

$$\begin{array}{ccc} \Gamma\mathbb{C} & \xrightarrow{F^{\#}_{\mathbb{C}}} & \Gamma F\mathbb{C} \\ {\scriptstyle \Gamma\mathbb{H}}\downarrow & & \downarrow{\scriptstyle \Gamma F\mathbb{H}} \\ \Gamma\mathbb{D} & \xrightarrow{F^{\#}_{\mathbb{D}}} & \Gamma F\mathbb{D} \end{array}$$

*Proof.* This is straightforward from the definitions. □

### 2.4.4 Terminal objects and global section functors

In this and the next subsections, we provide an elementary description of terminal objects and pullbacks in internal categories. For the general notion of limits in internal categories, we refer to Borceux [4, Section 8.3] or Jacobs [14, Section 7.4].

**Definition 2.46 (Terminal object)** For a $\mathcal{B}$-internal category $\mathbb{C}$, a *terminal object* of $\mathbb{C}$ consists of a global object $1_{\mathbb{C}} \in_1 \mathbb{C}$ and a morphism $!_{\mathbb{C}} : |\mathbb{C}| \to |\mathbb{C}^{\to}|$ such that the following diagrams commute:

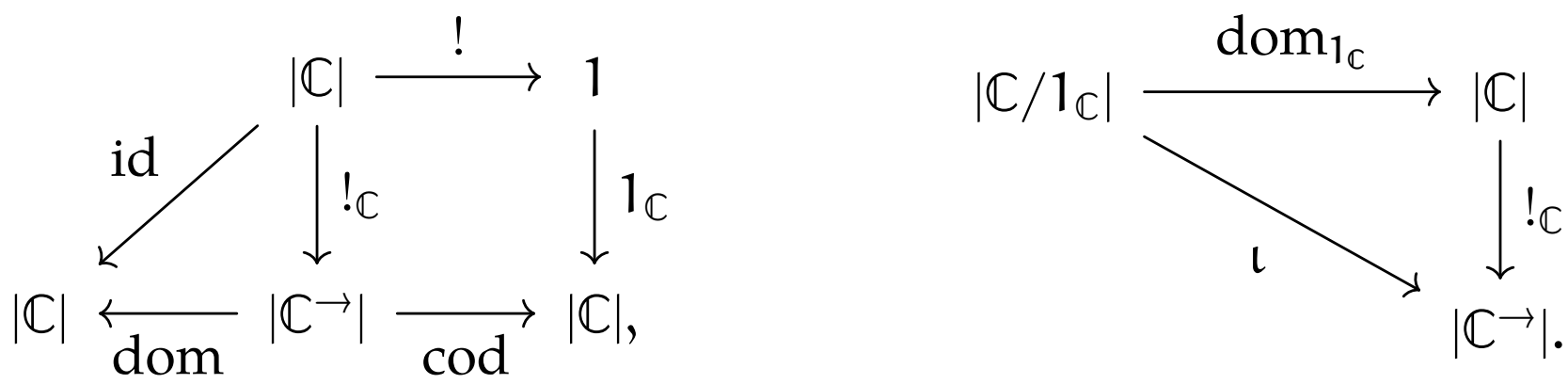

To avoid confusion with ordinary categories, we always explicitly include the subscripts in the notations $1_{\mathbb{C}}$ and $!_{\mathbb{C}}$.

As remarked in Section 2.1.2, we regard a terminal object in an internal category as a chosen structure. To avoid confusion, we may use the explicit phrase *internal categories with chosen terminal objects.* An internal functor *strictly preserves* the terminal object if it maps the morphisms constituting the terminal object in the domain to the corresponding morphisms in the codomain.

**Proposition 2.47** *Let $\mathcal{B}, \mathcal{B}'$ be categories with finite limits and $F : \mathcal{B} \to \mathcal{B}'$ be a functor preserving finite limits.*

1. *If a $\mathcal{B}$-internal category $\mathbb{C}$ has a chosen terminal object, then its change of base $F\mathbb{C}$ is naturally equipped with a terminal object induced by $F$.*
2. *If a $\mathcal{B}$-internal functor $\mathbb{H} : \mathbb{C} \to \mathbb{D}$ strictly preserves the terminal object, then its change of base $F\mathbb{H} : F\mathbb{C} \to F\mathbb{D}$ also strictly preserves it.*

*Proof.* (1) Since $F$ preserves the terminal object and the commutativity of diagrams, the pair $(F(1_{\mathbb{C}}), F(!_{\mathbb{C}}))$ satisfies the required commutativity for a terminal object in $F\mathbb{C}$ (after precomposing with the isomorphism $1_{\mathcal{B}'} \cong F(1_{\mathcal{B}})$). Thus, these data yield a terminal object of $F\mathbb{C}$.

(2) The terminal objects of $F\mathbb{C}$ and $F\mathbb{D}$ are those induced by (1). Since $\mathbb{H}$ strictly preserves the terminal object, we have

$$(F\mathbb{H})(1_{F\mathbb{C}}) = (F\mathbb{H})(F(1_{\mathbb{C}})) = F(\mathbb{H}(1_{\mathbb{C}})) = F(1_{\mathbb{D}}) = 1_{F\mathbb{D}},$$

and the similar equality holds for $!_{F\mathbb{C}}$. Thus, $F\mathbb{H}$ strictly preserves the terminal object. □

**Proposition 2.48** *For a $\mathcal{B}$-internal category $\mathbb{C}$ with a chosen terminal object, any hom-functor $\mathbb{C}(X, -) : \Gamma\mathbb{C} \to \mathcal{B}$ preserves the terminal object.*

*Proof.* We prove $\mathbb{C}(X, 1_\mathbb{C}) \cong 1$ for any global object $X \in_1 \mathbb{C}$. First, composing $!_\mathbb{C} : |\mathbb{C}| \to |\mathbb{C}^\to|$ with $X$ yields a global element of $|\mathbb{C}^\to|$, which corresponds to a global morphism $!_X \in_1 \mathbb{C}(X, 1_\mathbb{C})$. It suffices to show that the following commutes:

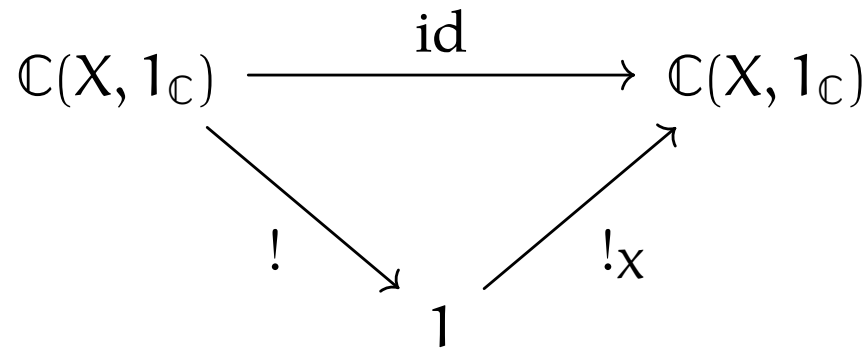

Since $\iota : \mathbb{C}(X, 1_\mathbb{C}) \to |\mathbb{C}^\to|$ is a monomorphism, it is enough to show that $\iota = \iota \circ {!_X} \circ {!}$. This follows from the equality $\iota \circ \iota' = \iota$ and the commutativity of the following diagram:

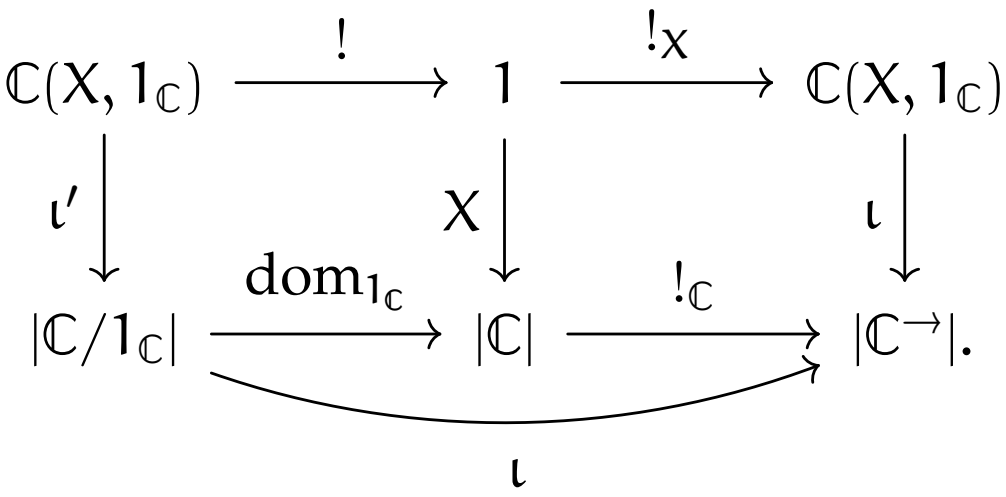

This completes the proof. ☐

**Definition 2.49 (Global section functor)** Let $\mathbb{C}$ be a $\mathcal{B}$-internal category with a chosen terminal object $1_\mathbb{C} \in_1 \mathbb{C}$. We call the hom-functor $\mathbb{C}(1_\mathbb{C}, -) : \Gamma\mathbb{C} \to \mathcal{B}$ the *global section functor* of $\mathbb{C}$ and write it as $\gamma_\mathbb{C}$ (or simply $\gamma$).

**Proposition 2.50 ($\Gamma \circ \gamma \cong \Gamma$)** *For a $\mathcal{B}$-internal category $\mathbb{C}$ with a chosen terminal object, there is a natural isomorphism $\Gamma_\mathcal{B} \circ \gamma_\mathbb{C} \cong \Gamma_{\Gamma\mathbb{C}}$.*

$$\begin{array}{ccc} \Gamma\mathbb{C} & & \\ {\scriptstyle \gamma_\mathbb{C}}\downarrow & \searrow^{\Gamma_{\Gamma\mathbb{C}}} & \\ \mathcal{B} & \xrightarrow{\Gamma_\mathcal{B}} & \textbf{Set} \end{array}$$

*Proof.* The global section functor $\Gamma_{\Gamma\mathbb{C}}$ sends each global object $X \in_1 \mathbb{C}$ to the set $(\Gamma\mathbb{C})(1_\mathbb{C}, X) \cong \mathcal{B}(1, \mathbb{C}(1_\mathbb{C}, X)) = \Gamma_\mathcal{B}(\gamma_\mathbb{C}(X))$. This isomorphism is natural in $X$. ☐

**Proposition 2.51** *Let $\mathcal{B}, \mathcal{B}'$ be categories with finite limits and $F : \mathcal{B} \to \mathcal{B}'$ be a functor preserving finite limits. For a $\mathcal{B}$-internal category $\mathbb{C}$, there is a natural isomorphism $\gamma_{F\mathbb{C}} \circ F^\#_\mathbb{C} \cong F \circ \gamma_\mathbb{C}$.*

$$\begin{array}{ccc} \Gamma\mathbb{C} & \xrightarrow{F^\#} & \Gamma F\mathbb{C} \\ {\scriptstyle \gamma_\mathbb{C}}\downarrow & & \downarrow{\scriptstyle \gamma_{F\mathbb{C}}} \\ \mathcal{B} & \xrightarrow{F} & \mathcal{B}' \end{array}$$

*Proof.* For any global object $X \in_1 \mathbb{C}$, we have

$$\gamma_{F\mathbb{C}}\big(F^\# X\big) = (F\mathbb{C})\big(1_{F\mathbb{C}}, F^\# X\big) = (F\mathbb{C})\big(F^\# 1_{\mathbb{C}}, F^\# X\big) \cong F(\mathbb{C}(1_{\mathbb{C}}, X)) = F(\gamma_{\mathbb{C}}(X)).$$

This isomorphism is natural in $X$. □

### 2.4.5 Pullbacks

We start with some preparations to define pullbacks in internal categories.

**Definition 2.52** Let $\mathbb{C}$ be a $\mathcal{B}$-internal category.

(1) The *object of cospans* in $\mathbb{C}$, written $|\mathbb{C}^{\to\leftarrow}|$, is defined by the following pullback:

$$\begin{array}{ccc} |\mathbb{C}^{\to\leftarrow}| & \xrightarrow{\pi_2} & |\mathbb{C}^{\to}| \\ {\scriptstyle \pi_1}\downarrow & \lrcorner & \downarrow{\scriptstyle \mathrm{cod}} \\ |\mathbb{C}^{\to}| & \xrightarrow{\mathrm{cod}} & |\mathbb{C}|. \end{array}$$

(2) The *object of spans* in $\mathbb{C}$, written $|\mathbb{C}^{\leftarrow\to}|$, is defined by the following pullback:

$$\begin{array}{ccc} |\mathbb{C}^{\leftarrow\to}| & \xrightarrow{\pi_2} & |\mathbb{C}^{\to}| \\ {\scriptstyle \pi_1}\downarrow & \lrcorner & \downarrow{\scriptstyle \mathrm{dom}} \\ |\mathbb{C}^{\to}| & \xrightarrow{\mathrm{dom}} & |\mathbb{C}|. \end{array}$$

(3) The *object of commutative squares,* written $|\mathbb{C}^{sq}|$, is defined by the following pullback:

$$\begin{array}{ccc} |\mathbb{C}^{sq}| & \xrightarrow{\pi_2} & |\mathbb{C}^{\to\to}| \\ {\scriptstyle \pi_1}\downarrow & \lrcorner & \downarrow{\scriptstyle \circ} \\ |\mathbb{C}^{\to\to}| & \xrightarrow{\circ} & |\mathbb{C}^{\to}|. \end{array}$$

Additionally, we define four projections $|\mathbb{C}^{sq}| \to |\mathbb{C}|$ by

$$\pi_a = \mathrm{dom} \circ \pi_1 \circ \pi_1 = \mathrm{dom} \circ \pi_1 \circ \pi_2, \qquad \pi_b = \mathrm{cod} \circ \pi_1 \circ \pi_1 = \mathrm{dom} \circ \pi_2 \circ \pi_1,$$
$$\pi_c = \mathrm{cod} \circ \pi_1 \circ \pi_2 = \mathrm{dom} \circ \pi_2 \circ \pi_2, \qquad \pi_d = \mathrm{cod} \circ \pi_2 \circ \pi_1 = \mathrm{cod} \circ \pi_2 \circ \pi_2,$$

and also four projections $|\mathbb{C}^{sq}| \to |\mathbb{C}^{\to}|$ by

$$\pi_{ab} = \pi_1 \circ \pi_1, \qquad \pi_{bd} = \pi_2 \circ \pi_1, \qquad \pi_{ac} = \pi_1 \circ \pi_2, \qquad \pi_{cd} = \pi_2 \circ \pi_2.$$

An intuitive picture is as follows:

$$\begin{array}{ccc} a & \longrightarrow & c \\ \downarrow & & \downarrow \\ b & \longrightarrow & d. \end{array}$$

**Definition 2.53 (Pullbacks)** For a $\mathcal{B}$-internal category $\mathbb{C}$, *pullbacks* in $\mathbb{C}$ consist of morphisms $\theta : |\mathbb{C}^{\to\leftarrow}| \to |\mathbb{C}^{sq}|$ and $\varphi : |\mathbb{C}^{sq}| \to |\mathbb{C}^{\to}|$ such that the following diagrams commute:

$$
\begin{array}{ccc}
|\mathbb{C}^{\to\leftarrow}| & \xrightarrow{\ \theta\ } & |\mathbb{C}^{sq}| \\
& \underset{\mathrm{id}}{\searrow} & \Big\downarrow{\scriptstyle \langle \pi_{bd}, \pi_{cd} \rangle} \\
& & |\mathbb{C}^{\to\leftarrow}|,
\end{array}
\qquad
\begin{array}{ccccc}
& & |\mathbb{C}^{sq}| & \xrightarrow{\ \overline{\theta}\ } & |\mathbb{C}^{sq}| \\
& \overset{\pi_a}{\swarrow} & \Big\downarrow{\scriptstyle \varphi} & & \Big\downarrow{\scriptstyle \pi_a} \\
|\mathbb{C}| & \xleftarrow[\mathrm{dom}]{} & |\mathbb{C}^{\to}| & \xrightarrow[\mathrm{cod}]{} & |\mathbb{C}|,
\end{array}
$$

$$
\begin{array}{ccccc}
|\mathbb{C}^{sq}| & \xrightarrow{\ \langle \varphi, \pi_{ab} \circ \overline{\theta} \rangle\ } & |\mathbb{C}^{\to\to}| & \xleftarrow{\ \langle \varphi, \pi_{ac} \circ \overline{\theta} \rangle\ } & |\mathbb{C}^{sq}| \\
& \underset{\pi_{ab}}{\searrow} & \Big\downarrow{\scriptstyle \circ} & \underset{\pi_{ac}}{\swarrow} & \\
& & |\mathbb{C}^{\to}|, & &
\end{array}
$$

$$
\begin{array}{ccc}
|\mathbb{C}^{\to}| \times_{|\mathbb{C}|} |\mathbb{C}^{\to\leftarrow}| & \xrightarrow{\ \langle \langle f_{ab}, f_{bd} \rangle, \langle f_{ac}, f_{cd} \rangle \rangle\ } & |\mathbb{C}^{sq}| \\
& \underset{\pi_1}{\searrow} & \Big\downarrow{\scriptstyle \varphi} \\
& & |\mathbb{C}^{\to}|.
\end{array}
$$

Here, the object $|\mathbb{C}^{\to}| \times_{|\mathbb{C}|} |\mathbb{C}^{\to\leftarrow}|$ is defined by the following pullbacks:

$$
\begin{array}{ccc}
|\mathbb{C}^{\to}| \times_{|\mathbb{C}|} |\mathbb{C}^{\to\leftarrow}| & \xrightarrow{\ \pi_2\ } & |\mathbb{C}^{\to\leftarrow}| \\
{\scriptstyle \pi_1}\Big\downarrow \ \lrcorner & & \Big\downarrow{\scriptstyle \pi_a \circ \theta} \\
|\mathbb{C}^{\to}| & \xrightarrow{\ \mathrm{cod}\ } & |\mathbb{C}|.
\end{array}
$$

The morphisms $\overline{\theta}, f_{ab}, f_{bd}, f_{ac}, f_{cd}$ are defined by the following composites:

$$
\overline{\theta} = \left( |\mathbb{C}^{sq}| \xrightarrow{\ \langle \pi_{bd}, \pi_{cd} \rangle\ } |\mathbb{C}^{\to\leftarrow}| \xrightarrow{\ \theta\ } |\mathbb{C}^{sq}| \right),
$$

$$
f_{ab} = \left( |\mathbb{C}^{\to}| \times_{|\mathbb{C}|} |\mathbb{C}^{\to\leftarrow}| \xrightarrow{\ \langle \pi_1, \pi_{ab} \circ \theta \circ \pi_2 \rangle\ } |\mathbb{C}^{\to\to}| \xrightarrow{\ \circ\ } |\mathbb{C}^{\to}| \right),
$$

$$
f_{bd} = \left( |\mathbb{C}^{\to}| \times_{|\mathbb{C}|} |\mathbb{C}^{\to\leftarrow}| \xrightarrow{\ \pi_2\ } |\mathbb{C}^{\to\leftarrow}| \xrightarrow{\ \pi_1\ } |\mathbb{C}^{\to}| \right),
$$

$$
f_{ac} = \left( |\mathbb{C}^{\to}| \times_{|\mathbb{C}|} |\mathbb{C}^{\to\leftarrow}| \xrightarrow{\ \langle \pi_1, \pi_{ac} \circ \theta \circ \pi_2 \rangle\ } |\mathbb{C}^{\to\to}| \xrightarrow{\ \circ\ } |\mathbb{C}^{\to}| \right),
$$

$$
f_{cd} = \left( |\mathbb{C}^{\to}| \times_{|\mathbb{C}|} |\mathbb{C}^{\to\leftarrow}| \xrightarrow{\ \pi_2\ } |\mathbb{C}^{\to\leftarrow}| \xrightarrow{\ \pi_2\ } |\mathbb{C}^{\to}| \right).
$$

Intuitively, the morphism $\theta : |\mathbb{C}^{\to\leftarrow}| \to |\mathbb{C}^{sq}|$ corresponds to the following operation:

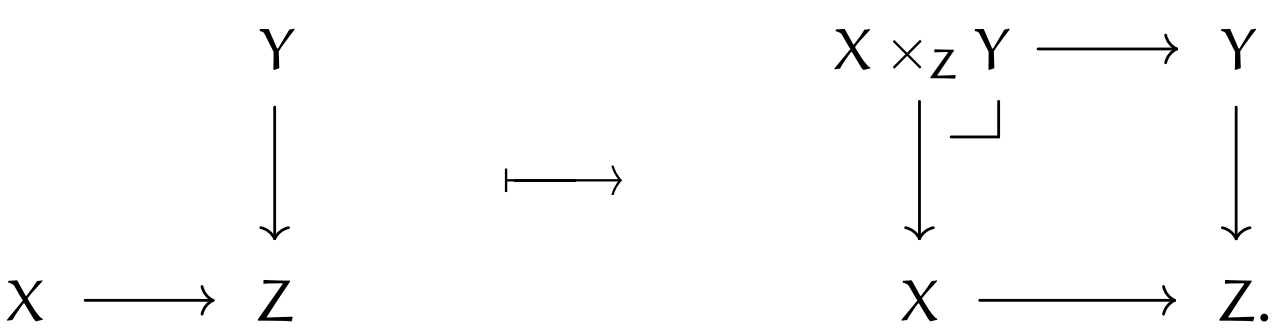

The morphism $\varphi : |\mathbb{C}^{sq}| \to |\mathbb{C}^{\to}|$ corresponds to the operation which maps a commutative square to the uniquely determined morphism $W \to X \times_Z Y$, depicted as follows:

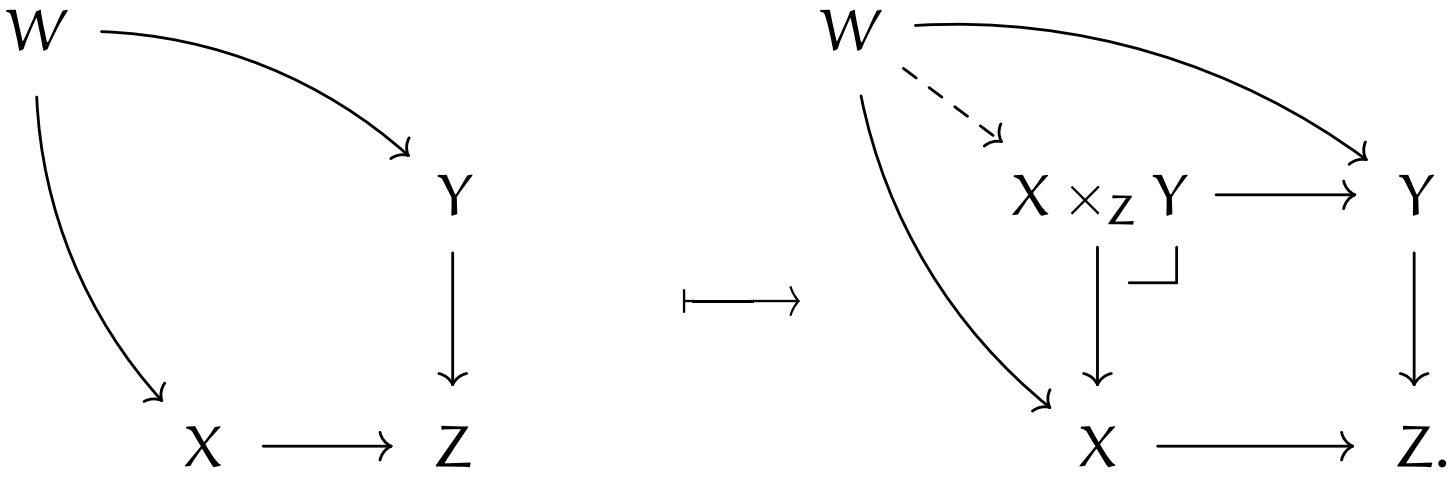

As before, we always regard pullbacks in internal categories as specifically chosen structures. We say that $\mathbb{C}$ has chosen *finite limits* if it is equipped with a terminal object and pullbacks.

The following definition provides an internal analog of the reindexing for codomain fibrations, $(-) \cdot (-) : |\mathcal{B}/A| \times \mathcal{B}(B, A) \to |\mathcal{B}/B|$.

**Definition 2.54** Let $\mathbb{C}$ be a $\mathcal{B}$-internal category with chosen pullbacks. For global objects $X, Y \in_1 \mathbb{C}$, the morphism $\mathrm{pb}_{X,Y} : |\mathbb{C}/Y| \times \mathbb{C}(X, Y) \to |\mathbb{C}/X|$ is defined as the unique morphism making the following commute:

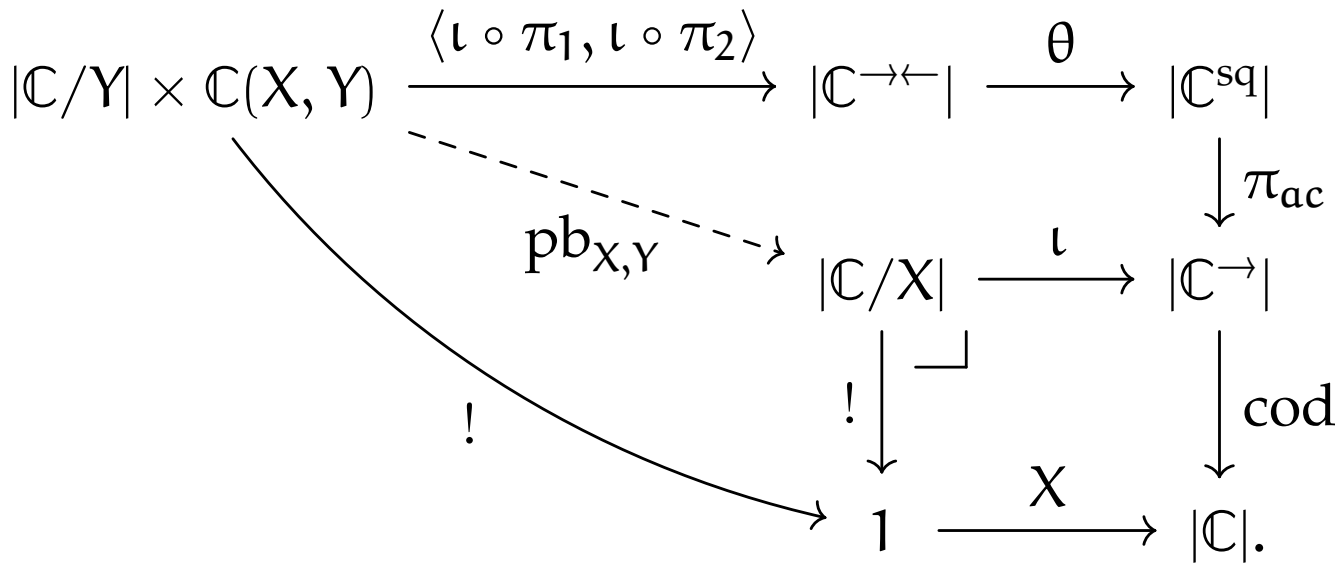

**Proposition 2.55** *Let $\mathcal{B}, \mathcal{B}'$ be categories with finite limits and $F : \mathcal{B} \to \mathcal{B}'$ be a functor preserving finite limits.*

(1) *If a $\mathcal{B}$-internal category $\mathbb{C}$ has chosen pullbacks, then its change of base $F\mathbb{C}$ is naturally equipped with pullbacks induced by $F$.*

(2) *If a $\mathcal{B}$-internal functor $\mathbb{H} : \mathbb{C} \to \mathbb{D}$ strictly preserves pullbacks, then its change of base $F\mathbb{H} : F\mathbb{C} \to F\mathbb{D}$ also strictly preserves them.*

*Proof.* The proof is similar to that of Proposition 2.47. □

In particular, if $\mathbb{C}$ has finite limits, then its underlying category $\Gamma\mathbb{C}$ also has finite limits.

**Proposition 2.56** *If a $\mathcal{B}$-internal category $\mathbb{C}$ has chosen pullbacks, then any hom-functor $\mathbb{C}(X, -) : \Gamma\mathbb{C} \to \mathcal{B}$ preserves pullbacks.*

*Proof.* This is shown by an argument similar to that of Proposition 2.48, though the details are straightforward but tedious. We omit the proof. □

The following is the internal analog of Proposition 2.29:

**Proposition 2.57** *Let $\mathbb{C}$ be a $\mathcal{B}$-internal category with chosen finite limits. For any global objects $W, X, Y, Z \in_1 \mathbb{C}$, the following diagram commutes:*

$$\begin{array}{ccc} \mathbb{C}(W,X) \times \mathbb{C}(X, Y \times Z) & \xrightarrow{\cong} & \mathbb{C}(W,X) \times (\mathbb{C}(X,Y) \times \mathbb{C}(X,Z)) \\ & & \downarrow \langle\langle \pi_1, \pi_1 \circ \pi_2\rangle, \langle \pi_1, \pi_2 \circ \pi_2\rangle\rangle \\ \downarrow \circ_{\mathbb{C}} & & (\mathbb{C}(W,X) \times \mathbb{C}(X,Y)) \times (\mathbb{C}(W,X) \times \mathbb{C}(X,Z)) \\ & & \downarrow \circ_{\mathbb{C}} \times \circ_{\mathbb{C}} \\ \mathbb{C}(W, Y \times Z) & \xrightarrow{\cong} & \mathbb{C}(W,Y) \times \mathbb{C}(W,Z). \end{array}$$

*Here, the global object $Y \times Z \in_1 \mathbb{C}$ is the binary product of $Y \in_1 \mathbb{C}$ and $Z \in_1 \mathbb{C}$ in the underlying category $\Gamma\mathbb{C}$.*

*Proof.* We define the $\mathcal{B}$-enriched category $\mathbb{C}_{\mathrm{enr}}$ as follows:

**Objects** Objects of $\mathbb{C}_{\mathrm{enr}}$ are global objects $X \in_1 \mathbb{C}$.

**Hom-objects** For $X, Y \in_1 \mathbb{C}$, we define $\mathbb{C}_{\mathrm{enr}}(X, Y) = \mathbb{C}(X, Y)$.

Then, it can be checked straightforwardly that $\mathbb{C}_{\mathrm{enr}}$ has finite products induced by the finite limits in $\mathbb{C}$. The proposition follows from applying Proposition 2.29 to the $\mathcal{B}$-enriched category $\mathbb{C}_{\mathrm{enr}}$. □

### 2.4.6 Doubly internal categories

We consider the notion of an internal category within another internal category. For an $\mathcal{A}$-internal category $\mathbb{B}$, the notion of an "internal category in $\mathbb{B}$" is simply identified with that of an ordinary $\Gamma\mathbb{B}$-internal category. In what follows, we provide analogs of some notions for internal categories in this context.

**Definition 2.58 ($\mathbb{F}^{\#}$)** Let $\mathbb{B}, \mathbb{B}'$ be an $\mathcal{A}$-internal category with chosen finite limits and $\mathbb{F} : \mathbb{B} \to \mathbb{B}'$ be an $\mathcal{A}$-internal functor strictly preserving finite limits. For a $\Gamma\mathbb{B}$-internal category $\mathbb{C}$, we define the $\mathcal{A}$-internal functor $\mathbb{F}^{\#}_{\mathbb{C}} : \gamma_{\mathbb{B}}\mathbb{C} \to \gamma_{\mathbb{B}'}(\Gamma\mathbb{F})\mathbb{C}$ (or simply $\mathbb{F}^{\#}$) as follows.

**On the object of objects** We define $|\mathbb{F}^{\#}|$ to be $\mathbb{F}_{1_{\mathbb{B}},|\mathbb{C}|}$. Note that $|\gamma_{\mathbb{B}}\mathbb{C}| = \mathbb{B}(1_{\mathbb{B}}, |\mathbb{C}|)$ and $|\gamma_{\mathbb{B}'}(\Gamma\mathbb{F})\mathbb{C}| = \mathbb{B}'(1_{\mathbb{B}'}, \mathbb{F}|\mathbb{C}|)$.

**On the object of morphisms** Similarly, $\left|(\mathbb{F}^{\#})^{\to}\right|$ is defined to be $\mathbb{F}_{1_{\mathbb{B}},|\mathbb{C}^{\to}|}$.

$$
\begin{array}{ccccccc}
 & \mathcal{A} & \ni & \mathbb{B} \xrightarrow{\mathbb{F}} \mathbb{B}' & & \gamma_{\mathbb{B}}\mathbb{C} \xrightarrow{\mathbb{F}^{\#}} \gamma_{\mathbb{B}'}(\Gamma\mathbb{F})\mathbb{C} \\
 & \gamma_{\mathbb{B}}\uparrow & & & & \gamma_{\mathbb{B}}\uparrow \qquad\qquad \uparrow \\
\gamma_{\mathbb{B}'} & \Gamma\mathbb{B} & \ni & & & \mathbb{C} \qquad\qquad \gamma_{\mathbb{B}'} \\
 & \Gamma\mathbb{F}\downarrow & & & & \Gamma\mathbb{F} \searrow \qquad \\
 & \Gamma\mathbb{B}' & \ni & & & (\Gamma\mathbb{F})\mathbb{C}
\end{array}
$$

**Proposition 2.59** *Let $\mathsf{K} : \mathcal{A} \to \mathcal{A}'$ be a functor preserving finite limits and $\mathbb{F} : \mathbb{B} \to \mathbb{B}'$ be an $\mathcal{A}$-internal functor strictly preserving finite limits. For any $\Gamma\mathbb{B}$-internal category $\mathbb{C}$, the identity $\mathsf{K}\mathbb{F}^{\#}_{\mathbb{C}} = (\mathsf{K}\mathbb{F})^{\#}_{\mathsf{K}^{\#}\mathbb{C}}$ holds as $\mathcal{A}$-internal functors, provided that their domains and codomains are identified along canonical isomorphisms.*

$$
\begin{array}{ccccc}
\begin{array}{cccc}
 & \mathcal{A} & \ni & \gamma_{\mathbb{B}}\mathbb{C} \xrightarrow{\mathbb{F}^{\#}} \gamma_{\mathbb{B}'}(\Gamma\mathbb{F})\mathbb{C} \\
 & \uparrow\gamma_{\mathbb{B}} & & \gamma_{\mathbb{B}}\uparrow \qquad\qquad \uparrow \\
\gamma_{\mathbb{B}'} & \Gamma\mathbb{B} & \ni & \mathbb{C} \qquad\qquad \gamma_{\mathbb{B}'} \\
 & \downarrow\Gamma\mathbb{F} & & \Gamma\mathbb{F}\searrow \qquad \\
 & \Gamma\mathbb{B}' & \ni & (\Gamma\mathbb{F})\mathbb{C}
\end{array}
& \overset{\mathsf{K}}{\dashrightarrow} &
\begin{array}{cccc}
 & \mathcal{A}' & \ni & \gamma_{\mathsf{K}\mathbb{B}}\mathsf{K}^{\#}\mathbb{C} \xrightarrow{(\mathsf{K}\mathbb{F})^{\#}} \gamma_{\mathsf{K}\mathbb{B}}(\Gamma\mathsf{K}\mathbb{F})\mathsf{K}^{\#}\mathbb{C} \\
 & \uparrow\gamma_{\mathsf{K}\mathbb{B}} & & \gamma_{\mathsf{K}\mathbb{B}}\uparrow \qquad\qquad \uparrow \\
\gamma_{\mathsf{K}\mathbb{B}'} & \Gamma\mathsf{K}\mathbb{B} & \ni & \mathsf{K}^{\#}\mathbb{C} \qquad\qquad \gamma_{\mathsf{K}\mathbb{B}'} \\
 & \downarrow\Gamma\mathsf{K}\mathbb{F} & & \Gamma\mathsf{K}\mathbb{F}\searrow \qquad \\
 & \Gamma\mathsf{K}\mathbb{B}' & \ni & (\Gamma\mathsf{K}\mathbb{F})\mathsf{K}^{\#}\mathbb{C}
\end{array}
\end{array}
$$

*Proof.* First, note that the domains and codomains of both sides, namely

$$\mathsf{K}\mathbb{F}^{\#} : \mathsf{K}\gamma_{\mathbb{B}}\mathbb{C} \to \mathsf{K}\gamma_{\mathbb{B}'}(\Gamma\mathbb{F})\mathbb{C}, \qquad (\mathsf{K}\mathbb{F})^{\#} : \gamma_{\mathsf{K}\mathbb{B}}\mathsf{K}^{\#}\mathbb{C} \to \gamma_{\mathsf{K}\mathbb{B}}(\Gamma\mathsf{K}\mathbb{F})\mathsf{K}^{\#}\mathbb{C},$$

are respectively isomorphic. Indeed, by Proposition 2.45 and Proposition 2.51, we have

$$\mathsf{K}\gamma_{\mathbb{B}} \cong \gamma_{\mathsf{K}\mathbb{B}}\mathsf{K}^{\#}, \qquad \mathsf{K}\gamma_{\mathbb{B}}(\Gamma\mathbb{F}) \cong \gamma_{\mathsf{K}\mathbb{B}}\mathsf{K}^{\#}\Gamma\mathbb{F} = \gamma_{\mathsf{K}\mathbb{B}}(\Gamma\mathsf{K}\mathbb{F})\mathsf{K}^{\#}.$$

Under these isomorphisms, the given internal functors coincide on the object of objects by the following chain of identities:

$$\left|(\mathsf{K}\mathbb{F})^{\#}\right| = (\mathsf{K}\mathbb{F})_{1_{\mathsf{K}\mathbb{B}},|\mathsf{K}^{\#}\mathbb{C}|} = \mathsf{K}\big(\mathbb{F}_{1_{\mathbb{B}},|\mathbb{C}|}\big) = \mathsf{K}\left|\mathbb{F}^{\#}\right| = \left|\mathsf{K}\mathbb{F}^{\#}\right|.$$

A similar argument holds for the object of morphisms. $\square$

**Corollary 2.60** *Let $\mathbb{F} : \mathbb{B} \to \mathbb{B}'$ be an $\mathcal{A}$-internal functor strictly preserving finite limits. For any $\Gamma\mathbb{B}$-internal category $\mathbb{C}$, the identity $\Gamma\mathbb{F}^{\#}_{\mathbb{C}} = (\Gamma\mathbb{F})^{\#}_{\mathbb{C}}$ holds under the canonical isomorphisms.*

**Proposition 2.61 ($\gamma \circ \gamma \cong \gamma$)** *Let $\mathbb{B}$ be an $\mathcal{A}$-internal category with chosen finite limits and $\mathbb{C}$ be a $\Gamma\mathbb{B}$-internal category with a chosen terminal object. There is a natural isomorphism $\gamma_{\mathbb{B}} \circ \gamma_{\mathbb{C}} \cong \gamma_{\gamma_{\mathbb{B}}\mathbb{C}}$ under the canonical isomorphism.*

$$
\begin{array}{cccccc}
 & \mathcal{A} & \ni & \mathbb{B} & & \gamma_{\mathbb{B}}\mathbb{C} \\
 & \uparrow\gamma_{\mathbb{B}} & & & & \uparrow\gamma_{\mathbb{B}} \\
\gamma_{\gamma_{\mathbb{B}}\mathbb{C}} \nearrow & \Gamma\mathbb{B} & \ni & & & \mathbb{C} \\
 & \uparrow\gamma_{\mathbb{C}} & & & & \\
\Gamma\gamma_{\mathbb{B}}\mathbb{C} \xrightarrow{\cong} & \Gamma\mathbb{C} & & & &
\end{array}
$$

*Proof.* First, note that their domains $\Gamma\gamma_{\mathbb{B}}\mathbb{C}$ and $\Gamma\mathbb{C}$ are identified via the isomorphism established in Proposition 2.50. In fact, this isomorphism coincides with $(\gamma_{\mathbb{B}})^{\#}_{\mathbb{C}} : \Gamma\mathbb{C} \to \Gamma\gamma_{\mathbb{B}}\mathbb{C}$, as is easily verified.

Consider a global object $X \in_1 \mathbb{C}$. Under the isomorphism $(\gamma_{\mathbb{B}})^{\#}$, we have

$$\gamma_{\gamma_{\mathbb{B}}\mathbb{C}}(X) = (\gamma_{\mathbb{B}}\mathbb{C})\big(1_{\gamma_{\mathbb{B}}\mathbb{C}}, X\big) = (\gamma_{\mathbb{B}}\mathbb{C})\big((\gamma_{\mathbb{B}})^{\#}1_{\mathbb{C}}, (\gamma_{\mathbb{B}})^{\#}X\big) \cong \gamma_{\mathbb{B}}(\mathbb{C}(1_{\mathbb{C}}, X)) = \gamma_{\mathbb{B}}(\gamma_{\mathbb{C}}(X)).$$

This isomorphism is natural in $X$. □

## 2.5 Finite limit theories

Various syntactic formulations correspond to the class of categories with finite limits, such as *essentially algebraic theories* [1], *Cartesian theories* [16], *Cartesian logic* [9], and *partial Horn theories* [35]. All of these admit interpretations within any category with finite limits and the construction of syntactic categories. We do not fix any specific formulation and refer to them collectively as *finite limit theories.* The reader may view them as any of the aforementioned types of theories, or simply as categories with finite limits regarded as syntactic categories.

We summarize basic properties of finite limit theories used in this thesis. These theories subsume many-sorted equational theories. The primary examples of finite limit theories that cannot be expressed equationally are the theory of categories $\mathcal{T}_{\text{cat}}$ and the theory of categories with *chosen* finite limits $\mathcal{T}_{\text{FLcat}}$, which play a crucial role in this thesis. Indeed, many classes of categories can be expressed as finite limit theories, provided that property-like structures are treated as *chosen.* Such examples include the theories of arithmetic universes [27, 29, 31], elementary toposes, and those with natural number objects.

We proceed to the notion of models. Given a finite limit theory $\mathcal{T}$ and a category $\mathcal{B}$ with finite limits, there is a notion of $\mathcal{B}$*-models* of $\mathcal{T}$. Furthermore, we have a notion of $\mathcal{B}$*-model homomorphisms* $f : M \to N$ between $\mathcal{B}$-models. These constitute *the category of* $\mathcal{B}$*-models of* $\mathcal{T}$, which we write as $\mathbf{Mod}(\mathcal{T}, \mathcal{B})$. For instance, $\mathbf{Mod}(\mathcal{T}_{\text{FLcat}}, \mathcal{B})$ consists of $\mathcal{B}$-internal categories with chosen finite limits and $\mathcal{B}$-internal functors strictly preserving them.

Given a functor $F : \mathcal{B} \to \mathcal{B}'$ preserving finite limits, a functor $F_* : \mathbf{Mod}(\mathcal{T}, \mathcal{B}) \to \mathbf{Mod}(\mathcal{T}, \mathcal{B}')$ is induced. For any $M \in \mathbf{Mod}(\mathcal{T}, \mathcal{B})$, the $\mathcal{B}'$-model $F_*(M)$ is simply written as $FM$. When $\mathcal{T} = \mathcal{T}_{\text{cat}}$ or $\mathcal{T}_{\text{FLcat}}$, the functor $F_*$ coincides with the change of base construction discussed in the previous subsection. In particular, the global section functor $\Gamma : \mathcal{B} \to \mathbf{Set}$ induces $\Gamma_* : \mathbf{Mod}(\mathcal{T}, \mathcal{B}) \to \mathbf{Mod}(\mathcal{T}, \mathbf{Set})$, which generalizes the construction of the underlying category.

On the other hand, there is a notion of *theory morphisms* $\varphi : \mathcal{T} \to \mathcal{T}'$ between finite limit theories. A theory morphism $\varphi : \mathcal{T} \to \mathcal{T}'$ induces a functor $\varphi^* : \mathbf{Mod}(\mathcal{T}', \mathcal{B}) \to \mathbf{Mod}(\mathcal{T}, \mathcal{B})$ in the opposite direction between the category of models. For instance, there is a theory morphism $\varphi : \mathcal{T}_{\text{cat}} \to \mathcal{T}_{\text{FLcat}}$, which induces the forgetful functor $\varphi^* : \mathbf{Mod}(\mathcal{T}_{\text{FLcat}}, \mathcal{B}) \to \mathbf{Mod}(\mathcal{T}_{\text{cat}}, \mathcal{B})$. Functors $\varphi^*$ and $F_*$ are compatible in the sense that the following commutes:

$$\begin{array}{ccc} \mathbf{Mod}(\mathcal{T}', \mathcal{B}) & \xrightarrow{F_*} & \mathbf{Mod}(\mathcal{T}', \mathcal{B}') \\ {\scriptstyle \varphi^*}\downarrow & & \downarrow{\scriptstyle \varphi^*} \\ \mathbf{Mod}(\mathcal{T}, \mathcal{B}) & \xrightarrow{F_*} & \mathbf{Mod}(\mathcal{T}, \mathcal{B}'). \end{array}$$

The $F^\#$ construction (Definition 2.44) naturally generalizes to models of any finite limit theory. Given a model $M \in \mathbf{Mod}(\mathcal{T}, \mathcal{B})$ and a functor $F : \mathcal{B} \to \mathcal{B}'$ preserving finite limits, one can construct a morphism $F^\#_M : \Gamma M \to \Gamma F M$ (or simply $F^\#$) in the category $\mathbf{Mod}(\mathcal{T}, \mathbf{Set})$. When $\mathcal{T} = \mathcal{T}_{\mathrm{cat}}$ or $\mathcal{T}_{\mathrm{FLcat}}$, this construction coincides with Definition 2.44. Further, this construction is compatible with functors $\varphi^*$ in the sense that the identity $\varphi^*(F^\#_M) = F^\#_{\varphi^*(M)}$ holds.

A basic fact about finite limit theories is that the functor $\varphi^*$ always has a left adjoint. Among its consequences, we shall only use the following result:

**Theorem 2.62** *For any finite limit theory $\mathcal{T}$, there exists an initial* **Set***-model of $\mathcal{T}$. In other words, the category* $\mathbf{Mod}(\mathcal{T}, \mathbf{Set})$ *has an initial object.*

A detailed proof for partial Horn theories can be found in Palmgren and Vickers [35]. Note that, while the initial **Set**-model always exists, for a general category $\mathcal{B}$ with finite limits, an initial $\mathcal{B}$-model may not exist.

# 3 Code structures and the fixed point theorem

One of the most important lemmas in mathematical logic is the *diagonal lemma.* It guarantees the existence of a certain kind of fixed points for any predicate with respect to a given Gödel numbering system.

In this section, we introduce a structure on a fibration called a *code structure,* which serves as an abstraction of Gödel numbering. It is designed to possess the minimal properties required to prove the diagonal lemma. After developing several constructions of code structures, we establish the *fixed point theorem for fibrations with codes* (Theorem 3.14), which provides a direct generalization of the diagonal lemma from the perspective of fibrational semantics.

While our fixed point theorem is proven in the same way as the classical lemma, it unifies several well-known results proven via diagonal arguments: Cantor's theorem, Kleene's second recursion theorem, the diagonal lemma, and the derivation of Löb's theorem from the diagonal lemma. It may be viewed as a partial extension of Lawvere's fixed point theorem [26] to *intensional* settings.

Another, and perhaps more significant, purpose of introducing code structures is to conceptually reorganize the theory of geminal categories by Ramesh [38]. This will be addressed in the subsequent sections.

## 3.1 Code structures

We start with the definition of a code structure.

**Definition 3.1 (Code structure)** Let $p : \mathcal{E} \to \mathcal{B}$ be a cloven fibration over a category $\mathcal{B}$ with binary products and a chosen terminal object $1$. A *code structure* (or simply *codes*) on the fibration $p$ consists of the following data:

(i) **Internalization of fiber objects:**
- For each object $A \in \mathcal{B}$, an object $\ulcorner \mathcal{E}_A \urcorner \in \mathcal{B}$.
- For each object $A \in \mathcal{B}$ and $X \in \mathcal{E}_A$, a global element $\ulcorner X \urcorner \in_1 \ulcorner \mathcal{E}_A \urcorner$.

(ii) **Internalization of hom-sets:**
- For each pair $A, B \in \mathcal{B}$, an object $\ulcorner \mathcal{B}(A, B) \urcorner \in \mathcal{B}$.
- For each morphism $f : A \to B$ in $\mathcal{B}$, a global element $\ulcorner f \urcorner \in_1 \ulcorner \mathcal{B}(A, B) \urcorner$.

(iii) **Internalization of reindexing on objects:**
- For each pair $A, B \in \mathcal{B}$, a morphism $\mathsf{app}_{A,B} : \ulcorner \mathcal{E}_B \urcorner \times \ulcorner \mathcal{B}(A, B) \urcorner \to \ulcorner \mathcal{E}_A \urcorner$.
- For any $X \in \mathcal{E}_B$ and $f : A \to B$, we require $\mathsf{app}_{A,B}(\ulcorner X \urcorner, \ulcorner f \urcorner) \approx \ulcorner X \cdot f \urcorner$.

(iv) **Internalization of internalization of fiber objects:**
- For each object $A \in \mathcal{B}$, a morphism $\mathsf{code}_A : \ulcorner \mathcal{E}_A \urcorner \to \ulcorner \mathcal{B}(1, \ulcorner \mathcal{E}_A \urcorner) \urcorner$.
- For any $X \in \mathcal{E}_A$, we require $\mathsf{code}_A(\ulcorner X \urcorner) \approx \ulcorner \ulcorner X \urcorner \urcorner$.

A fibration equipped with a code structure is called a *fibration with codes.*

The definition of the weak equality $\approx$ is given in Definition 2.10. See also Remark 3.9 for a discussion on the use of weak equality.

It may be worth clarifying the typing of condition (iv). The object $\ulcorner \mathcal{B}(1, \ulcorner \mathcal{E}_A \urcorner) \urcorner$ is obtained by applying the first part of (ii) to $A = 1$ and $B = \ulcorner \mathcal{E}_A \urcorner$. Correspondingly, the right-hand side $\ulcorner \ulcorner X \urcorner \urcorner$ of the equation is obtained by applying the second part of (ii) to $\ulcorner X \urcorner : 1 \to \ulcorner \mathcal{E}_A \urcorner$. Thus, $\mathsf{code}_A(\ulcorner X \urcorner)$ and $\ulcorner \ulcorner X \urcorner \urcorner$ are both global elements of $\ulcorner \mathcal{B}(1, \ulcorner \mathcal{E}_A \urcorner) \urcorner$.

Before providing examples, we make two further remarks.

**Remark 3.2** All conditions (i)–(iv) follow a common pattern: for certain sets (or classes) and functions between them, there exist corresponding objects and morphisms in $\mathcal{B}$. These data in $\mathcal{B}$ are intended to *internalize* or *encode* those set-theoretic entities. Assuming $\mathcal{B}$ is locally small, these correspondences are summarized as below:

$$\frac{X : \{*\} \longrightarrow |\mathcal{E}_A| \text{ in } \mathbf{Set}}{\ulcorner X \urcorner : 1 \longrightarrow \ulcorner \mathcal{E}_A \urcorner \text{ in } \mathcal{B},} \qquad \frac{f : \{*\} \longrightarrow \mathcal{B}(A, B) \text{ in } \mathbf{Set}}{\ulcorner f \urcorner : 1 \longrightarrow \ulcorner \mathcal{B}(A, B) \urcorner \text{ in } \mathcal{B},}$$

$$\frac{(-) \cdot (-) : |\mathcal{E}_B| \times \mathcal{B}(A, B) \longrightarrow |\mathcal{E}_A| \text{ in } \mathbf{Set}}{\mathsf{app}_{A,B} : \ulcorner \mathcal{E}_B \urcorner \times \ulcorner \mathcal{B}(A, B) \urcorner \longrightarrow \ulcorner \mathcal{E}_A \urcorner \text{ in } \mathcal{B},} \quad \frac{\ulcorner - \urcorner : |\mathcal{E}_A| \longrightarrow \mathcal{B}(1, \ulcorner \mathcal{E}_A \urcorner) \text{ in } \mathbf{Set}}{\mathsf{code}_A : \ulcorner \mathcal{E}_A \urcorner \longrightarrow \ulcorner \mathcal{B}(1, \ulcorner \mathcal{E}_A \urcorner) \urcorner \text{ in } \mathcal{B}.}$$

This observation leads to the construction of code structures described in Section 4.1.

**Remark 3.3** There is a subtlety regarding the choice of terminal objects arising from condition (iv). For two different choices $1$ and $1'$ of terminal objects, there is no guarantee that $\ulcorner \mathcal{B}(1, \ulcorner \mathcal{E}_A \urcorner) \urcorner$ and $\ulcorner \mathcal{B}(1', \ulcorner \mathcal{E}_A \urcorner) \urcorner$ are isomorphic in any sense. Thus, one cannot induce $\mathtt{code}_A$ for $1'$ from a version defined for $1$ in the canonical manner. This is why we define a code structure on a category with a *chosen* terminal object.

To be sure, one could induce a code structure for $1'$ from that for $1$ by using the automorphism on the category swapping $1$ and $1'$ and leaving other objects unchanged. However, this relies on the nature of terminal objects and should not be regarded as a "canonical" way of transporting the structure. This becomes clear considering that, in general, there is no isomorphism $(\mathcal{B}, \times) \to (\mathcal{B}, \times')$ that strictly preserves different choices of binary products.

## 3.2 Codes via Gödel numbering

The standard Gödel numbering provides a primary example of a code structure on a syntactic fibration. We refer to Section 2.2.4 for the formulation of these fibrations. The construction is straightforward and applicable to various logical systems, not limited to the language of arithmetic.

Nonetheless, we start with arithmetic. We assume that a specification $(\Sigma, \Pi)$ contains Primitive Recursive Arithmetic (PRA) in the following sense. First, for each primitive recursive function $f(x_1, ..., x_n)$ (more precisely, each of its *definitions* or *programs*), there is a corresponding term

$$v_1 : N, ..., v_n : N \vdash \underline{f}(v_1, ..., v_n) : N.$$

We assume that the defining equations for such functions are derivable. For example, there are terms $\mathrm{Pred}(x)$ and $x \dot{-} y$ representing the predecessor and the modified subtraction respectively, such that:

$$\begin{aligned} () \mid () &\vdash \mathrm{Pred}(0) =_N 0, & x : N \mid () &\vdash x \dot{-} 0 =_N x, \\ x : N \mid () &\vdash \mathrm{Pred}(S(x)) =_N x, & x : N, y : N \mid () &\vdash x \dot{-} S(y) =_N \mathrm{Pred}(x \dot{-} y). \end{aligned}$$

Furthermore, any equation evaluating a primitive recursive function at specific values is assumed to be derivable. For a natural number $k$, we write $\underline{k}$ as its *numeral* (which may be defined as $S^k(0)$, or treated as the term $\underline{f}$ for the constant function $f = k$). Then, the following judgment is derivable:

$$() \mid () \vdash \underline{f(k_1, ..., k_n)} =_N \underline{f}(\underline{k_1}, ..., \underline{k_n}).$$

In standard signatures such as $(N; 0, S, +, \times)$, these terms $\underline{f}(v_1, ..., v_n)$ may not be directly available. This is resolved by the standard method of *extensions by definitions:* if a primitive recursive function is provably total, its function symbols can be conservatively added to the signature. For instance, any system containing the arithmetic $I\Sigma_1$ can be extended to an equivalent specification $(\Sigma, \Pi)$ that includes PRA in the sense described above [12]. For an equational presentation of PRA, we refer to Goodstein [10].

Traditionally, a Gödel numbering assigns a natural number to each term and formula. In our setting, we naturally extend this assignment to other syntactic objects, such as types, contexts, and tuples of terms (i.e., morphisms in $\mathrm{Cl}(\Sigma)$). We write $\llcorner X \lrcorner$ for the Gödel number of a syntactic object $X$.

**Proposition 3.4** *Let $(\Sigma, \Pi)$ be a specification containing Primitive Recursive Arithmetic in the sense described above. Given a Gödel numbering on the system, one can construct a code structure on the syntactic fibration $p : \mathcal{L}(\Sigma, \Pi) \to \mathrm{Cl}(\Sigma)$.*

*Proof.* Let $p : \mathcal{E} \to \mathcal{B}$ refer to $p : \mathcal{L}(\Sigma, \Pi) \to \mathrm{Cl}(\Sigma)$. We define each component (i)–(iv) of the code structure (Definition 3.1) as follows:

(i) For any $\Gamma \in \mathrm{Cl}(\Sigma)$, we define $\ulcorner \mathcal{E}_\Gamma \urcorner$ to be the context $(v_1 : \mathrm{N}) \in \mathrm{Cl}(\Sigma)$. For any $(\Gamma \mid \Phi) \in \mathcal{E}_\Gamma$, we define $\ulcorner(\Gamma \mid \Phi)\urcorner$ to be the unary list consisting of the term $\underline{\llcorner \Phi \lrcorner} : \mathrm{N}$. For brevity, we simply write $\ulcorner \mathcal{E}_\Gamma \urcorner = \mathrm{N}$ and $\ulcorner(\Gamma \mid \Phi)\urcorner = \underline{\llcorner \Phi \lrcorner}$. We apply similar notation below.

(ii) For any $\Gamma, \Delta \in \mathrm{Cl}(\Sigma)$, we define $\ulcorner \mathcal{B}(\Gamma, \Delta) \urcorner = \mathrm{N}$. For a morphism $(M_1, ..., M_n) : \Gamma \to \Delta$, we define $\ulcorner(M_1, ..., M_n)\urcorner$ to be the term $\underline{\llcorner(M_1, ..., M_n)\lrcorner} : \mathrm{N}$.

(iii) Let $\Gamma$ and $\Delta = (v_1 : \sigma_1, ..., v_n : \sigma_n)$ be type contexts. By the nature of the standard Gödel numbering, there exists a primitive recursive function $f(x, y)$ such that for any proposition context $\Phi$ in $\Delta$ and any morphism $(M_1, ..., M_n) : \Gamma \to \Delta$,

$$f(\llcorner \Phi \lrcorner, \llcorner(M_1, ..., M_n)\lrcorner) = \llcorner \Phi[M_1, ..., M_n / v_1, ..., v_n] \lrcorner.$$

We define $\mathsf{app}_{\Gamma,\Delta}$ to be the term

$$v_1 : \mathrm{N}, v_2 : \mathrm{N} \vdash \underline{f}(v_1, v_2) : \mathrm{N}.$$

For any $\Phi$ and $(M_1, ..., M_n)$, the required equation holds:

$$\begin{aligned}
&\mathsf{app}_{\Gamma,\Delta}(\ulcorner \Phi \urcorner, \ulcorner(M_1, ..., M_n)\urcorner) \\
&= (\underline{f}(v_1, v_2)) \circ \left(\underline{\llcorner \Phi \lrcorner}, \underline{\llcorner(M_1, ..., M_n)\lrcorner}\right) \\
&= \underline{f}(v_1, v_2)\left[\underline{\llcorner \Phi \lrcorner}, \underline{\llcorner(M_1, ..., M_n)\lrcorner} \,/\, v_1, v_2\right] \\
&= \underline{f}\left(\underline{\llcorner \Phi \lrcorner}, \underline{\llcorner(M_1, ..., M_n)\lrcorner}\right) \\
&\approx \underline{f(\llcorner \Phi \lrcorner, \llcorner(M_1, ..., M_n)\lrcorner)} && (*) \\
&= \underline{\llcorner \Phi[M_1, ..., M_n / v_1, ..., v_n] \lrcorner} && \text{by the choice of } f \\
&= \ulcorner(\Gamma \mid \Phi[M_1, ..., M_n / v_1, ..., v_n])\urcorner \\
&= \ulcorner(\Delta \mid \Phi) \cdot (M_1, ..., M_n)\urcorner.
\end{aligned}$$

Here, the weak equality $(*)$ follows from Proposition 2.18 and the derivable judgment

$$() \mid () \vdash \underline{f}\left(\underline{\llcorner \Phi \lrcorner}, \underline{\llcorner(M_1, ..., M_n)\lrcorner}\right) =_{\mathrm{N}} \underline{f(\llcorner \Phi \lrcorner, \llcorner(M_1, ..., M_n)\lrcorner)}.$$

(iv) For any type context $\Gamma$, there exists a primitive recursive function $g(x)$ such that for any proposition context $\Phi$ in $\Gamma$,

$$g(\llcorner \Phi \lrcorner) = \llcorner(\underline{\llcorner \Phi \lrcorner})\lrcorner.$$

Note that $(\underline{\llcorner \Phi \lrcorner})$ is the unary tuple consisting of the term $\underline{\llcorner \Phi \lrcorner} : \mathrm{N}$. We define $\mathsf{code}_\Gamma$ to be the term

$$v_1 : \mathrm{N} \vdash \underline{g}(v_1) : \mathrm{N}.$$

Then, for any $\Phi$, we have:

$$
\begin{aligned}
\mathsf{code}_\Gamma(\ulcorner(\Gamma \mid \Phi)\urcorner) &= \big(\underline{g}(\nu_1)\big) \circ (\underline{\llcorner\Phi\lrcorner}) \\
&= \underline{g}(\nu_1)[\underline{\llcorner\Phi\lrcorner}/\nu_1] \\
&= \underline{g}(\underline{\llcorner\Phi\lrcorner}) \\
&\approx \underline{g(\llcorner\Phi\lrcorner)} \\
&= \underline{\llcorner(\underline{\llcorner\Phi\lrcorner})\lrcorner} && \text{by the choice of } g \\
&= \ulcorner\ulcorner(\Gamma \mid \Phi)\urcorner\urcorner.
\end{aligned}
$$

Again, the weak equality follows from Proposition 2.18. □

This construction applies to various types of formal systems that can be organized into fibrations, not limited to traditional first-order arithmetic.

**Example 3.5 (Internal languages of free categories)** The construction applies to the internal languages of categories freely generated with certain structures including a (parameterized) natural number object, possibly subject to additional axioms. Examples include elementary toposes with natural number objects [24], as well as the arithmetic universes and locoses via their internal dependent type theories [28]. In these cases, a Gödel numbering induces a code structure on the subobject fibration or the splitting of the codomain fibration.

The construction also applies to systems without a basic type N, provided that natural numbers can be embedded into other types. Thus, even set theory can be handled directly.

**Example 3.6 (ZFC)** Consider a specification of ZFC set theory where the signature consists of one basic type $\mathsf{Set}$ and a few function symbols such as $\emptyset, \bigcup$ or $\{x, y\}$, obtained via extensions by definitions. Given a Gödel numbering, we obtain a code structure on the syntactic fibration by setting $\ulcorner\mathcal{E}_A\urcorner = \mathsf{Set}$ and $\ulcorner X\urcorner = \underline{\llcorner X\lrcorner}$, where $\underline{k}$ is the term of type $\mathsf{Set}$ representing the natural number $k$, e.g., via von Neumann ordinals.

Alternatively, if the meta-theory is set-theoretic, one may use the syntactic objects themselves as their own "Gödel codes" rather than encoding them as natural numbers. Such a direct approach is sometimes used to simplify the proofs of incompleteness theorems [36, 39].

Moreover, Gödel codes need not be countable. To illustrate this, we provide an easy example involving a language with uncountably many constants.

**Example 3.7 (Uncountable constants)** Consider the theory of real numbers $(\mathbb{R}; 0, 1, +, \times, <, \mathbb{Z})$ with parameters, i.e., each $r \in \mathbb{R}$ is represented by a distinct constant of type R. One can assign a "Gödel code" in $\mathbb{R}$ to each term and formula by using a suitable function $f : \mathbb{R}^{<\omega} \to \mathbb{R}$ that encodes finite sequences of reals.

Provided that f satisfies certain definability properties, the construction in Proposition 3.4 carries over by setting $\ulcorner\mathcal{E}_A\urcorner = \ulcorner\mathcal{B}(A, B)\urcorner = R$ with the morphisms $\mathsf{app}_{A,B}$ and $\mathtt{code}_A$ defined via the construction of $f$. While $f$ need not be injective in general, an injective encoding is also possible within this language.

As a final example, we consider a code structure derived from the untyped lambda calculus. Although this does not arise from syntactic fibrations in the strict sense of Section 2.2.4, it utilizes the same underlying idea.

We define a fibration arising from untyped lambda calculus as follows. Let $\mathrm{Cl}(\Omega)$ be a category consisting of:

**Objects** Natural numbers, denoted by formal powers $1, \Omega, \Omega^2, \Omega^3, \ldots$.

**Morphisms** $\Omega^m \to \Omega^n$ are $n$-tuples $(M_1, \ldots, M_n)$ of lambda terms where each $M_i$ has free variables only among $\{v_1, \ldots, v_m\}$.

Composition is given by substitution. Next, let $\mathcal{L}(\Omega)$ be a category consisting of:

**Objects** Pairs $(n \mid P)$ of a natural number $n$ and a lambda term $P$ such that free variables in $P$ are among $\{v_1, \ldots, v_n\}$.

**Morphisms** $(m \mid P) \to (n \mid Q)$ are morphisms $(M_1, \ldots, M_n) : \Omega^m \to \Omega^n$ in $\mathrm{Cl}(\Omega)$ such that $P$ and $Q[M_1, \ldots, M_n/v_1, \ldots, v_n]$ are $\beta\eta$-equivalent.

Then, the projection $\mathcal{L}(\Omega) \to \mathrm{Cl}(\Omega), (n \mid P) \mapsto \Omega^n$ forms a fibration.

**Example 3.8 (Untyped lambda calculus)** By fixing a Gödel numbering for untyped lambda terms and using *Church numerals* $\mathbf{k}$ in place of the numerals $\underline{k}$ in arithmetic, one can construct a code structure on the fibration $\mathcal{L}(\Omega) \to \mathrm{Cl}(\Omega)$. The objects $\ulcorner \mathcal{E}_A \urcorner$ and $\ulcorner \mathcal{B}(A, B) \urcorner$ are defined to be $\Omega$. The morphisms $\mathtt{app}_{A,B}$ and $\mathtt{code}_A$ are internalization of the recursive operations on Gödel numbers, represented as lambda terms acting on Church numerals.

Finally, we conclude this subsection with a remark on the use of weak equality.

**Remark 3.9** The reason for using the weak equality $\approx$ rather than strict equality in code structures is to ensure that the construction in Proposition 3.4 is well-defined. While the equation $\underline{f(k_1, \ldots, k_n)} =_N \underline{f}(\underline{k_1}, \ldots, \underline{k_n})$ is provable within the logic, the two sides are not identical as terms, and thus represent different morphisms in $\mathrm{Cl}(\Sigma)$.

One could alternatively regard such equivalence as *external equality* in the term calculus, rather than *internal equality* in logic (see Jacobs [14, p. 177] for a discussion). In this case, the base category would be a quotient of $\mathrm{Cl}(\Sigma)$, and the code structure would satisfy strict equality. For arithmetic, this approach leads to the *initial Skolem theory* — the category freely generated from finite products and a parameterized natural number object — which is categorically well-behaved. Indeed, Trimble [41] adopts this approach to organize the first incompleteness theorem categorically.

However, this quotienting introduces a new complication: since multiple terms correspond to a single morphism, $\ulcorner f \urcorner$ cannot be assigned to a morphism without choosing a representative term. This is conceptually somewhat unnatural since a fundamental characteristic of Gödel coding is its *intensionality* in the sense of Kavvos [18, 19, 21]: it does not preserve any equality derived within the system. Therefore, codes should be assigned to terms themselves, rather than their equivalence classes.

Our capture of intensionality relies on the specific fact that $\mathrm{Cl}(\Sigma)$ forms a category without quotienting terms. This dependency is evident in Example 3.8; while it leads to Kleene's second recursion theorem (Example 3.16), the construction is unavailable for typical programming languages where "composition" might fail associativity. As Kavvos [18, 19] argues, to fully model intensionality, ordinary categories are insufficient and more relaxed structures such as $\mathcal{P}$-categories are needed.

## 3.3 Representable fibrations with codes

In this subsection, we focus on code structures on representable fibrations and provide several examples. First, we present simple conditions sufficient to induce a code structure on any representable fibration.

**Proposition 3.10** *Let $\mathcal{B}$ be a category with binary products and a chosen terminal object. For any $A \in \mathcal{B}$, a code structure is induced on the representable fibration $\mathrm{dom}_A : \mathcal{B}/A \to \mathcal{B}$ provided that $\mathcal{B}$ is equipped with the following data:*[1]

(ii′) **Internalization of hom-sets:**
- For each pair $A, B \in \mathcal{B}$, an object $\ulcorner\mathcal{B}(A, B)\urcorner \in \mathcal{B}$.
- For each morphism $f : A \to B$ in $\mathcal{B}$, a global element $\ulcorner f\urcorner \in_1 \ulcorner\mathcal{B}(A, B)\urcorner$.

(iii′) **Internalization of composition:**
- For each $A, B, C \in \mathcal{B}$, a morphism $\mathsf{comp}_{A,B,C} : \ulcorner\mathcal{B}(A, B)\urcorner \times \ulcorner\mathcal{B}(B, C)\urcorner \to \ulcorner\mathcal{B}(A, C)\urcorner$.
- For any $f : A \to B$ and $g : B \to C$, we require $\mathsf{comp}_{A,B,C}(\ulcorner f\urcorner, \ulcorner g\urcorner) = \ulcorner g \circ f\urcorner$.

(iv′) **Internalization of internalization of hom-sets:**
- For each $A, B \in \mathcal{B}$, a morphism $\mathsf{code}_{A,B} : \ulcorner\mathcal{B}(A, B)\urcorner \to \ulcorner\mathcal{B}(1, \ulcorner\mathcal{B}(A, B)\urcorner)\urcorner$.
- For any $f : A \to B$, we require $\mathsf{code}_{A,B}(\ulcorner f\urcorner) = \ulcorner\ulcorner f\urcorner\urcorner$.

*Proof.* We define each component (i)–(iv) of the code structure (Definition 3.1) on $\mathrm{dom}_A : \mathcal{B}/A \to \mathcal{B}$ as follows:

(i) For $B \in \mathcal{B}$, we define $\ulcorner\mathcal{E}_B\urcorner = \ulcorner\mathcal{B}(B, A)\urcorner$. For each $x \in \mathcal{E}_B = \mathcal{B}(B, A)$, a global element $\ulcorner x\urcorner \in_1 \ulcorner\mathcal{E}_B\urcorner = \ulcorner\mathcal{B}(B, A)\urcorner$ is defined by the second part of (ii′).

(ii) These are inherited from (ii′).

(iii) For $B, C \in \mathcal{B}$, we define the morphism $\mathsf{app}_{B,C}$ by the following composite:

$$\ulcorner\mathcal{B}(C, A)\urcorner \times \ulcorner\mathcal{B}(B, C)\urcorner \xrightarrow{\cong} \ulcorner\mathcal{B}(B, C)\urcorner \times \ulcorner\mathcal{B}(C, A)\urcorner \xrightarrow{\mathsf{comp}_{B,C,A}} \ulcorner\mathcal{B}(B, A)\urcorner.$$

The second part of (iii′) exactly states the required equation for $\mathsf{app}_{B,C}$.

(iv) For $B \in \mathcal{B}$, we define $\mathsf{code}_B = \mathsf{code}_{B,A} : \ulcorner\mathcal{B}(B, A)\urcorner \to \ulcorner\mathcal{B}(1, \ulcorner\mathcal{B}(B, A)\urcorner)\urcorner$. The second part of (iv′) exactly states the required equation for $\mathsf{code}_B$. □

[1] We start the numbering with (ii′) to indicate the correspondence with Definition 3.1.

Note that in representable fibrations, the strict equality $=$ and the weak equality $\approx$ coincide. We now present some examples of code structures that do not arise from Gödel numbering.

**Example 3.11 (Cartesian closed categories)** Let $\mathcal{B}$ be a Cartesian closed category. Code structures on representable fibrations over $\mathcal{B}$ is constructed via Proposition 3.10 as follows:

(ii′) We define $\ulcorner \mathcal{B}(A, B) \urcorner = B^A$. For $f : A \to B$ in $\mathcal{B}$, the global element $\ulcorner f \urcorner \in_1 B^A$ is defined as the Currying of $f \circ \pi_2 : 1 \times A \to B$.

(iii′) For each $A, B, C \in \mathcal{B}$, we define $\mathtt{comp}_{A,B,C} : B^A \times C^B \to C^A$ as the Currying of the following composition:

$$(B^A \times C^B) \times A \xrightarrow{\cong} C^B \times (B^A \times A) \xrightarrow{\mathrm{id} \times \mathrm{ev}} C^B \times B \xrightarrow{\mathrm{ev}} C.$$

The required equation $\mathtt{comp}_{A,B,C}(\ulcorner f \urcorner, \ulcorner g \urcorner) = \ulcorner g \circ f \urcorner$ follows immediately.

(iv′) For each $A, B \in \mathcal{B}$, we define $\mathtt{code}_{A,B} : B^A \to (B^A)^1$ to be the canonical isomorphism, i.e., the Currying of the projection $\pi_1 : B^A \times 1 \to B^A$. The required equation holds because both $\mathtt{code}_{A,B}(\ulcorner f \urcorner)$ and $\ulcorner \ulcorner f \urcorner \urcorner$ are the Currying of $\ulcorner f \urcorner \circ ! : 1 \times 1 \to B^A$.

Notice that the notation $\ulcorner f \urcorner : 1 \to B^A$ for $f : A \to B$ coincides with the convention introduced by Lawvere [26]. For the connection with Lawvere's fixed point theorem, see Example 3.17.

Next, we see another example involving arithmetic, but constructed in a different manner from Section 3.2. Let $T$ be a (possibly intuitionistic) theory of arithmetic. We define the *Lindenbaum–Tarski algebra* $\mathcal{L}_T$ of $T$ as the preordered set consisting of:

**Elements** Closed formulas $\varphi, \psi, \ldots$.

**Order** $\varphi \leqslant \psi$ if $\vdash_T \varphi \to \psi$ (i.e., $\varphi \to \psi$ is provable in T).

The preorder $\mathcal{L}_T$ forms a Heyting algebra (after quotienting by provable equivalence), and in particular, it is Cartesian closed as a category. Hence, the construction in Example 3.11 can be applied to $\mathcal{L}_T$. However, there is an alternative way to define a code structure via Proposition 3.10, that utilizes a provability predicate.

**Example 3.12 (Lindenbaum–Tarski algebra)** Fix a Gödel numbering $\llcorner - \lrcorner$ and a provability predicate $\mathrm{Pr}(x)$ for the theory T. For any closed formula $\varphi$, we write $\Box\varphi$ for $\mathrm{Pr}(\underline{\llcorner \varphi \lrcorner})$. We assume that $\mathrm{Pr}(x)$ satisfies the Hilbert–Bernays derivability conditions; namely, for any closed formulas $\varphi$ and $\psi$,

(D1) $\vdash_T \varphi$ implies $\vdash_T \Box\varphi$.

(D2) $\vdash_T \Box(\varphi \to \psi) \to \Box\varphi \to \Box\psi$.

(D3) $\vdash_T \Box\varphi \to \Box\Box\varphi$.

Then, a code structure is induced via Proposition 3.10 as follows:

(ii′) For $\varphi, \psi \in \mathcal{L}_T$, we define $\ulcorner \mathcal{L}_T(\varphi, \psi) \urcorner = \Box(\varphi \to \psi)$. A global element $\ulcorner f \urcorner$ for a morphism $f : \varphi \leqslant \psi$ is induced since $\vdash_T \varphi \to \psi$ implies $\vdash_T \top \to \Box(\varphi \to \psi)$ by (D1).

(iii′) The existence of a morphism $\mathsf{code}_{\varphi,\psi,\chi}$ requires the provability of the formula $\Box(\varphi \to \psi) \wedge \Box(\psi \to \chi) \to \Box(\varphi \to \chi)$. This is shown by the following calculation:

1. $(\varphi \to \psi) \to (\psi \to \chi) \to (\varphi \to \chi)$ — tautology
2. $\Box((\varphi \to \psi) \to (\psi \to \chi) \to (\varphi \to \chi))$ — 1, (D1)
3. $\Box(\varphi \to \psi) \to \Box((\psi \to \chi) \to (\varphi \to \chi))$ — 2, (D2)
4. $\Box(\varphi \to \psi) \to \Box(\psi \to \chi) \to \Box(\varphi \to \chi)$ — 3, (D2)
5. $\Box(\varphi \to \psi) \wedge \Box(\psi \to \chi) \to \Box(\varphi \to \chi)$ — 4.

(iv′) To establish a morphism $\mathsf{app}_{\varphi,\psi}$, it suffices to show that $\Box(\varphi \to \psi) \to \Box(\top \to \Box(\varphi \to \psi))$ is provable. This follows from (D3).

## 3.4 The fixed point theorem for fibrations with codes

In this subsection, we prove the *fixed point theorem for fibrations with codes*, which serves as a direct generalization of the diagonal lemma to fibrations. The proof is essentially the same as that of the classical diagonal lemma; indeed, the definition of a code structure is specifically designed to be sufficient for this proof.

We begin with a lemma that corresponds to the construction of a primitive recursive function representing the *diagonalization* of formulas, $\varphi(x) \mapsto \varphi(\underline{\llcorner \varphi(x) \lrcorner})$, a technique originally due to Gödel.

**Lemma 3.13 (Diagonalization)** *Let $p : \mathcal{E} \to \mathcal{B}$ be a fibration with codes. Given an object $A \in \mathcal{B}$ and a morphism $i : \ulcorner \mathcal{E}_A \urcorner \to A$, there exists a morphism $d_i : \ulcorner \mathcal{E}_A \urcorner \to \ulcorner \mathcal{E}_1 \urcorner$ such that, for any $X \in \mathcal{E}_A$, the following holds:*

$$d_i(\ulcorner X \urcorner) \approx \ulcorner X \cdot i \cdot \ulcorner X \urcorner \urcorner.$$

Intuitively, $d_i$ internalizes the mapping that sends each "predicate" $X$ over $A$ to its diagonalization:

$$|\mathcal{E}_A| \to |\mathcal{E}_1|, \qquad X \mapsto X \cdot i \cdot \ulcorner X \urcorner.$$

Note that this map is well-defined since $X \cdot i \in \mathcal{E}_{\ulcorner \mathcal{E}_A \urcorner}$ and $\ulcorner X \urcorner : 1 \to \ulcorner \mathcal{E}_A \urcorner$.

*Proof.* We define $d_i$ as the following composite:

$$\begin{aligned}
\ulcorner \mathcal{E}_A \urcorner &\xrightarrow{\langle \mathrm{id}, \ulcorner i \urcorner \circ !, \mathsf{code}_A \rangle} \ulcorner \mathcal{E}_A \urcorner \times \ulcorner \mathcal{B}(\ulcorner \mathcal{E}_A \urcorner, A) \urcorner \times \ulcorner \mathcal{B}(1, \ulcorner \mathcal{E}_A \urcorner) \urcorner \\
&\xrightarrow{\mathsf{app}_{\ulcorner \mathcal{E}_A \urcorner, A} \times \mathrm{id}} \ulcorner \mathcal{E}_{\ulcorner \mathcal{E}_A \urcorner} \urcorner \times \ulcorner \mathcal{B}(1, \ulcorner \mathcal{E}_A \urcorner) \urcorner \\
&\xrightarrow{\mathsf{app}_{1, \ulcorner \mathcal{E}_A \urcorner}} \ulcorner \mathcal{E}_1 \urcorner.
\end{aligned}$$

Then, the required equation is verified straightforwardly. Indeed, for any $X \in \mathcal{E}_A$,

$$
\begin{aligned}
d_i(\ulcorner X \urcorner) &= \mathsf{app}_{1,\ulcorner \mathcal{E}_A \urcorner} \circ \left( \mathsf{app}_{\ulcorner \mathcal{E}_A \urcorner, A} \times \mathrm{id} \right) \circ \langle \mathrm{id}, \ulcorner i \urcorner \circ !, \mathsf{code}_A \rangle \circ \ulcorner X \urcorner \\
&= \mathsf{app}_{1,\ulcorner \mathcal{E}_A \urcorner} \circ \left( \mathsf{app}_{\ulcorner \mathcal{E}_A \urcorner, A} \times \mathrm{id} \right) \circ \langle \ulcorner X \urcorner, \ulcorner i \urcorner, \mathsf{code}_A(\ulcorner X \urcorner) \rangle \\
&\approx \mathsf{app}_{1,\ulcorner \mathcal{E}_A \urcorner} \circ \left( \mathsf{app}_{\ulcorner \mathcal{E}_A \urcorner, A} \times \mathrm{id} \right) \circ \langle \ulcorner X \urcorner, \ulcorner i \urcorner, \ulcorner \ulcorner X \urcorner \urcorner \rangle && (*) \\
&= \mathsf{app}_{1,\ulcorner \mathcal{E}_A \urcorner} \circ \left\langle \mathsf{app}_{\ulcorner \mathcal{E}_A \urcorner, A}(\ulcorner X \urcorner, \ulcorner i \urcorner), \ulcorner \ulcorner X \urcorner \urcorner \right\rangle \\
&\approx \mathsf{app}_{1,\ulcorner \mathcal{E}_A \urcorner}(\ulcorner X \cdot i \urcorner, \ulcorner \ulcorner X \urcorner \urcorner) && (*) \\
&\approx \ulcorner X \cdot i \cdot \ulcorner X \urcorner \urcorner,
\end{aligned}
$$

where we use Proposition 2.11 (2), (3) at the weak equalities $(*)$. □

**Theorem 3.14 (Fixed point theorem for fibrations with codes)** *Let $p : \mathcal{E} \to \mathcal{B}$ be a fibration with codes. Suppose that there exist an object $A \in \mathcal{B}$ and a morphism $i : \ulcorner \mathcal{E}_A \urcorner \to A$ such that the reindexing functor along $i$ is essentially surjective. Then, for any $P \in \mathcal{E}_{\ulcorner \mathcal{E}_1 \urcorner}$, there exists $X \in \mathcal{E}_1$ such that $P \cdot \ulcorner X \urcorner \cong X$ holds in $\mathcal{E}_1$.*

*Proof.* Let $d_i : \ulcorner \mathcal{E}_A \urcorner \to \ulcorner \mathcal{E}_1 \urcorner$ be the morphism constructed in Lemma 3.13. For any given $P \in \mathcal{E}_{\ulcorner \mathcal{E}_1 \urcorner}$, we consider its reindexing along $d_i$, which gives $P \cdot d_i \in \mathcal{E}_{\ulcorner \mathcal{E}_A \urcorner}$.

Since $i^* : \mathcal{E}_A \to \mathcal{E}_{\ulcorner \mathcal{E}_A \urcorner}$ is essentially surjective, there exists an object $Q \in \mathcal{E}_A$ such that $Q \cdot i \cong P \cdot d_i$ in $\mathcal{E}_{\ulcorner \mathcal{E}_A \urcorner}$. We define $X = Q \cdot i \cdot \ulcorner Q \urcorner \in \mathcal{E}_1$. Then, we obtain the following chain of isomorphisms in the fiber $\mathcal{E}_1$:

$$
\begin{aligned}
X &= Q \cdot i \cdot \ulcorner Q \urcorner \\
&\cong P \cdot d_i \cdot \ulcorner Q \urcorner && \text{since } Q \cdot i \cong P \cdot d_i \\
&\cong P \cdot (d_i \circ \ulcorner Q \urcorner) && \text{by Proposition 2.8 (2)} \\
&\cong P \cdot (\ulcorner Q \cdot i \cdot \ulcorner Q \urcorner \urcorner) && \text{by the construction of } d_i \text{ and Proposition 2.11 (1)} \\
&= P \cdot \ulcorner X \urcorner.
\end{aligned}
$$

Thus, $X$ satisfies the required condition $P \cdot \ulcorner X \urcorner \cong X$. □

Note that the object $A$, the morphism $i$, and even the components (ii)–(iv) of a code structure (Definition 3.1) do not appear in the conclusion; it is their mere existence that is required. Additionally, the morphisms $\mathsf{app}$ and $\mathsf{code}$ are used only for the construction of $d_i$ in Lemma 3.13.

We present several applications of the theorem to various fibrations with codes constructed in previous sections. These illustrate the unification nature of the theorem.

**Example 3.15 (Diagonal lemma)** Consider the syntactic fibration with codes $\mathcal{L}(\Sigma, \Pi) \to \mathrm{Cl}(\Sigma)$ from Proposition 3.4. Since $\ulcorner \mathcal{E}_N \urcorner = N$ for the type of natural numbers $N$, the identity morphism $\mathrm{id}_N : \ulcorner \mathcal{E}_N \urcorner \to N$ satisfies the assumption of Theorem 3.14. Unpacking the conclusion of the theorem, we recover the classical diagonal lemma. Specifically, for any formula $\varphi$ with one free variable $v_1 : N$, there exists a closed formula $\psi$ such that $\varphi[\llcorner \psi \lrcorner / v_1] \vdash \psi$ and $\psi \vdash \varphi[\llcorner \psi \lrcorner / v_1]$ are both derivable.

**Example 3.16 (Kleene's second recursion theorem)** Consider the fibration with codes $\mathcal{L}(\Omega) \to \mathrm{Cl}(\Omega)$ from Example 3.8. Again, since $\ulcorner \mathcal{E}_\Omega \urcorner = \Omega$, the identity morphism on $\Omega$ satisfies the assumption of Theorem 3.14. Applying the theorem, we obtain the following: for any lambda term $M$ with one free variable $v_1$, there exists a closed lambda term $N$ with its Gödel number $k$ such that $M[\mathbf{k}/v_1]$ and $N$ are $\beta\eta$-equivalent. This can be viewed as *Kleene's second recursion theorem* formulated in the untyped lambda calculus (albeit in a non-standard manner).

**Example 3.17** Consider the representable fibration $\mathrm{dom}_A : \mathcal{B}/A \to \mathcal{B}$ over a Cartesian closed category $\mathcal{B}$ with codes as defined in Example 3.11. By the fixed point theorem, we obtain: if there exist $B \in \mathcal{B}$ and $i : A^B \to B$ such that $(-) \circ i : \mathcal{B}(B, A) \to \mathcal{B}(A^B, A)$ is surjective, then any $P : A \to A$ has a fixed point $X \in_1 A$ such that $P \circ X = X$.

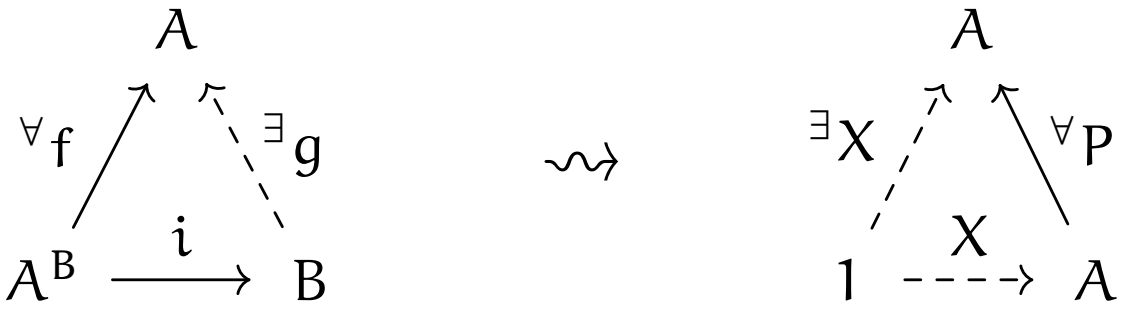

This example can be viewed as a variant of Lawvere's fixed point theorem [26]. While Lawvere's theorem typically assumes a morphism $e : B \to A^B$, ours assumes a morphism $i : A^B \to B$ in the opposite direction. Although there seems to be no immediate implication between their general forms, Lawvere's theorem for a split epimorphism $e$ follows from our result. (See the end of this subsection for a comparison between Lawvere's theorem and Theorem 3.14 itself.)

As a simple application, we can derive Cantor's theorem:

**Example 3.18 (Cantor's theorem)** Consider the subobject fibration $\mathrm{Sub}(\mathbf{Set}) \to \mathbf{Set}$. Since this fibration is isomorphic to the representable fibration $\mathbf{Set}/\{0,1\} \to \mathbf{Set}$, the previous example can be applied. As the conclusion of the theorem — that any function $\{0,1\} \to \{0,1\}$ has a fixed point — is clearly false, the assumption must fail. Thus, for any set $B$ and any map $i : \mathcal{P}(B) \to B$, the induced map $i^{-1} : \mathcal{P}(B) \to \mathcal{P}(\mathcal{P}(B))$ cannot be surjective. This is equivalent to the statement that there is no injective function from $\mathcal{P}(B)$ to $B$.

As a last example, we apply our theorem to recover the derivation of Löb's theorem:

**Example 3.19 (Löb's theorem)** Consider the code structure on the representable fibration $\mathcal{L}_T/\varphi \to \mathcal{L}_T$ from Example 3.12. The assumption of Theorem 3.14 requires a closed formula $\chi$ such that:

- $\vdash_T \Box(\chi \to \varphi) \to \chi$.
- $\vdash_T \Box(\chi \to \varphi) \to \varphi$ implies $\vdash_T \chi \to \varphi$.

In fact, such a formula $\chi$ can be constructed by the classical diagonal lemma so as to satisfy $\vdash_T \Box(\chi \to \varphi) \leftrightarrow \chi$. Thus, the fixed point theorem yields that $\vdash_T \Box\varphi \to \varphi$ implies $\vdash_T \varphi$. This is exactly *Löb's theorem.*

**Remark 3.20** Owing to the fact that $\mathcal{L}_T$ is a preorder, the example above does not constitute a profound result itself. Indeed, what the fixed point theorem performs here is to derive Löb's theorem from the diagonal lemma and the derivability conditions, which is nothing more than an elementary calculation in propositional modal logic. Rather, what this example demonstrates is that the derivation of Löb's theorem from the diagonal lemma can be viewed as a construction of a fixed point via the diagonal method. This perspective has been previously recognized in the context of modal calculi [21, 23].

Interestingly, our formulation reveals that the standard proof of Löb's theorem employs the fixed point theorem *twice:* once in the diagonal lemma (Example 3.15), and once in the derivation of Löb's theorem from it (Example 3.19).

Finally, we compare our fixed point theorem with the well-known Lawvere's fixed point theorem [26]. While both unify classical results based on diagonal arguments, they differ conceptually in their treatment of *intensionality,* which is explicitly discussed by Kavvos [18, 19, 21].

Following Kavvos' terminology, Lawvere's theorem constructs an *extensional* fixed point, whereas ours is *intensional.* Specifically, our fixed point equation $F \cdot \ulcorner X \urcorner \cong X$ involves the coding operation $\ulcorner - \urcorner$. The crucial point here is that this operation need not preserve the very notion of "equality" used in the equation — namely, the isomorphism in the fiber. In contrast, the fixed point equation in Lawvere's theorem is a strict equality between global elements, and any endomorphism must respect that equality. By incorporating an intensional notion of codes, our proof "stratifies" Lawvere's argument, which already reminds us of multi-staged computation in metaprogramming [7, 40].

As we have discussed in Remark 3.9 following Kavvos, our formulation does not yet fully capture the essence of intensionality and more relaxed structures may be needed. However, it should be noted that our formulation provides a categorical treatment of the diagonal lemma that Kavvos deferred to future work [18].

We should also remark that Lawvere's theorem can handle Gödel coding, as seen in his work [26] on Tarski's undefinability of truth and the first incompleteness theorem. However, in his framework, Gödel codes are merely assigned via surjectivity; intensionality is not employed within the fixed point equation itself. While Lawvere's theorem is effective for deriving *negative* results such as Tarski's theorem, it does not seem to address the construction of *positive* intensional fixed points, such as those in the diagonal lemma or Kleene's second recursion theorem.

# 4 Codes via functors

All examples of code structures discussed so far, with the exception of Example 3.11, were constructed using Gödel coding. Due to the very nature of such coding, these structures are highly non-canonical; the operations $\ulcorner - \urcorner$ do not respect the categorical structure of the original fibration, such as isomorphisms.

While this intensional nature seems essential to "coding", there are several ways to induce code structures from data that are more categorically natural. In this section, we present examples of codes induced by simple categorical structures involving functors.

In particular, we focus on constructing codes on representable and codomain fibrations. Since codomain fibrations serve as a "universe" over the base category, codes on them enable the *bootstrapping argument*, which we describe in Section 5. A successful strategy for constructing codes on codomain fibrations leads to the definition of *pre-geminal categories*, a slightly weaker version of Ramesh's *geminal categories*. This approach motivates the definition of geminal categories in the perspective of Gödel coding.

## 4.1 Codes via functors from sets

Recall Remark 3.2: the definition of a code structure requires that there exist objects and morphisms in $\mathcal{B}$ corresponding to certain sets and functions. This observation leads to the following construction:

> **Proposition 4.1** *Let $p : \mathcal{E} \to \mathcal{B}$ be a cloven fibration over a category $\mathcal{B}$ with binary products and a chosen terminal object $1$. Assume that $\mathcal{B}$ is locally small and the set of objects $|\mathcal{E}_A|$ of each fiber is small. Then, any functor $F : \mathbf{Set} \to \mathcal{B}$ preserving finite products induces a code structure on $p$.*

*Proof.* Note that an element of a set $S$ corresponds to a global element $\{*\} \to S$ in the category **Set**. Since $F$ preserves the terminal object, any element $x \in S$ yields the global element $F(x) \in_1 F(S)$. We define the components (i)–(iv) of the code structure as follows:

(i) We define $\ulcorner \mathcal{E}_A \urcorner = F|\mathcal{E}_A|$. For each $X \in |\mathcal{E}_A|$, we define $\ulcorner X \urcorner = F(X) \in_1 F|\mathcal{E}_A|$.

(ii) We define $\ulcorner \mathcal{B}(A, B) \urcorner = F(\mathcal{B}(A, B))$. For each $f : A \to B$, we regard it as $f \in \mathcal{B}(A, B)$ and define $\ulcorner f \urcorner = F(f) \in_1 F(\mathcal{B}(A, B))$.

(iii) For each $A, B \in \mathcal{B}$, the cleavage on $p$ determines the reindexing function $(-) \cdot (-) : |\mathcal{E}_B| \times \mathcal{B}(A, B) \to |\mathcal{E}_A|$, a morphism in **Set**. We define $\mathsf{app}_{A,B}$ by the following composite:

$$F|\mathcal{E}_B| \times F(\mathcal{B}(A, B)) \xrightarrow{\cong} F(|\mathcal{E}_B| \times \mathcal{B}(A, B)) \xrightarrow{F((-) \cdot (-))} F|\mathcal{E}_A|.$$

The required equality for $\mathsf{app}_{A,B}$ follows directly from the functoriality of $F$. (In this case, strict equality holds.)

(iv) From the construction in (i), there is a function $\ulcorner - \urcorner : |\mathcal{E}_A| \to \mathcal{B}(1, F|\mathcal{E}_A|)$ in **Set** for each $A \in \mathcal{B}$. We define $\mathsf{code}_A$ as the morphism $F(\ulcorner - \urcorner) : F|\mathcal{E}_A| \to F(\mathcal{B}(1, F|\mathcal{E}_A|))$. Again, the required equality follows from the functoriality of $F$. $\square$

**Example 4.2** For a small category $\mathcal{C}$, the identity functor $\mathrm{id} : \mathbf{Set} \to \mathbf{Set}$ induces a code structure on the family fibration $\mathrm{Fam}(\mathcal{C}) \to \mathbf{Set}$.

In particular, the proposition may yield code structures on representable fibrations:

**Corollary 4.3** *For a locally small category $\mathcal{B}$ with finite products, any functor $F : \mathbf{Set} \to \mathcal{B}$ preserving finite products induces a code structure on any representable fibration over $\mathcal{B}$.*

What about the codomain fibration? To apply Proposition 4.1 to the codomain fibration, each $|\mathcal{B}/A|$ must be small, which implies that $\mathcal{B}$ itself must be a small category. Thus, we would need a functor $F : \mathbf{Set} \to \mathcal{B}$ from the large category $\mathbf{Set}$ to the small category $\mathcal{B}$. Unfortunately, such functors are very rare.

**Lemma 4.4** *Let $F : \mathbf{Set} \to \mathcal{B}$ be a functor from* **Set** *to a* small *category $\mathcal{B}$. Then, for any nonempty set $X$, the morphism $F(!_X) : F(X) \to F(\{*\})$ is an isomorphism.*

*Proof.* We write $\#X$ for the cardinality of a set $X$. Since the set of morphisms $|\mathcal{B}^{\to}|$ is small, there exists a set $I$ such that $\#I > \#|\mathcal{B}^{\to}|$. Consider the set of functions $X^I$ and let $\pi_i : X^I \to X$ be the projection at $i \in I$. Applying the functor $F$, we obtain a family of morphisms $\big(F(\pi_i) : F(X^I) \to F(X)\big)_{i \in I}$ in $\mathcal{B}$. Since $\#\mathcal{B}\big(F(X^I), F(X)\big) \leqslant \#|\mathcal{B}^{\to}| < \#I$, by the pigeonhole principle, there must exist distinct indices $i_0, i_1 \in I$ such that $F(\pi_{i_0}) = F(\pi_{i_1})$.

Since $X$ is nonempty, we can fix an element $x_0 \in X$. Let $s : X \to X^I$ be the function defined as follows:

$$s(x)(i) = \begin{cases} x_0 & \text{if } i = i_0, \\ x & \text{if } i \neq i_0. \end{cases}$$

Then, we have a commutative diagram in **Set**,

$$\begin{array}{ccccc} \{*\} & \xleftarrow{!_X} & X & & \\ {\scriptstyle x_0}\downarrow & & {\scriptstyle s}\downarrow & \searrow^{\mathrm{id}_X} & \\ X & \xleftarrow{\pi_{i_0}} & X^I & \xrightarrow{\pi_{i_1}} & X. \end{array}$$

Applying $F$ to this diagram and using $F(\pi_{i_0}) = F(\pi_{i_1})$, we obtain $F(x_0) \circ F(!_X) = \mathrm{id}_{F(X)}$. On the other hand, $F(!_X) \circ F(x_0) = F(\mathrm{id}_{\{*\}}) = \mathrm{id}_{F(\{*\})}$ holds. Thus, $F(!_X)$ has an inverse $F(x_0)$. □

**Proposition 4.5** *Let $\mathcal{B}$ be a small category. Consider the inclusion functor*

$$j : \mathbf{2} \to \mathbf{Set}, \qquad (\cdot \to \cdot) \mapsto \big(\emptyset \xrightarrow{!} \{*\}\big).$$

*Then, the precomposition functor $j^* = (-) \circ j : [\mathbf{Set}, \mathcal{B}] \to [\mathbf{2}, \mathcal{B}] \cong \mathcal{B}^{\to}$ gives an equivalence between the functor categories $[\mathbf{Set}, \mathcal{B}]$ and $\mathcal{B}^{\to}$. In other words, any functor $F : \mathbf{Set} \to \mathcal{B}$ is determined, up to natural isomorphism, only by its values on the morphism $! : \emptyset \to \{*\}$.*

*Proof.* We construct a (weak) inverse $\Phi : \mathcal{B}^{\to} \to [\mathbf{Set}, \mathcal{B}]$ of the functor $j^*$. For an object $f : A \to B$ in $\mathcal{B}^{\to}$, we define $\Phi(f) : \mathbf{Set} \to \mathcal{B}$ as follows:

$$\Phi(f)(X) = \begin{cases} A & \text{if } X = \emptyset, \\ B & \text{if } X \neq \emptyset, \end{cases}$$

$$\Phi(f)\left(X \xrightarrow{s} Y\right) = \begin{cases} \mathrm{id}_A & \text{if } X = Y = \emptyset, \\ f & \text{if } X = \emptyset, Y \neq \emptyset, \\ \mathrm{id}_B & \text{if } X \neq \emptyset, Y \neq \emptyset. \end{cases}$$

The action of $\Phi$ on morphisms in $\mathcal{B}^{\rightarrow}$ is defined naturally.

By definition, $j^* \circ \Phi \cong \mathrm{id}$ holds. For the other direction, we define a natural transformation $\sigma : \mathrm{id} \Rightarrow \Phi \circ j^*$ such that for any $F : \mathbf{Set} \to \mathcal{B}$ and $X \in \mathbf{Set}$,

$$\sigma_{F,X} : F(X) \to \Phi(F \circ j)(X),$$

$$\sigma_{F,X} = \begin{cases} \mathrm{id}_{F(\emptyset)} & \text{if } X = \emptyset, \\ F(!_X) & \text{if } X \neq \emptyset. \end{cases}$$

When $\mathcal{B}$ is small, then $\sigma : \mathrm{id} \Rightarrow \Phi \circ j^*$ is a natural isomorphism by Lemma 4.4. Therefore, $j^*$ is an equivalence. □

**Corollary 4.6** *For a small category $\mathcal{B}$ with finite products, there is a one-to-one correspondence between the following, up to natural isomorphism:*

- *Functors $F : \mathbf{Set} \to \mathcal{B}$ preserving finite products.*
- *Subterminals in $\mathcal{B}$ (i.e., subobjects of the terminal object).*

*Proof.* By Proposition 4.5, any functor $F : \mathbf{Set} \to \mathcal{B}$ is uniquely determined by its values on $\emptyset, \{*\}$ and $! : \emptyset \to \{*\}$. When we require $F$ to preserve finite products, the values at $\{*\}$ and $!$ are uniquely determined from the other; thus, the functor depends only on the value at $\emptyset$. The condition imposed on $A = F(\emptyset)$ is that the two diagrams on the right-hand side of the following form product cones:

$$\emptyset \xleftarrow{\mathrm{id}} \emptyset \xrightarrow{\mathrm{id}} \emptyset \quad \stackrel{F}{\longmapsto} \quad A \xleftarrow{\mathrm{id}} A \xrightarrow{\mathrm{id}} A,$$

$$\underset{(X \neq \emptyset)}{\emptyset \xleftarrow{\mathrm{id}} \emptyset \xrightarrow{\exists!} X} \quad \stackrel{F}{\longmapsto} \quad A \xleftarrow{\mathrm{id}} A \xrightarrow{!} 1.$$

The second diagram always forms a product. The first diagram forms a product if and only if $!_A : A \to 1$ is a monic, namely, $A$ corresponds to a subterminal. □

Unfortunately, applying the fixed point theorem (Theorem 3.14) to the code structures arising from these functors $F : \mathbf{Set} \to \mathcal{B}$ results in only trivial facts. Hence, we conclude that the construction in Proposition 4.1 fails to produce any useful codes on codomain fibrations. In the subsequent subsections, we present a different, successful strategy to construct codes on codomain fibrations.

**Remark 4.7** While the obstacle here is the size problem, the construction of code structures in the proof of Proposition 4.1 does not necessarily require functors from all sets. For instance, to induce codes on the codomain fibration over $\mathcal{B}$, it suffices to have a functor $\mathbf{Set}_{\leqslant\kappa} \to \mathcal{B}$ from the category of sets with cardinality at most $\kappa$, where $\kappa = \#|\mathcal{B}^{\to}|$. We do not pursue the possibility of such constructions here; indeed, this might be a special case of the following remark with $\mathcal{A} = \mathbf{Set}_{\leqslant\kappa}$.

**Remark 4.8** Despite the failure discussed above, it is worth noting that the proofs of Lemma 4.4 and Proposition 4.5 heavily rely on the non-constructive principle that every non-empty set is inhabited. This suggests that non-trivial codes on codomain fibrations might be constructible via an analog of Proposition 4.1, if one works internally to a category $\mathcal{A}$ where such classical principles fail.

In fact, this idea relates to the concept of *introspective theories* by Ramesh [38]. When working internally to $\mathcal{A}$, the role of **Set** is played by the codomain fibration over $\mathcal{A}$; see Bénabou [5]. Consequently, a functor $\mathbf{Set} \to \mathcal{B}$ corresponds to a fibered functor from $\mathrm{cod}_{\mathcal{A}}$ to a small fibration over $\mathcal{A}$. If this functor preserves finite limits, we recover the definition of an introspective theory. Indeed, Ramesh's original proof of *Löb's theorem for introspective theories* [38, Theorem 4.19] is virtually identical to applying the bootstrapping argument (Section 5) to the "$\mathcal{A}$-internal code structures" induced by introspective theories. A detailed treatment of introspective theories is, however, beyond the scope of this thesis.

## 4.2 Codes on representable fibrations via enriched categories

Recalling Proposition 3.10, condition (ii′) naively states that there is an internalization of hom-sets of $\mathcal{B}$ within $\mathcal{B}$ itself. Thus, it is natural to consider $\mathcal{B}$-enriched categories, whose hom-objects are within $\mathcal{B}$. Indeed, given a $\mathcal{B}$-enriched category $\mathcal{C}$ and a mapping of objects from $\mathcal{B}$ to $\mathcal{C}$, one can obtain such an internalization by mapping a pair of objects of $\mathcal{B}$ into $\mathcal{C}$ and taking their hom-object. This idea provides another method for equipping representable fibrations with codes. It also serves as a first step toward the concept of Ramesh's geminal categories. We refer to Section 2.3 for the notation.

**Lemma 4.9** *Let $\mathcal{B}$ be a category with finite products and $\mathcal{C}$ be a $\mathcal{B}$-enriched category. Then, any functor $F : \mathcal{B} \to \Gamma_*\mathcal{C}$ from $\mathcal{B}$ to the underlying category of $\mathcal{C}$ induces the components* (i)–(iii) *of a code structure (Definition 3.1) on any representable fibration over $\mathcal{B}$.*

*Proof.* It is enough to construct components (ii′) and (iii′) in Proposition 3.10.

(ii′) We define $\ulcorner\mathcal{B}(A, B)\urcorner = \mathcal{C}(FA, FB) \in \mathcal{B}$. For any morphism $f : A \to B$, $\ulcorner f\urcorner$ is defined as the global element $F(f) \in_1 \mathcal{C}(FA, FB)$, which is a morphism in $\Gamma_*\mathcal{C}$.

(iii′) We define $\mathtt{comp}_{A,B,C}$ to be the composition morphism in $\mathcal{C}$,

$$\circ_{FA,FB,FC} : \mathcal{C}(FA, FB) \times \mathcal{C}(FB, FC) \to \mathcal{C}(FA, FC).$$

The required equality $\mathtt{comp}_{A,B,C}(\ulcorner f\urcorner, \ulcorner g\urcorner) = \ulcorner g \circ f\urcorner$ reduces to $Fg \circ Ff = F(g \circ f)$, which follows from the functoriality of $F$. □

However, these data are not enough to yield the component (iv′) in Proposition 3.10 (and thus (iv) in Definition 3.1). To complete the structure, we must construct a morphism

$$\mathsf{code}_{A,B} : \mathcal{C}(FA, FB) \to \mathcal{C}(F1, F(\mathcal{C}(FA, FB)))$$

for each $A, B \in \mathcal{B}$. Since this is a morphism from a hom-object of $\mathcal{C}$, it is natural to further assume the existence of a $\mathcal{B}$-enriched functor from $\mathcal{C}$. To induce such a morphism from a $\mathcal{B}$-enriched functor $H : \mathcal{C} \to \mathcal{D}$, the $\mathcal{B}$-enriched category $\mathcal{D}$ must satisfy

$$\mathcal{D}(HFA, HFB) \cong \mathcal{C}(F1, F(\mathcal{C}(FA, FB))).$$

Such an enriched category $\mathcal{D}$ can be constructed by taking a change of base of $\mathcal{C}$.

**Definition 4.10 (The endofunctor $\Box$)** Let $\mathcal{B}$ be a category with finite products, $\mathcal{C}$ be a $\mathcal{B}$-enriched category with finite products, and $F : \mathcal{B} \to \Gamma_*\mathcal{C}$ be a functor preserving finite products. We define the endofunctor $\Box_{\mathcal{B}} : \mathcal{B} \to \mathcal{B}$ (or simply $\Box$) by the composite $\gamma_{\mathcal{C}} \circ F$, where $\gamma_{\mathcal{C}}$ is the global section functor of $\mathcal{C}$ (see Definition 2.30).

$$\mathcal{B} \xrightarrow{F} \Gamma_*\mathcal{C} \xrightarrow{\gamma_{\mathcal{C}}} \mathcal{B} \quad (\text{composite } \Box : \mathcal{B} \to \mathcal{B})$$

**Remark 4.11** By Proposition 2.28, the endofunctor $\Box$ preserves finite products. Consequently, when $\mathcal{B}$ is Cartesian closed, the pair $(\mathcal{B}, \Box)$ forms a categorical model for the modal logic K. Such a pair is called a *Kripke category* [20, 34]. The relationship between the endofunctor $\Box$ and modality will be further discussed in Remark 4.18.

Considering the change of base of $\mathcal{C}$ along $\Box$, we observe that

$$(\Box_*\mathcal{C})(X, Y) = \Box(\mathcal{C}(X, Y)) = \mathcal{C}(1_{\mathcal{C}}, F(\mathcal{C}(X, Y))) \cong \mathcal{C}(F1, F(\mathcal{C}(X, Y))),$$

which is close to the condition for $\mathcal{D}$ discussed above. Indeed, if we assume the existence of a $\mathcal{B}$-enriched functor $H : \mathcal{C} \to \Box_*\mathcal{C}$ satisfying a certain condition that corresponds to the required equation in Proposition 3.10 (iv′), we can successfully equip any representable fibration with a code structure.

**Theorem 4.12** *Given the following data, one can construct a code structure on any representable fibration over $\mathcal{B}$:*

(I′) A category $\mathcal{B}$ with finite products.

(II′) A $\mathcal{B}$-enriched category $\mathcal{C}$ with finite products.

(III′) A functor $F : \mathcal{B} \to \Gamma_*\mathcal{C}$ preserving finite products.

(IV′) A $\mathcal{B}$-enriched functor $H : \mathcal{C} \to \Box_*\mathcal{C}$ satisfying $\Gamma_*H \circ F = F^{\#e}_{\mathcal{C}} \circ F$.

$$\mathcal{B} \xrightarrow{F} \Gamma_*\mathcal{C} \underset{\Gamma_*H}{\overset{F^{\#e}}{\rightrightarrows}} \Gamma_*\Box_*\mathcal{C}$$

*Proof.* First, observe that $\Gamma_* \Box_* \mathcal{C} = \Gamma_* \gamma_* F_* \mathcal{C} = \Gamma_* F_* \mathcal{C}$ by Proposition 2.31, and hence the functor $F^{\#e} : \Gamma_* \mathcal{C} \to \Gamma_* \Box_* \mathcal{C}$ is induced according to Definition 2.26.

By Lemma 4.9, it suffices to construct the component (iv′) in Proposition 3.10. For each $A, B \in \mathcal{B}$, the $\mathcal{B}$-enriched functor $H$ provides a morphism between hom-objects,

$$H_{FA,FB} : \mathcal{C}(FA, FB) \to (\Box_* \mathcal{C})(HFA, HFB).$$

By $\Gamma_* H \circ F = F^{\#e} \circ F$, we have

$$HFA = (\Gamma_* H)FA = F^{\#e} FA = FA,$$

and similarly for $B$. Thus,

$$(\Box_* \mathcal{C})(HFA, HFB) = (\Box_* \mathcal{C})(FA, FB) = \Box(\mathcal{C}(FA, FB)) \cong \mathcal{C}(F1, F(\mathcal{C}(FA, FB))).$$

Hence, we may define $\mathbf{code}_{A,B}$ as the following composite:

$$\mathcal{C}(FA, FB) \xrightarrow{H_{FA,FB}} (\Box_* \mathcal{C})(HFA, HFB) \xrightarrow{\cong} \mathcal{C}(F1, F(\mathcal{C}(FA, FB))).$$

Next, we verify the equation $\mathbf{code}_{A,B}(\ulcorner f \urcorner) = \ulcorner\ulcorner f \urcorner\urcorner$ for $f : A \to B$. By assumption, we have $(\Gamma_* H)Ff = F^{\#e} Ff$, which shows that the following commutes:

$$\begin{array}{ccccc} 1 & \xrightarrow{Ff} & \mathcal{C}(FA, FB) & \xrightarrow{H_{FA,FB}} & (\Box_* \mathcal{C})(HFA, HFB) \\ & \searrow^{F(Ff)} & & & \downarrow \cong \\ & & & & \mathcal{C}(F1, F(\mathcal{C}(FA, FB))). \end{array}$$

This precisely means $\mathbf{code}_{A,B} \circ \ulcorner f \urcorner = \ulcorner\ulcorner f \urcorner\urcorner$, as required. □

**Remark 4.13** The structure consisting of (I′)–(IV′) is an enriched analog of what we call a pre-geminal category (Definition 4.15). In a similar way, an enriched analog of a geminal category (Definition 6.5) could be considered. Such data induce a triple $(\mathcal{B}, \Box, \delta)$ as in Theorem 6.16, a natural generalization of Kavvos' *Kripke-4 category* [20] to a setting that requires only finite products rather than Cartesian closedness. We leave a development of this enriched theory for future work. We thank Ramesh for suggesting this connection.

## 4.3 Codes on codomain fibrations via internal categories

While the previous construction using enriched categories succeeds in yielding code structures on representable fibrations, it is insufficient for constructing those on codomain fibrations. To achieve the latter, we need the internalization of not only hom-sets but also the collection of objects of slice categories, $|\mathcal{B}/A|$. For this purpose, it is appropriate to use the notion of internal categories rather than enriched categories. Fortunately, the arguments from the previous subsection can be extended to the internal setting in almost the same manner.

**Definition 4.14 (The endofunctor $\Box$)** Let $\mathcal{B}$ be a category with finite limits, $\mathbb{C}$ be a $\mathcal{B}$-internal category with chosen finite limits, and $F : \mathcal{B} \to \Gamma\mathbb{C}$ be a functor preserving finite limits. We define the endofunctor $\Box_{\mathcal{B}} : \mathcal{B} \to \mathcal{B}$ (or simply $\Box$) by the composite $\gamma_{\mathbb{C}} \circ F$, where $\gamma_{\mathbb{C}}$ is the global section functor of $\mathbb{C}$ (see Definition 2.49).

$$\mathcal{B} \xrightarrow{F} \Gamma\mathbb{C} \xrightarrow{\gamma_{\mathbb{C}}} \mathcal{B}, \quad \Box = \gamma_{\mathbb{C}} \circ F$$

By Proposition 2.48 and Proposition 2.56, the functor $\Box$ preserves finite limits. In particular, one may consider the change of base of $\mathbb{C}$ along $\Box$, which leads to an internal analog of the structure presented in Theorem 4.12. We call this structure a *pre-geminal category,*[2] as it is a slightly weaker version of the geminal category introduced by Ramesh [38].

**Definition 4.15 (Pre-geminal category)** A *pre-geminal category* consists of the following data:

(I) A category $\mathcal{B}$ with chosen finite limits.
(II) A $\mathcal{B}$-internal category $\mathbb{C}$ with chosen finite limits.
(III) A functor $F : \mathcal{B} \to \Gamma\mathbb{C}$ strictly preserving finite limits.
(IV) A $\mathcal{B}$-internal functor $\mathbb{H} : \mathbb{C} \to \Box\mathbb{C}$ satisfying $\Gamma\mathbb{H} \circ F = F^{\#}_{\mathbb{C}} \circ F$.

$$\mathcal{B} \xrightarrow{F} \Gamma\mathbb{C} \underset{\Gamma\mathbb{H}}{\overset{F^{\#}}{\rightrightarrows}} \Gamma\Box\mathbb{C}$$

In condition (IV), we identify the codomains of $\Gamma\mathbb{H} : \Gamma\mathbb{C} \to \Gamma\Box\mathbb{C}$ and $F^{\#} : \Gamma\mathbb{C} \to \Gamma F\mathbb{C}$ along the isomorphism $\Gamma\Box\mathbb{C} = \Gamma\gamma F\mathbb{C} \cong \Gamma F\mathbb{C}$, which is obtained by Proposition 2.50.

**Lemma 4.16** *Any pre-geminal category* $(\mathcal{B}, \mathbb{C}, F, \mathbb{H})$ *induces a quadruple* $(\mathcal{B}, \mathcal{C}', F', H')$ *that satisfies the requirements for components* (I′)–(IV′) *in Theorem 4.12.*

*Proof.* We define each component of the quadruple $(\mathcal{B}, \mathcal{C}', F', H')$ as follows:

(I′) The first component is $\mathcal{B}$.

(II′) The $\mathcal{B}$-enriched category $\mathcal{C}'$ is defined as follows:

**Objects** Objects of $\mathcal{B}$, i.e., we set $|\mathcal{C}'| = |\mathcal{B}|$.
**Hom-objects** For $A, B \in \mathcal{B}$, we define the hom-object $\mathcal{C}'(A, B)$ as $\mathbb{C}(FA, FB) \in \mathcal{B}$.
**Identities and composition** These are inherited from $\mathbb{C}$.

It is straightforward to check that $\mathcal{C}'$ has finite products inherited from those in $\mathcal{B}$.

(III′) The functor $F' : \mathcal{B} \to \Gamma_*\mathcal{C}'$ is defined as follows:

**On objects** For an object $A \in \mathcal{B}$, we define $F'A = A$.

[2] This use of the prefix *pre-* is not related to the notion of *pre-introspective theories* in Ramesh [38]. Perhaps a more appropriate name for what we call a pre-geminal category would be an *almost geminal category.*

**On morphisms** For a morphism $f : A \to B$ in $\mathcal{B}$, we define $F'f = Ff \in_1 \mathbb{C}(FA, FB)$.

It is straightforward to check that $F'$ preserves finite products.

(IV′) First, note that the endofunctors $\Box : \mathcal{B} \to \mathcal{B}$ defined via $(\mathcal{B}, \mathbb{C}, F)$ and $(\mathcal{B}, \mathcal{C}', F')$ coincide: indeed, for any object $A \in \mathcal{B}$, we have $\mathcal{C}'(1_{\mathcal{C}'}, F'A) = \mathbb{C}(1_{\mathbb{C}}, FA)$, and similarly for morphisms. We define a $\mathcal{B}$-enriched functor $H' : \mathcal{C}' \to \Box_* \mathcal{C}'$ as follows:

**On objects** For any object $A \in \mathcal{C}'$ (i.e., $A \in \mathcal{B}$), we set $H'A = A$.

**On hom-objects** Let $A, B \in \mathcal{B}$. The internal functor $\mathbb{H}$ provides a morphism in $\mathcal{B}$,

$$\mathbb{H}_{FA,FB} : \mathbb{C}(FA, FB) \to (\Box\mathbb{C})(\mathbb{H}FA, \mathbb{H}FB).$$

The global object $\mathbb{H}FA \in_1 \Box\mathbb{C}$ is obtained by applying the functor $(\Gamma\mathbb{H}) \circ F : \mathcal{B} \to \Gamma\Box\mathbb{C}$ to $A \in \mathcal{B}$. Since $(\Gamma\mathbb{H}) \circ F = F^\# \circ F$, this global object $\mathbb{H}FA \in_1 \Box\mathbb{C}$ coincides with $F^\# FA \in_1 F\mathbb{C}$ under the canonical isomorphism $\Gamma\Box\mathbb{C} \cong \Gamma F\mathbb{C}$. As discussed in Proposition 2.61, this isomorphism is provided by the functor

$$(\gamma_{\mathbb{C}})^\#_{F\mathbb{C}} : \Gamma F\mathbb{C} \to \Gamma\gamma_{\mathbb{C}} F\mathbb{C} = \Gamma\Box\mathbb{C}.$$

Thus, we have the identity

$$\mathbb{H}FA = (\gamma_{\mathbb{C}})^\# F^\# FA = \Box^\# FA \in_1 \Box\mathbb{C}.$$

Based on this, we define $H'_{A,B}$ as the following composite:

$$\mathbb{C}(FA, FB) \xrightarrow{\mathbb{H}_{FA,FB}} (\Box\mathbb{C})(\mathbb{H}FA, \mathbb{H}FB) = (\Box\mathbb{C})\big(\Box^\# FA, \Box^\# FB\big) \xrightarrow{\cong} \Box(\mathbb{C}(FA, FB)).$$

Finally, we verify the required equation $\Gamma_* H' \circ F' = (F')^{\#e} \circ F'$. On objects, we have $(\Gamma_* H')F'A = (F')^{\#e} F'A = A$ for any $A \in \mathcal{B}$. On morphisms, for any $f : A \to B$ in $\mathcal{B}$, the required equation $(\Gamma_* H')F'f = (F')^{\#e} F'f$ follows directly from the assumption $(\Gamma\mathbb{H})Ff = F^\# Ff$ in the definition of a pre-geminal category. □

We proceed to the construction of codes on the codomain fibration.

> **Theorem 4.17** *Any pre-geminal category $(\mathcal{B}, \mathbb{C}, F, \mathbb{H})$ induces code structures on any representable fibration and the codomain fibration over $\mathcal{B}$.*

*Proof.* Codes on representable fibrations are induced by Theorem 4.12 and Lemma 4.16. We define each component (i)–(iv) of the code structure (Definition 3.1) on the codomain fibration over $\mathcal{B}$.

(i) We define $\ulcorner \mathcal{E}_A \urcorner = |\mathbb{C}/FA|$. For any object $x : B \to A$ in the fiber $\mathcal{B}/A$, we have a global element $Fx \in_1 \mathbb{C}(FB, FA)$. We define $\ulcorner x \urcorner = \iota'(Fx) \in_1 |\mathbb{C}/FA|$, where $\iota' : \mathbb{C}(FB, FA) \to |\mathbb{C}/FA|$ is the inclusion morphism.

(ii) We define $\ulcorner \mathcal{B}(A, B) \urcorner = \mathbb{C}(FA, FB)$. For $f : A \to B$, we define $\ulcorner f \urcorner = Ff \in_1 \mathbb{C}(FA, FB)$.

(iii) For $A, B \in \mathcal{B}$, we define $\mathsf{app}_{A,B}$ to be the morphism

$$\mathrm{pb}_{FA,FB} : |\mathbb{C}/FB| \times \mathbb{C}(FA, FB) \to |\mathbb{C}/FA|,$$

which is induced from pullbacks in $\mathbb{C}$ (see Definition 2.54). Since $F : \mathcal{B} \to \Gamma\mathbb{C}$ is assumed to *strictly* preserve pullbacks, for any $x \in \mathcal{B}/B$ and $f : A \to B$, we have

$$F(x \cdot f) = \mathrm{pb}_{FA,FB}(F(x), F(f)) \in_1 |\mathbb{C}/FA|.$$

This shows the required equation for $\mathsf{app}_{A,B}$.

(iv) First, note that

$$\ulcorner\mathcal{B}(1, \ulcorner\mathcal{E}_A\urcorner)\urcorner = \mathbb{C}(F1, F|\mathbb{C}/FA|) = \mathbb{C}(1_{\mathbb{C}}, F|\mathbb{C}/FA|) = \Box|\mathbb{C}/FA| \cong \big|(\Box\mathbb{C})/\big(\Box^{\#}FA\big)\big|.$$

On the other hand, the internal functor $\mathbb{H} : \mathbb{C} \to \Box\mathbb{C}$ induces a morphism $\mathbb{H}_A : |\mathbb{C}/FA| \to |(\Box\mathbb{C})/(\mathbb{H}FA)|$, which is defined as the unique morphism making the following commute:

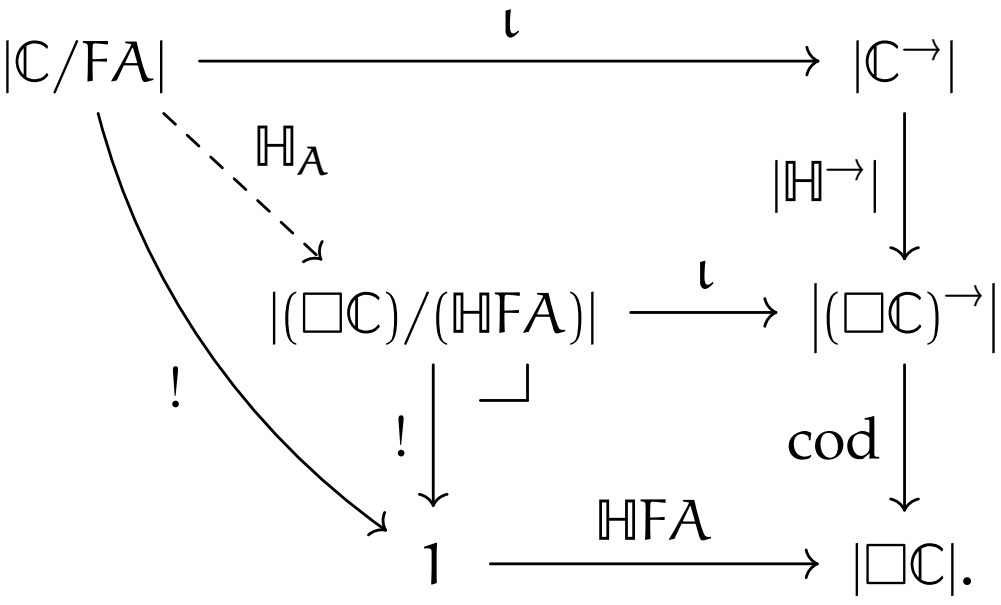

As observed in Lemma 4.16, we have $\Box^{\#}FA = \mathbb{H}FA \in_1 \Box\mathbb{C}$. Based on this, we define $\mathsf{code}_A$ by the following composite:

$$|\mathbb{C}/FA| \xrightarrow{\mathbb{H}_A} |(\Box\mathbb{C})/(\mathbb{H}FA)| = \big|(\Box\mathbb{C})/\big(\Box^{\#}FA\big)\big| \xrightarrow{\cong} \Box|\mathbb{C}/FA|.$$

Finally, we verify the required equation $\mathsf{code}_A(\ulcorner x\urcorner) = \ulcorner\ulcorner x\urcorner\urcorner$. Let $x : B \to A$ be an object of the fiber $\mathcal{B}/A$. Then, the equation follows from the following commutative diagram:

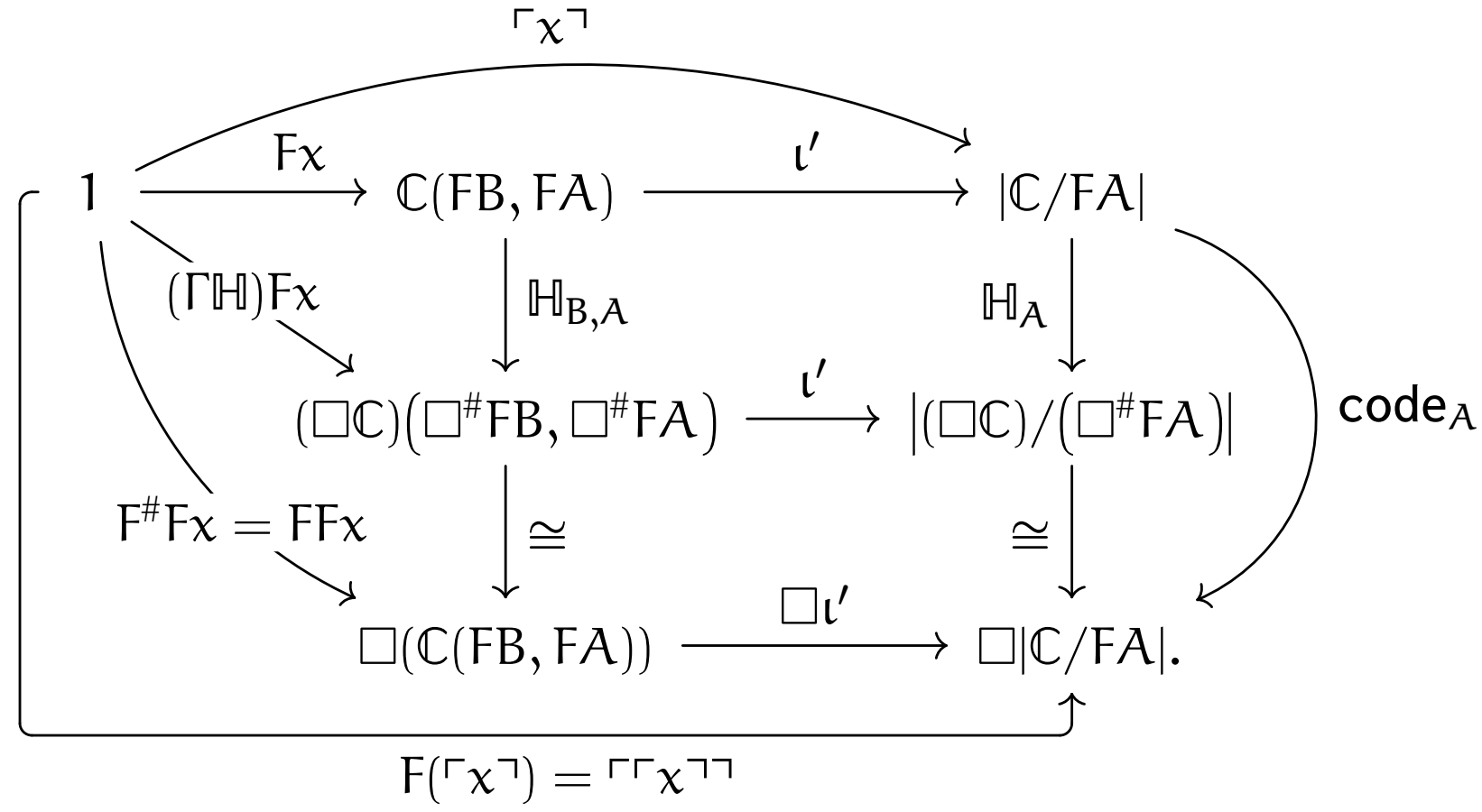

Here, the assumption $(\Gamma\mathbb{H})Fx = F^{\#}Fx$ is used to commute the left-hand triangle. $\square$

**Remark 4.18** Beyond Remark 4.11, the endofunctor $\Box = \gamma_{\mathbb{C}} \circ F : \mathcal{B} \to \mathcal{B}$ can be interpreted as a *provability modality* in the setting of internal categories.

Suppose that $\mathcal{B}$ is a syntactic category of some theory. Then, a $\mathcal{B}$-internal category $\mathbb{C}$ represents another theory formalized within $\mathcal{B}$. This situation is analogous to the internalization required for defining provability predicates. Under this interpretation, the components of $\square$ correspond to the following logical notions:

- The functor $F : \mathcal{B} \to \Gamma\mathbb{C}$ acts as *Gödel coding* (and taking its numerals), as suggested by the construction of code structures in Theorem 4.17. It maps formulas $\varphi$ in $\mathcal{B}$ to their internal representations $\llcorner\varphi\lrcorner$ in $\mathbb{C}$.
- The global section functor $\gamma_{\mathbb{C}} : \Gamma\mathbb{C} \to \mathcal{B}$ acts as the *provability predicate* $\mathrm{Pr}(x)$. In a syntactic category, a closed formula corresponds to a subterminal $X \rightarrowtail 1$. For an ordinary category $\mathcal{C}$, the set $\Gamma_{\mathcal{C}}(X) = \mathcal{C}(1, X)$ is inhabited if and only if there exists a morphism $1 \to X$, i.e., $X$ is provable. Internalizing this, $\gamma_{\mathbb{C}}(X)$ becomes an object in $\mathcal{B}$ expressing "$X$ is provable in $\mathbb{C}$." Notably, this interpretation extends the provability predicate from subterminals to the entire category.

Consequently, the composite $\square = \gamma_{\mathbb{C}} \circ F$ maps a formula $\varphi$ to the formula $\mathrm{Pr}(\llcorner\varphi\lrcorner)$, exactly capturing the provability modality in provability logic. This interpretation appears to originate from Joyal's work, as discussed in certain accounts [8, 13].

## 4.4 Pre-geminal models of finite limit theories

Recall that a category with chosen finite limits can be viewed as a **Set**-model of the finite limit theory $\mathcal{T}_{\mathrm{FLcat}}$ (see Section 2.5). The notion of a pre-geminal category is naturally generalized to models of any finite limit theory containing $\mathcal{T}_{\mathrm{FLcat}}$. This generalization is particularly useful for constructing examples of both pre-geminal and geminal categories (for the latter, see Example 6.6). An alternative construction of geminal categories based on presheaf categories can be found in Ramesh [38, Section 6.5].

Let $\mathcal{T}$ be a finite limit theory equipped with a theory morphism $\varphi : \mathcal{T}_{\mathrm{FLcat}} \to \mathcal{T}$. We write $\|{-}\|$ for the functor $\varphi^* : \mathbf{Mod}(\mathcal{T}, \mathcal{B}) \to \mathbf{Mod}(\mathcal{T}_{\mathrm{FLcat}}, \mathcal{B})$ induced by $\varphi$. For instance, if $M$ is a **Set**-model of $\mathcal{T}$, then $\|M\|$ is a category with chosen finite limits.

**Definition 4.19 (Pre-geminal model)** Let $\mathcal{T}$ be a finite limit theory equipped with a theory morphism $\varphi : \mathcal{T}_{\mathrm{FLcat}} \to \mathcal{T}$. A *pre-geminal model of* $\mathcal{T}$ consists of the following data:

- A model $M_B \in \mathbf{Mod}(\mathcal{T}, \mathbf{Set})$.
- A model $M_C \in \mathbf{Mod}(\mathcal{T}, \|M_B\|)$.
- A morphism $f : M_B \to \Gamma M_C$ in $\mathbf{Mod}(\mathcal{T}, \mathbf{Set})$.
- A morphism $h : M_C \to \square M_C$ in $\mathbf{Mod}(\mathcal{T}, \|M_B\|)$ satisfying $\Gamma h \circ f = \|f\|^{\#}_{M_C} \circ f$ in $\mathbf{Mod}(\mathcal{T}, \mathbf{Set})$.

$$M_B \xrightarrow{\ f\ } \Gamma M_C \underset{\Gamma h}{\overset{\|f\|^{\#}}{\rightrightarrows}} \Gamma\square M_C$$

Here, the endofunctor $\square : \|M_B\| \to \|M_B\|$ is defined from the $\|M_B\|$-internal category $\|M_C\|$ and the functor $\|f\| : \|M_B\| \to \Gamma\|M_C\|$ via Definition 4.14.

As before, we identify the codomains of $\Gamma h : \Gamma M_C \to \Gamma \Box M_C$ and $\|f\|^\# : \Gamma M_C \to \Gamma F M_C$ via the canonical isomorphism. Note that a pre-geminal model of $\mathcal{T}_{\mathrm{FLcat}}$ (equipped with the identity theory morphism) is precisely an ordinary pre-geminal category.

> **Example 4.20 (Initial pre-geminal models)** Let $M_B$ be the initial **Set**-model of $\mathcal{T}$, which always exists by Theorem 2.62. If there also exists an initial $\|M_B\|$-model $M_C$ of $\mathcal{T}$, one can construct a pre-geminal model of $\mathcal{T}$ from them. Indeed, there are unique homomorphisms $f : M_B \to \Gamma M_C$ and $h : M_C \to \Box M_C$ due to the initiality of $M_B$ and $M_C$, respectively. The required equation follows from the initiality of $M_B$.

While this construction is trivial, it is useful for providing examples of pre-geminal categories via the following proposition:

> **Proposition 4.21** *If* $(M_B, M_C, f, h)$ *is a pre-geminal model of* $\mathcal{T}$*, then* $(\|M_B\|, \|M_C\|, \|f\|, \|h\|)$ *forms a pre-geminal category.*

*Proof.* The result follows directly from the fact that the functor $\|-\| = \varphi^*$ commutes with $\Gamma_*$, $\|f\|_*$ and $\|f\|^\#$, as described in Section 2.5. □

An essential task in constructing a pre-geminal model as in Example 4.20 is finding a theory $\mathcal{T}$ such that, for its initial **Set**-model $M$, there exists an initial $\|M\|$-model of $\mathcal{T}$. In other words, we need to find a categorical structure that is rich enough to internalize the construction of its own initial model.

Joyal's *arithmetic universes* were introduced for this very reason. Proofs of the fact that the initial arithmetic universe contains an internal initial arithmetic universe can be found in Morrison [31] or Vickers [42]. Similarly, one could take $\mathcal{T}$ to be a theory with even stronger structures, such as elementary toposes with a natural number object. Additional axioms can also be imposed on the theory.

In principle, these constructions suggest that pre-geminal (and geminal) categories can be obtained from any sufficiently expressive logical system. We do not pursue a detailed, systematic description of these constructions here. For the relationship between the geminal category arising from the initial arithmetic universe and those arising from first-order arithmetic, we refer the reader to Ramesh [38, Section 6.4].

# 5 Bootstrapping and Löb's theorem

In the previous section, we showed that any pre-geminal category induces code structures on both representable and codomain fibrations. For these code structures, we employ the following trick: we apply the fixed point theorem to codomain fibrations to obtain the data required to satisfy the premises of the fixed point theorem for other fibrations. Following the usage of the term "bootstrapping" by Ramesh [38], we call this trick the *bootstrapping argument*.

As a result, we obtain *Löb's theorem for pre-geminal categories*, which slightly generalizes the original theorem for geminal categories [38, Observation 5.24]. By minimizing internal reasoning, our formulation clarifies the essence of Ramesh's bootstrapping argument and provides a simplified proof of the theorem.

## 5.1 An idea of the bootstrapping argument

Before turning to the proof, we first provide an informal overview of the bootstrapping argument.

The construction of a fixed point based on the diagonal argument typically requires a form of surjection from some objects to functions acting on them. For example, Lawvere's fixed point theorem [26] requires a morphism $X \to Y^X$ that induces a surjection between their global elements. Similarly, our fixed point theorem (Theorem 3.14) requires a morphism that induces an essentially surjective reindexing functor. This surjectivity condition is crucial to the theorem; indeed, its contrapositive is also frequently used to prove the non-existence of such morphisms, as seen in Cantor's theorem.

To use the theorem *positively*, one must construct such a surjection. In the case of the diagonal lemma (Example 3.15), it is achieved via Gödel numbering, where formulas on natural numbers N are represented by natural numbers itself, namely, $\ulcorner \mathcal{E}_N \urcorner = N$. In the derivation of Löb's theorem (Example 3.19), this surjection is constructed by the diagonal lemma itself.

The latter example may be similar to the idea of the bootstrapping argument. The essential idea is as follows: to obtain a surjection from X to $Y^X$, it suffices to construct an isomorphism between them. Thus, we could try to apply the fixed point theorem itself to the operation $X \mapsto Y^X$.

Naively, this operation $X \mapsto Y^X$ acts on the collection of all objects, or type-theoretically, on the universe $\mathcal{U}$. To apply the fixed point theorem to this operation, one requires a surjection $\mathcal{X} \to \mathcal{U}^{\mathcal{X}}$ for some $\mathcal{X}$. If $\mathcal{U}$ were to contain *all* types literally, such a surjection could be constructed as follows: let $\mathcal{W}$ be the universe of all pointed types (i.e., pairs $(X, x)$ of a type X and its object $x : X$). Then, an inclusion $i : \mathcal{U}^{\mathcal{W}} \to \mathcal{W}, a \mapsto (\mathcal{U}^{\mathcal{W}}, a)$ induces the surjection $(-) \circ i : \mathcal{U}^{\mathcal{W}} \to \mathcal{U}^{\mathcal{U}^{\mathcal{W}}}$. In the fibrational reading, $\mathcal{W}$ and $\mathcal{U}$ correspond to the codomain fibration of $\mathcal{B}$, which may make the argument seem simpler.

If this works, a fixed point for any action on any type is constructed by utilizing the diagonal argument *twice*. This strategy appears implicitly in Pitts and Taylor [37], who establish a negative result — the non-existence of a generic object in the codomain fibration over any non-trivial locally Cartesian closed category — using a type-theoretic analog of Russell's paradox. In fact, Russell's paradox itself can be viewed as a degenerate instance of this double diagonal argument:

1. The first application of the diagonal argument constructs the set $\{x \mid x \notin x\}$, leading to the logical equivalence $\varphi \leftrightarrow \neg\varphi$, i.e., an isomorphism between $\varphi$ and $\varphi \to \bot$.
2. The second application derives a contradiction $\bot$ from this equivalence.

Under the Curry–Howard correspondence, the proposition $\varphi \to \bot$ corresponds to the type of functions from $\varphi$ to $\bot$. This explains the second derivation as a degenerate instance of the diagonal argument. Because Pitts and Taylor [37] lift the propositional logic to type theory, they employ the diagonal argument twice in an explicit manner.

Ramesh [38] employs this bootstrapping argument in a *positive* manner to derive *Löb's theorem for introspective theories.* The proof we provide below also follows a similar development. We argue that this success is due to applying the method to *intensional* settings. In such a framework, the construction of fixed points yields *modal* or *intensional* fixed points in the sense of Kavvos [18–21], the existence of which is a key property of the type-theoretic counterpart of the provability logic GL. We will discuss this correspondence further in Section 6.3.

Note that while there exist positive and extensional applications of the diagonal argument, such as Kleene's first recursion theorem [18], we do not know whether there is a positive application of the bootstrapping argument in an extensional setting.

## 5.2 Löb's theorem for pre-geminal categories

The bootstrapping argument for a pre-geminal category proceeds in two stages. First, we apply the fixed point theorem to the codomain fibration with codes. In this case, the premises of the theorem are always satisfied. This first application yields a specific isomorphism, which we then use to satisfy the premises for a second application. In the second step, we apply the fixed point theorem once more, but this time to the representable fibrations with codes.

> **Lemma 5.1** *Let $\mathcal{B}$ be a category with pullbacks, and consider its codomain fibration. For any monomorphism $i : A \to B$ in $\mathcal{B}$, the reindexing functor $i^* : \mathcal{B}/B \to \mathcal{B}/A$ is essentially surjective.*

*Proof.* Let $x : X \to A$ be any object in the slice category $\mathcal{B}/A$. If we define $y = i \circ x \in \mathcal{B}/B$, an isomorphism $i^*(y) \cong x$ in $\mathcal{B}/A$ is obtained by the following pullback diagram:

$$\begin{array}{ccc} X & \xrightarrow{\mathrm{id}} & X \\ {\scriptstyle x}\downarrow & \lrcorner & \downarrow{\scriptstyle x} \\ A & \xrightarrow{\mathrm{id}} & A \\ {\scriptstyle \mathrm{id}}\downarrow & \lrcorner & \downarrow{\scriptstyle i} \\ B & \xrightarrow{i} & B. \end{array} \qquad y = i \circ x : X \to B$$

Thus, the functor $i^*$ is essentially surjective. □

> **Lemma 5.2** *Let $(\mathcal{B}, \mathbb{C}, F, \mathbb{H})$ be a pre-geminal category. For any $P \in \mathcal{B}/|\mathbb{C}|$ in $\mathcal{B}$, there exists $A \in \mathcal{B}$ such that the following forms a pullback diagram:*

$$\begin{array}{ccc} A & \longrightarrow & X \\ {\scriptstyle !}\downarrow & \lrcorner & \downarrow{\scriptstyle P} \\ 1 & \xrightarrow{FA} & |\mathbb{C}|. \end{array}$$

*Proof.* We apply the fixed point theorem (Theorem 3.14) to the codomain fibration over $\mathcal{B}$ equipped with the code structure induced by the pre-geminal category. Consider the object $|\mathbb{C}^{\to}| \in \mathcal{B}$ and the inclusion morphism $\iota : |\mathbb{C}/F|\mathbb{C}^{\to}|| \to |\mathbb{C}^{\to}|$. As noted in Remark 2.35, this $\iota$ is a monomorphism. Thus, by the previous lemma, the reindexing functor $\iota^*$ is essentially surjective, which shows that the pair $(|\mathbb{C}^{\to}|, \iota)$ satisfies the premises of the fixed point theorem. The resulting isomorphism directly yields the assertion of the lemma. □

**Corollary 5.3** *Let* $(\mathcal{B}, \mathbb{C}, F, \mathbb{H})$ *be a pre-geminal category. For any object* $A \in \mathcal{B}$*, there exists an object* $B \in \mathcal{B}$ *such that* $\mathbb{C}(FB, FA) \cong B$.

*Proof.* We apply Lemma 5.2 to the morphism $\mathrm{dom}_{FA} : |\mathbb{C}/FA| \to |\mathbb{C}|$. This yields an object $B \in \mathcal{B}$ such that the following square is a pullback:

$$\begin{array}{ccc} B & \longrightarrow & |\mathbb{C}/FA| \\ {\scriptstyle !}\downarrow & \lrcorner & \downarrow{\scriptstyle \mathrm{dom}_{FA}} \\ 1 & \xrightarrow{FB} & |\mathbb{C}|. \end{array}$$

On the other hand, it is straightforward to check that $\mathbb{C}(FB, FA)$ is obtained by the same pullback as follows:

$$\begin{array}{ccc} \mathbb{C}(FB, FA) & \xrightarrow{\iota'} & |\mathbb{C}/FA| \\ {\scriptstyle !}\downarrow & \lrcorner & \downarrow{\scriptstyle \mathrm{dom}_{FA}} \\ 1 & \xrightarrow{FB} & |\mathbb{C}|. \end{array}$$

Thus, it follows that $B \cong \mathbb{C}(FB, FA)$. □

This corollary provides the data satisfying the premises of the fixed point theorem for representable fibrations. Now we perform the second application of the fixed point theorem:

**Theorem 5.4 (Löb's theorem for pre-geminal categories)** *Let* $(\mathcal{B}, \mathbb{C}, F, \mathbb{H})$ *be a pre-geminal category. For any object* $A \in \mathcal{B}$ *and any morphism* $f : \Box A \to A$*, there exists a global element* $a \in_1 A$ *such that* $f(Fa) = a$.

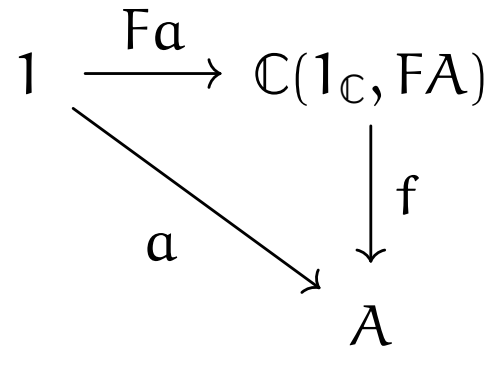

*Proof.* First, note that for any $a \in_1 A$, we have $Fa \in_1 \mathbb{C}(F1, FA) = \mathbb{C}(1_{\mathbb{C}}, FA) = \Box A$, and thus $f(Fa)$ is a well-defined global element of $A$.

Consider the representable fibration $\mathrm{dom}_A : \mathcal{B}/A \to \mathcal{B}$ with codes induced by the pre-geminal category. By Corollary 5.3, there exists $B \in \mathcal{B}$ such that $\mathbb{C}(FB, FA) \cong B$, i.e., $\ulcorner \mathcal{E}_B \urcorner \cong B$. This isomorphism ensures that the premises of the fixed point theorem (Theorem 3.14) are satisfied. The theorem then yields that, for any $f \in \mathcal{B}(\mathbb{C}(F1, FA), A)$, there exists $a \in \mathcal{B}(1, A)$ such that $f \circ Fa = a$. This is precisely what we wanted. $\Box$

It can also be shown that the global element $a \in_1 A$ is unique, which follows from a general argument using the fact that $\mathcal{B}$ has finite limits. For the proof, we refer the reader to Ramesh [38, Theorem 4.4].

Theorem 5.4 states that any morphism $\Box A \to A$ yields a global element $1 \to A$, which can be viewed as a categorical analog of Löb's theorem. Indeed, according to the interpretation in Remark 4.18, this result recovers Löb's theorem for the internal logic of the category, particularly when $A$ is a subterminal. Furthermore, the fixed point equation $f(Fa) = a$ may be interpreted as a form of *intensional recursion*; see Kavvos [21] for a discussion.

The core idea of our proof remains identical to that of Ramesh's main theorem, *Löb's theorem for introspective theories* [38, Theorem 4.19]. However, by focusing on geminal categories rather than introspective theories, our formulation bypasses any informal internal reasoning, which significantly simplifies the proof. Notably, Ramesh's main theorem follows from our formulation, as any introspective theory induces a geminal category [38, Construction 5.18]. Additionally, our Theorem 5.4 offers a slight generalization from geminal to pre-geminal categories, though this should be regarded as a byproduct of our reorganization rather than a practical extension.

**Remark 5.5** The bootstrapping technique is not limited to representable fibrations. For an arbitrary fibration $p : \mathcal{E} \to \mathcal{B}$, if there is an "internalization" of $p$, i.e., a morphism $\ulcorner p \urcorner : \ulcorner \mathcal{E} \urcorner \to \ulcorner \mathcal{B} \urcorner$ in $\mathcal{B}$ with appropriate properties, then one can apply Lemma 5.2 to $P = \ulcorner p \urcorner$. This application is expected to yield $A \in \mathcal{B}$ such that $\ulcorner \mathcal{E}_A \urcorner \cong A$ in $\mathcal{B}$, which provides necessary data to apply the fixed point theorem to the fibration $p$.

# 6 Geminal categories and the Gödel–Löb axiom

Since an analog of Löb's theorem holds for pre-geminal categories, it is natural to investigate how the endofunctor $\Box$ corresponds to the modality of the provability logic GL. To establish this correspondence, one would expect to derive not only the ordinary Löb's theorem but also its internalized version: *the Gödel–Löb axiom*, $\Box(\Box A \to A) \to \Box A$.

Naively, the Gödel–Löb axiom should be provable by internalizing the proof of Löb's theorem within a category. To perform this idea, we introduce the notion of an *internal pre-geminal category*. Notably, the original notion of a *geminal category* [38, Definition 5.15] emerges as a natural extension of a pre-geminal category that induces an internal pre-geminal category within itself. Our reorganization makes it manageable to provide a formal proof of a categorical counterpart of the Gödel–Löb axiom in any geminal category, not limited to those with exponentials (see Remark 6.13). We also discuss the structure arising from the endofunctor $\Box : \mathcal{B} \to \mathcal{B}$ associated with geminal categories and compare it with structures studied in the context of modal calculi.

## 6.1 Geminal categories

We introduce the notion of *internal pre-geminal categories.* In this internal analog, the role of a category $\mathcal{B}$ in a pre-geminal category is played by an internal category $\mathbb{B}$ within some ambient category $\mathcal{A}$. Furthermore, as discussed in Section 2.4.6, a $\mathcal{B}$-internal category $\mathbb{C}$ corresponds to a $\Gamma\mathbb{B}$-internal category. We begin by defining an analog of the endofunctor $\Box : \mathcal{B} \to \mathcal{B}$ in this setting.

**Definition 6.1 (The endofunctor $\Box_{\mathbb{B}}$)** Let $\mathbb{B}$ be an $\mathcal{A}$-internal category with chosen finite limits, $\mathbb{C}$ be a $\Gamma\mathbb{B}$-internal category with chosen finite limits, and $\mathbb{F} : \mathbb{B} \to \gamma_{\mathbb{B}}\mathbb{C}$ be an $\mathcal{A}$-internal functor strictly preserving finite limits. We define the endofunctor $\Box_{\mathbb{B}} : \Gamma\mathbb{B} \to \Gamma\mathbb{B}$ by the composite $\gamma_{\mathbb{C}} \circ (\Gamma\mathbb{F})$.

$$\Gamma\mathbb{B} \xrightarrow{\Gamma\mathbb{F}} \Gamma\mathbb{C} \xrightarrow{\gamma_{\mathbb{C}}} \Gamma\mathbb{B}, \qquad \Box_{\mathbb{B}} : \Gamma\mathbb{B} \to \Gamma\mathbb{B}$$

Since $\Box_{\mathbb{B}}$ preserves finite limits, we may consider the change of base along $\Box_{\mathbb{B}}$ for $\Gamma\mathbb{B}$-internal categories. Note that $\Box_{\mathbb{B}}$ is an ordinary endofunctor, not an internal one.

**Definition 6.2 (Internal pre-geminal category)** Let $\mathcal{A}$ be a category with finite limits. An $\mathcal{A}$-*internal pre-geminal category* consists of the following data:

(I*) An $\mathcal{A}$-internal category $\mathbb{B}$ with chosen finite limits.
(II*) A $\Gamma\mathbb{B}$-internal category $\mathbb{C}$ with chosen finite limits.
(III*) An $\mathcal{A}$-internal functor $\mathbb{F} : \mathbb{B} \to \gamma_{\mathbb{B}}\mathbb{C}$ strictly preserving finite limits.
(IV*) A $\Gamma\mathbb{B}$-internal functor $\mathbb{H} : \mathbb{C} \to \Box_{\mathbb{B}}\mathbb{C}$ satisfying $\gamma_{\mathbb{B}}\mathbb{H} \circ \mathbb{F} = \mathbb{F}^{\#}_{\mathbb{C}} \circ \mathbb{F}$.

$$\mathbb{B} \xrightarrow{\mathbb{F}} \gamma_{\mathbb{B}}\mathbb{C} \underset{\gamma_{\mathbb{B}}\mathbb{H}}{\overset{\mathbb{F}^{\#}}{\rightrightarrows}} \gamma_{\mathbb{B}}\Box_{\mathbb{B}}\mathbb{C}$$

In condition (IV*), we identify the codomains of $\gamma_{\mathbb{B}}\mathbb{H} : \gamma_{\mathbb{B}}\mathbb{C} \to \gamma_{\mathbb{B}}\Box_{\mathbb{B}}\mathbb{C}$ and $\mathbb{F}^{\#} : \gamma_{\mathbb{B}}\mathbb{C} \to \gamma_{\gamma_{\mathbb{B}}\mathbb{C}}(\Gamma\mathbb{F})\mathbb{C}$ via the canonical isomorphism obtained from Proposition 2.61,

$$\gamma_{\mathbb{B}}\Box_{\mathbb{B}}\mathbb{C} = \gamma_{\mathbb{B}}\gamma_{\mathbb{C}}(\Gamma\mathbb{F})\mathbb{C} \cong \gamma_{\gamma_{\mathbb{B}}\mathbb{C}}(\Gamma\mathbb{F})\mathbb{C}.$$

Note that a **Set**-internal pre-geminal category is precisely an ordinary pre-geminal category. The whole structures of ordinary and $\mathcal{A}$-internal geminal categories are illustrated as follows:

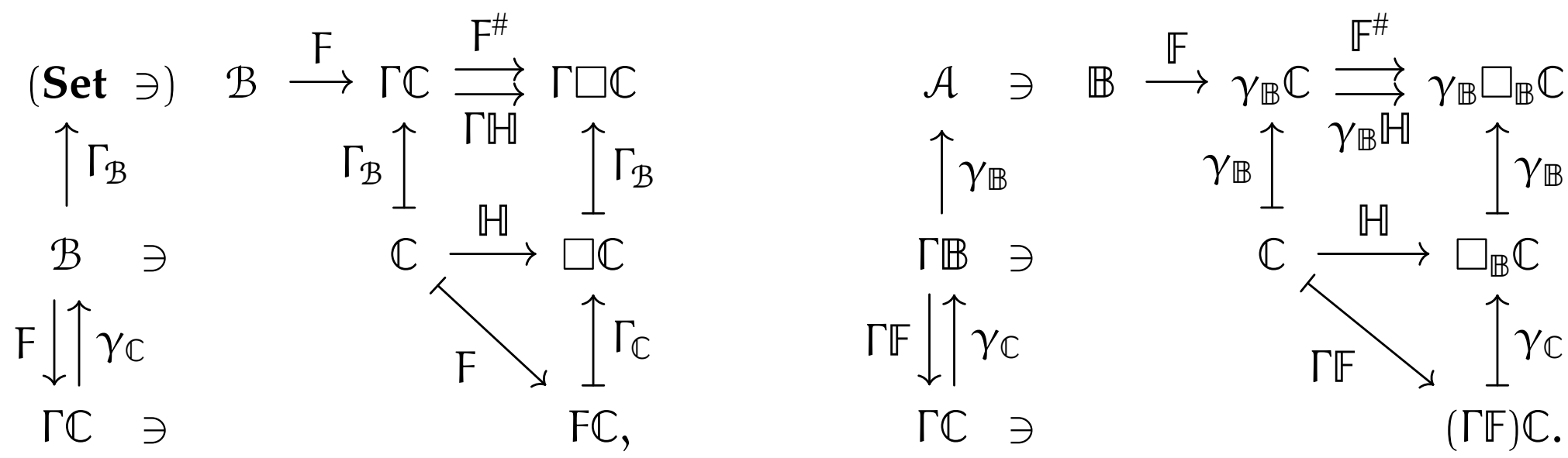

> **Proposition 6.3** *For an $\mathcal{A}$-internal pre-geminal category* $(\mathbb{B}, \mathbb{C}, \mathbb{F}, \mathbb{H})$, *the quadruple* $(\Gamma\mathbb{B}, \mathbb{C}, \Gamma\mathbb{F}, \mathbb{H})$ *forms a pre-geminal category.*

*Proof.* We verify that the components satisfy conditions (I)–(IV) in the definition of pre-geminal categories (Definition 4.15). It is immediate that $(\Gamma\mathbb{B}, \mathbb{C}, \Gamma\mathbb{F})$ provides components (I)–(III). Since the endofunctor $\Box_{\mathbb{B}}$ coincides with $\Box_{\Gamma\mathbb{B}} : \Gamma\mathbb{B} \to \Gamma\mathbb{B}$ defined from this triple, we may take $\mathbb{H} : \mathbb{C} \to \Box_{\Gamma\mathbb{B}}\mathbb{C}$ as component (IV). The required identity for $\mathbb{H}$ is obtained by applying $\Gamma_{\mathcal{A}}$ to the assumption $\gamma_{\mathbb{B}}\mathbb{H} \circ \mathbb{F} = \mathbb{F}^{\#} \circ \mathbb{F}$, combined with Proposition 2.50 and Corollary 2.60. $\square$

A key observation is that, given an ordinary pre-geminal category $(\mathcal{B}, \mathbb{C}, \mathsf{F}, \mathbb{H})$, we can induce most of the data for a $\mathcal{B}$-internal pre-geminal category. This is not surprising, as the definition of a pre-geminal category arises from requiring $\mathcal{B}$ to internalize sets and functions concerning its own structure.

> **Lemma 6.4** *Let* $(\mathcal{B}, \mathbb{C}, \mathsf{F}, \mathbb{H})$ *be a pre-geminal category and assume that* $\mathbb{H}$ *strictly preserves finite limits. Then, the quadruple* $(\mathbb{C}, \mathsf{F}\mathbb{C}, \mathbb{H}, \mathsf{F}\mathbb{H})$ *satisfies all the conditions for a $\mathcal{B}$-internal pre-geminal category, except possibly for the identity in* (IV*) *in Definition 6.2.*

In terms of the previous figure, this indicates that the diagram for a pre-geminal category can be extended downwards as follows:

$$
\begin{array}{ccccccccc}
(\mathbf{Set} & \ni) & \mathcal{B} & \xrightarrow{\mathsf{F}} & \Gamma\mathbb{C} & \overset{\mathsf{F}^{\#}}{\underset{\Gamma\mathbb{H}}{\rightrightarrows}} & \Gamma\Box\mathbb{C} & & \\
\uparrow \Gamma_{\mathcal{B}} & & & & \uparrow & & \uparrow & & \\
\mathcal{B} & \ni & & & \mathbb{C} & \xrightarrow{\mathbb{H}} & \Box\mathbb{C} & \overset{\mathbb{H}^{\#}}{\underset{\Box\mathbb{H}}{\rightrightarrows}} & \Box\Box\mathbb{C} \\
\mathsf{F}\downarrow\uparrow\gamma_{\mathbb{C}} & & & & & \searrow & \uparrow & \searrow & \uparrow \\
\Gamma\mathbb{C} & \ni & & & & & \mathsf{F}\mathbb{C} & \xrightarrow{\mathsf{F}\mathbb{H}} & \mathsf{F}\Box\mathbb{C} \\
\Gamma\mathbb{H}\downarrow\uparrow\gamma_{\mathsf{F}\mathbb{C}} & & & & & & & \searrow & \uparrow \\
(\Gamma\mathsf{F}\mathbb{C} \cong)\ \Gamma\Box\mathbb{C} & \ni & & & & & & & (\Gamma\mathbb{H})\mathsf{F}\mathbb{C}.
\end{array}
$$

*Proof.* Immediately, $(\mathbb{C}, \mathsf{F}\mathbb{C}, \mathbb{H})$ constitutes components (I*)–(III*). For component (IV*), notice that $\mathsf{F}\Box : \mathcal{B} \to \Gamma\mathbb{C}$ can be rearranged as

$$
\begin{aligned}
\mathsf{F}\square &= \mathsf{F}\gamma_{\mathbb{C}}\mathsf{F} \\
&\cong \gamma_{\mathsf{F}\mathbb{C}}\mathsf{F}^{\#}\mathsf{F} && \text{by Proposition 2.51} \\
&= \gamma_{\mathsf{F}\mathbb{C}}(\Gamma\mathbb{H})\mathsf{F} && \text{since } \Gamma\mathbb{H}\circ\mathsf{F} = \mathsf{F}^{\#}\circ\mathsf{F} \\
&= \square_{\mathbb{C}}\mathsf{F}.
\end{aligned}
$$

Thus, the $\Gamma\mathbb{C}$-internal functor $\mathsf{F}\mathbb{H} : \mathsf{F}\mathbb{C} \longrightarrow \mathsf{F}\square\mathbb{C}$ can be identified with a functor of the form $\mathsf{F}\mathbb{C} \longrightarrow \square_{\mathbb{C}}\mathsf{F}\mathbb{C}$, which is required for component (IV*). □

However, the identity $\square\mathbb{H} \circ \mathbb{H} = \mathbb{H}^{\#} \circ \mathbb{H}$, which corresponds to the required equation in (IV*), cannot generally be expected to hold. This is because the definition of a pre-geminal category only imposes an identity between ordinary functors, and there is nothing that induces an identity between internal functors. The requirement that this identity holds precisely coincides with the notion of a (*compactly presented*) *geminal category* introduced by Ramesh [38, Definition 5.15].

**Definition 6.5 (Geminal category)** A *geminal category* is a pre-geminal category $(\mathcal{B}, \mathbb{C}, \mathsf{F}, \mathbb{H})$ such that $\mathbb{H}$ strictly preserves finite limits and $\square\mathbb{H} \circ \mathbb{H} = \mathbb{H}^{\#}_{\mathsf{F}\mathbb{C}} \circ \mathbb{H}$ holds.

$$
\mathbb{C} \xrightarrow{\ \mathbb{H}\ } \square\mathbb{C} \underset{\square\mathbb{H}}{\overset{\mathbb{H}^{\#}}{\rightrightarrows}} \square\square\mathbb{C}
$$

**Example 6.6** The notion of pre-geminal models of a finite limit theory $\mathcal{T}$ (Definition 4.19) naturally extends to the notion of *geminal models of* $\mathcal{T}$, and most arguments also extend to the latter. In particular, geminal categories can be constructed from geminal models, and any initial pre-geminal model extends to a geminal model. Consequently, one can construct a geminal category from any finite limit theory $\mathcal{T}$ satisfying the condition stated in Example 4.20. For instance, the initial arithmetic universe and the initial topos with natural numbers object both give rise to geminal categories.

By its design, any geminal category induces an internal pre-geminal category. In fact, a stronger result holds: any geminal category induces an *internal geminal category*. This implies that the notion of a geminal category constitutes a kind of fixed point for the internalization process, requiring no further extension for self-internalization. This closure suggests that a geminal category is not merely an ad hoc notion, but rather a robust concept that may serve as a universal model for self-internalization.

**Definition 6.7 (Internal geminal category)** For a category $\mathcal{A}$ with finite limits, an $\mathcal{A}$-*internal geminal category* is an $\mathcal{A}$-internal pre-geminal category $(\mathbb{B}, \mathbb{C}, \mathbb{F}, \mathbb{H})$ such that $\mathbb{H}$ strictly preserves finite limits and $\square_{\mathbb{B}}\mathbb{H} \circ \mathbb{H} = \mathbb{H}^{\#}_{(\Gamma\mathbb{F})\mathbb{C}} \circ \mathbb{H}$ holds.

**Theorem 6.8 (Ramesh [38, Theorem 5.16])** *For a geminal category* $(\mathcal{B}, \mathbb{C}, \mathsf{F}, \mathbb{H})$, *the quadruple* $(\mathbb{C}, \mathsf{F}\mathbb{C}, \mathbb{H}, \mathsf{F}\mathbb{H})$ *forms a* $\mathcal{B}$-*internal geminal category.*

*Proof.* By Lemma 6.4 and the definition of geminal categories, we have already established that $(\mathbb{C}, \mathsf{F}\mathbb{C}, \mathbb{H}, \mathsf{F}\mathbb{H})$ forms a $\mathcal{B}$-internal pre-geminal category. Further, the $\Gamma\mathbb{C}$-internal functor $\mathsf{F}\mathbb{H} : \mathsf{F}\mathbb{C} \to \mathsf{F}\Box\mathbb{C}$ strictly preserves finite limits, as $\mathbb{H}$ does so (Proposition 2.55 (2)). Thus, it suffices to verify the identity $(\Box_{\mathbb{C}}\mathsf{F}\mathbb{H}) \circ (\mathsf{F}\mathbb{H}) = (\mathsf{F}\mathbb{H})^{\#} \circ (\mathsf{F}\mathbb{H})$. This is obtained by applying $\mathsf{F} : \mathcal{B} \to \Gamma\mathbb{C}$ to the identity $\Box\mathbb{H} \circ \mathbb{H} = \mathbb{H}^{\#} \circ \mathbb{H}$. Indeed, the isomorphism $\mathsf{F}\Box \cong \Box_{\mathbb{C}}\mathsf{F}$ described in the proof of Lemma 6.4 yields the identity $\mathsf{F}\Box\mathbb{H} = \Box_{\mathbb{C}}\mathsf{F}\mathbb{H}$ (noting that we identify the domains and codomains of them under the canonical isomorphisms). Moreover, we also have

$$\begin{aligned} (\mathsf{F}\mathbb{H})^{\#}_{(\Gamma\mathbb{H})\mathsf{F}\mathbb{C}} &= (\mathsf{F}\mathbb{H})^{\#}_{\mathsf{F}^{\#}\mathsf{F}\mathbb{C}} && \text{since } \Gamma\mathbb{H} \circ \mathsf{F} = \mathsf{F}^{\#} \circ \mathsf{F} \\ &= \mathsf{F}\mathbb{H}^{\#}_{\mathsf{F}\mathbb{C}} && \text{by Proposition 2.59.} \end{aligned}$$

Combining these, we obtain the required identity. □

Diagrammatically, this means that the diagram of a geminal category can be extended downwards to create a structure consisting of exactly the same kind of data at the next level, as shown below. By repeating this process, the diagram can be extended infinitely downwards. Indeed, Ramesh's primary definition of a geminal category [38, Definition 5.11] is based on this picture of infinite downward extension.

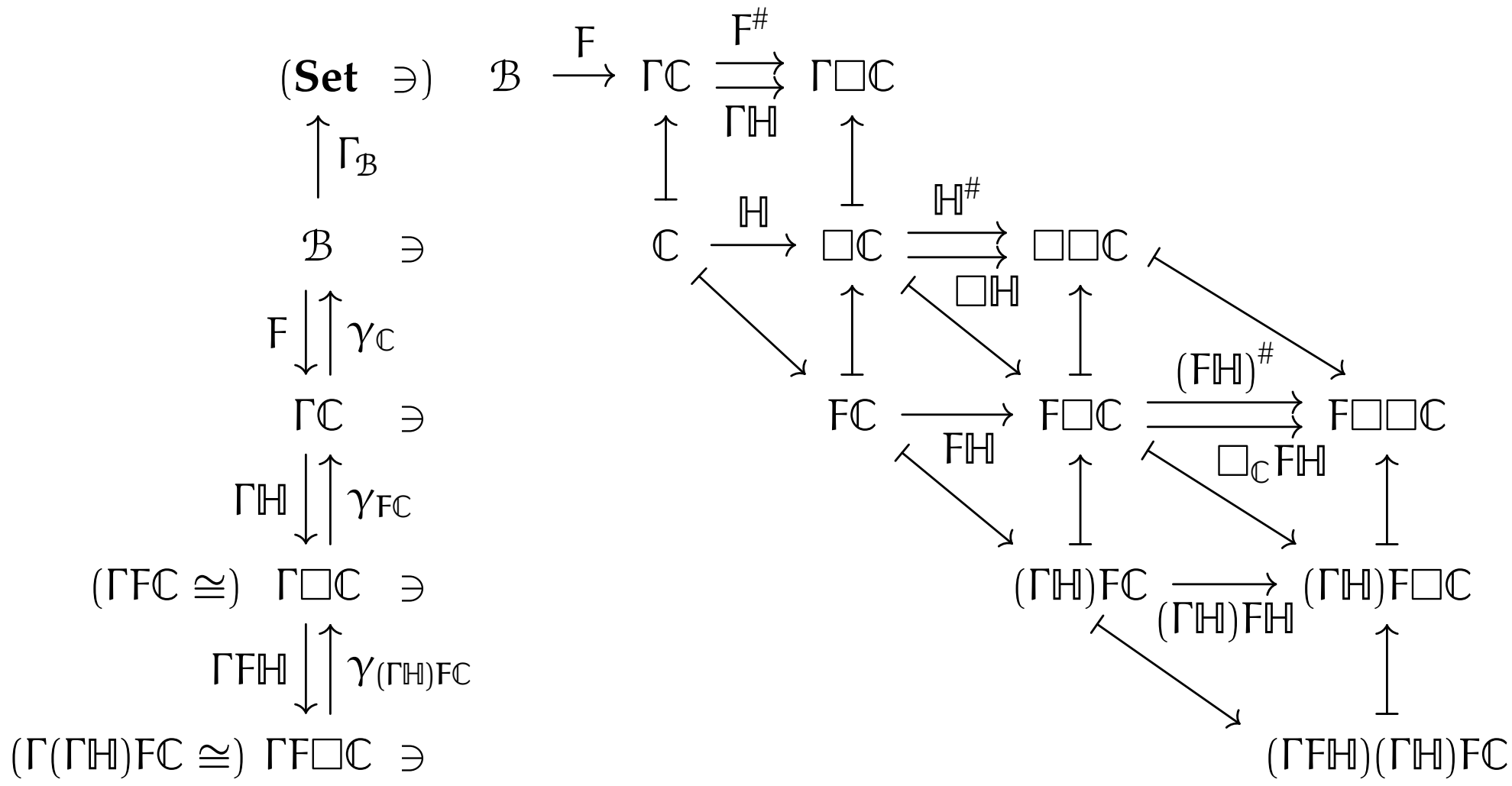

## 6.2 The Gödel–Löb axiom

We prove the Gödel–Löb axiom for geminal categories by internalizing the proof of Löb's theorem into the setting of internal pre-geminal categories. While such an internalization may appear tedious at first glance, the actual task is smaller than it first seems. For instance, the counterpart to our first fixed point result (Corollary 5.3) requires no internal reasoning; it suffices to apply the corollary to the ordinary pre-geminal category arising from the internal one. Internal arguments are only necessary when the fixed point theorem is applied for the second time to representable fibrations with codes. This observation allows us to provide a rigorous proof of the Gödel–Löb axiom for geminal categories.

Our first lemma is the internalized version of Lemma 3.13, specialized to (internal) representable fibrations with codes.

**Lemma 6.9** *Let* $(\mathbb{B}, \mathbb{C}, \mathbb{F}, \mathbb{H})$ *be an* $\mathcal{A}$*-internal pre-geminal category. For any* $A, B \in_1 \mathbb{B}$ *and any* $i \in_1 \mathbb{B}(\mathbb{C}(\mathbb{F}B, \mathbb{F}A), B)$, *one can construct a global element* $d_i \in_1 \mathbb{B}(\mathbb{C}(\mathbb{F}B, \mathbb{F}A), \mathbb{C}(1_{\mathbb{C}}, \mathbb{F}A))$ *such that the following commutes:*

$$\begin{array}{ccc} \mathbb{B}(B, A) & \xrightarrow{D_i} & \mathbb{B}(1_{\mathbb{B}}, A) \\ {\scriptstyle \langle \mathbb{F}_{B,A},\, d_i \circ ! \rangle}\downarrow & & \downarrow{\scriptstyle \mathbb{F}_{1_{\mathbb{B}},A}} \\ \mathbb{B}(1_{\mathbb{B}}, \mathbb{C}(\mathbb{F}B, \mathbb{F}A)) \times \mathbb{B}(\mathbb{C}(\mathbb{F}B, \mathbb{F}A), \mathbb{C}(1_{\mathbb{C}}, \mathbb{F}A)) & \xrightarrow{\circ_{\mathbb{B}}} & \mathbb{B}(1_{\mathbb{B}}, \mathbb{C}(1_{\mathbb{C}}, \mathbb{F}A)). \end{array}$$

*Here,* $D_i : \mathbb{B}(B, A) \to \mathbb{B}(1_{\mathbb{B}}, A)$ *is defined as the following composite:*

$$\begin{aligned} \mathbb{B}(B, A) &\xrightarrow{\langle \mathbb{F}_{B,A},\, i \circ !,\, \mathrm{id} \rangle} \mathbb{B}(1_{\mathbb{B}}, \mathbb{C}(\mathbb{F}B, \mathbb{F}A)) \times \mathbb{B}(\mathbb{C}(\mathbb{F}B, \mathbb{F}A), B) \times \mathbb{B}(B, A) \\ &\xrightarrow{\circ_{\mathbb{B}} \times \mathrm{id}} \mathbb{B}(1_{\mathbb{B}}, B) \times \mathbb{B}(B, A) \\ &\xrightarrow{\circ_{\mathbb{B}}} \mathbb{B}(1_{\mathbb{B}}, A). \end{aligned}$$

*Proof.* We define $d_i$ as the morphism obtained by applying Lemma 3.13 to the pre-geminal category $(\Gamma\mathbb{B}, \mathbb{C}, \Gamma\mathbb{F}, \mathbb{H})$. We give its explicit description below.

Let $E = \mathbb{C}(\mathbb{F}B, \mathbb{F}A) \in_1 \mathbb{B}$. Note that $i \in_1 \mathbb{B}(E, B)$. We define three global morphisms $f_1, f_2, f_3$ from $E$ in $\mathbb{B}$ as follows:

$$f_1 = \mathbb{H}_{\mathbb{F}B, \mathbb{F}A} \in_1 \mathbb{B}(\mathbb{C}(\mathbb{F}B, \mathbb{F}A), \Box(\mathbb{C}(\mathbb{F}B, \mathbb{F}A))) = \mathbb{B}(E, \mathbb{C}(1_{\mathbb{C}}, \mathbb{F}E)),$$

$$f_2 = \left(1 \xrightarrow{\langle !_{\mathbb{B}}(E), \mathbb{F}(i) \rangle} \mathbb{B}(E, 1_{\mathbb{B}}) \times \mathbb{B}(1_{\mathbb{B}}, \mathbb{C}(\mathbb{F}E, \mathbb{F}B)) \xrightarrow{\circ_{\mathbb{B}}} \mathbb{B}(E, \mathbb{C}(\mathbb{F}E, \mathbb{F}B))\right),$$

$$f_3 = e(E) \in_1 \mathbb{B}(E, E) = \mathbb{B}(E, \mathbb{C}(\mathbb{F}B, \mathbb{F}A)).$$

The typing of $\mathbb{H}_{\mathbb{F}B, \mathbb{F}A}$ follows from an argument similar to one in the proof of Lemma 4.16. In the above, we write $!_{\mathbb{B}}(E) \in_1 \mathbb{B}(E, 1_{\mathbb{B}})$ for the global morphism obtained from $!_{\mathbb{B}}(E) \in_1 |\mathbb{B}^{\to}|$, which is the application of $!_{\mathbb{B}} : |\mathbb{B}| \to |\mathbb{B}^{\to}|$ to $E \in_1 |\mathbb{B}|$. We define $e(E) \in_1 \mathbb{B}(E, E)$ similarly.

These global morphisms determine

$$\begin{aligned} \langle f_1, f_2, f_3 \rangle &\in_1 \mathbb{B}(E, \mathbb{C}(1_{\mathbb{C}}, \mathbb{F}E)) \times \mathbb{B}(E, \mathbb{C}(\mathbb{F}E, \mathbb{F}B)) \times \mathbb{B}(E, \mathbb{C}(\mathbb{F}B, \mathbb{F}A)) \\ &\cong \mathbb{B}(E, \mathbb{C}(1_{\mathbb{C}}, \mathbb{F}E) \times \mathbb{C}(\mathbb{F}E, \mathbb{F}B) \times \mathbb{C}(\mathbb{F}B, \mathbb{F}A)). \end{aligned}$$

On the other hand, we have the global morphism $c$ as follows, which is naturally induced by the composition morphism of $\mathbb{C}$:

$$c \in_1 \mathbb{B}(\mathbb{C}(1_{\mathbb{C}}, \mathbb{F}E) \times \mathbb{C}(\mathbb{F}E, \mathbb{F}B) \times \mathbb{C}(\mathbb{F}B, \mathbb{F}A), \mathbb{C}(1_{\mathbb{C}}, \mathbb{F}A)).$$

Then, the global morphism $d_i$ is defined as the following composite:

$$1 \xrightarrow{\langle\langle f_1, f_2, f_3\rangle, c\rangle} \mathbb{B}(E, \mathbb{C}(1_{\mathbb{C}}, \mathbb{F}E) \times \mathbb{C}(\mathbb{F}E, \mathbb{F}B) \times \mathbb{C}(\mathbb{F}B, \mathbb{F}A)) \times \mathbb{B}(\mathbb{C}(1_{\mathbb{C}}, \mathbb{F}E) \times \mathbb{C}(\mathbb{F}E, \mathbb{F}B) \times \mathbb{C}(\mathbb{F}B, \mathbb{F}A), \mathbb{C}(1_{\mathbb{C}}, \mathbb{F}A)) \xrightarrow{\circ_{\mathbb{B}}} \mathbb{B}(E, \mathbb{C}(1_{\mathbb{C}}, \mathbb{F}A)).$$

Next, we verify the required commutativity of the diagram. Let $g$ be the lower path in the given diagram, i.e., the following morphism:

$$\mathbb{B}(B, A) \xrightarrow{\langle \mathbb{F}_{B,A}, d_i \circ !\rangle} \mathbb{B}(1_{\mathbb{B}}, E) \times \mathbb{B}(E, \mathbb{C}(1_{\mathbb{C}}, \mathbb{F}A)) \xrightarrow{\circ_{\mathbb{B}}} \mathbb{B}(1_{\mathbb{B}}, \mathbb{C}(1_{\mathbb{C}}, \mathbb{F}A)).$$

By the definition of $d_i$ and the associativity of $\circ_{\mathbb{B}}$, this $g$ is rearranged into

$$\mathbb{B}(B, A) \xrightarrow{\langle \mathbb{F}_{B,A}, \langle f_1, f_2, f_3\rangle \circ !\rangle} \mathbb{B}(1_{\mathbb{B}}, E) \times \mathbb{B}(E, \mathbb{C}(1_{\mathbb{C}}, \mathbb{F}E) \times \mathbb{C}(\mathbb{F}E, \mathbb{F}B) \times \mathbb{C}(\mathbb{F}B, \mathbb{F}A))$$
$$\xrightarrow{\langle \circ_{\mathbb{B}}, c \circ !\rangle} \mathbb{B}(1_{\mathbb{B}}, \mathbb{C}(1_{\mathbb{C}}, \mathbb{F}E) \times \mathbb{C}(\mathbb{F}E, \mathbb{F}B) \times \mathbb{C}(\mathbb{F}B, \mathbb{F}A)) \times \mathbb{B}(\mathbb{C}(1_{\mathbb{C}}, \mathbb{F}E) \times \mathbb{C}(\mathbb{F}E, \mathbb{F}B) \times \mathbb{C}(\mathbb{F}B, \mathbb{F}A), \mathbb{C}(1_{\mathbb{C}}, \mathbb{F}A))$$
$$\xrightarrow{\circ_{\mathbb{B}}} \mathbb{B}(1_{\mathbb{B}}, \mathbb{C}(1_{\mathbb{C}}, \mathbb{F}A)).$$

Let us focus on the following morphism:

$$\mathbb{B}(B, A) \xrightarrow{\langle \mathbb{F}_{B,A}, \langle f_1, f_2, f_3\rangle \circ !\rangle} \mathbb{B}(1_{\mathbb{B}}, E) \times \mathbb{B}(E, \mathbb{C}(1_{\mathbb{C}}, \mathbb{F}E) \times \mathbb{C}(\mathbb{F}E, \mathbb{F}B) \times \mathbb{C}(\mathbb{F}B, \mathbb{F}A))$$
$$\xrightarrow{\circ_{\mathbb{B}}} \mathbb{B}(1_{\mathbb{B}}, \mathbb{C}(1_{\mathbb{C}}, \mathbb{F}E) \times \mathbb{C}(\mathbb{F}E, \mathbb{F}B) \times \mathbb{C}(\mathbb{F}B, \mathbb{F}A))$$
$$\xrightarrow{\cong} \mathbb{B}(1_{\mathbb{B}}, \mathbb{C}(1_{\mathbb{C}}, \mathbb{F}E)) \times \mathbb{B}(1_{\mathbb{B}}, \mathbb{C}(\mathbb{F}E, \mathbb{F}B)) \times \mathbb{B}(1_{\mathbb{B}}, \mathbb{C}(\mathbb{F}B, \mathbb{F}A)).$$

Let $\langle g_1, g_2, g_3\rangle$ be this morphism. We compute each component as follows. By Proposition 2.57, the morphism $g_1$ is obtained by the following composite:

$$\mathbb{B}(B, A) \xrightarrow{\langle \mathbb{F}_{B,A}, f_1 \circ !\rangle} \mathbb{B}(1_{\mathbb{B}}, E) \times \mathbb{B}(E, \mathbb{C}(1_{\mathbb{C}}, \mathbb{F}E)) \xrightarrow{\circ_{\mathbb{B}}} \mathbb{B}(1_{\mathbb{B}}, \mathbb{C}(1_{\mathbb{C}}, \mathbb{F}E)).$$

Since $f_1 = \mathbb{H}_{\mathbb{F}B,\mathbb{F}A}$, this is rewritten as

$$\mathbb{B}(B, A) \xrightarrow{\mathbb{F}_{B,A}} \mathbb{B}(1_{\mathbb{B}}, E) \xrightarrow{\mathbb{B}\left(1_{\mathbb{B}}, \mathbb{H}_{\mathbb{F}B,\mathbb{F}A}\right)} \mathbb{B}(1_{\mathbb{B}}, \mathbb{C}(1_{\mathbb{C}}, \mathbb{F}E)).$$

By the identity $\gamma_{\mathbb{B}}\mathbb{H} \circ \mathbb{F} = \mathbb{F}^{\#} \circ \mathbb{F}$ in Definition 6.2 (IV*), we have

$$g_1 = \left(\mathbb{B}(B, A) \xrightarrow{\mathbb{F}_{B,A}} \mathbb{B}(1_{\mathbb{B}}, E) \xrightarrow{\mathbb{F}_{1_{\mathbb{B}},E}} \mathbb{B}(1_{\mathbb{B}}, \mathbb{C}(1_{\mathbb{C}}, \mathbb{F}E))\right).$$

Similarly, the morphism $g_2$ is written as

$$\mathbb{B}(B, A) \xrightarrow{\langle \mathbb{F}_{B,A}, f_2 \circ !\rangle} \mathbb{B}(1_{\mathbb{B}}, E) \times \mathbb{B}(E, \mathbb{C}(\mathbb{F}E, \mathbb{F}B)) \xrightarrow{\circ_{\mathbb{B}}} \mathbb{B}(1_{\mathbb{B}}, \mathbb{C}(\mathbb{F}E, \mathbb{F}B)).$$

Substituting the definition of $f_2$ and using the associativity of $\circ_{\mathbb{B}}$, it is rearranged into

$$\mathbb{B}(B,A) \xrightarrow{\langle \mathbb{F}_{B,A}, !_{\mathbb{B}}(E)\circ !, \mathbb{F}(i)\circ ! \rangle} \mathbb{B}(1_{\mathbb{B}},E)\times\mathbb{B}(E,1_{\mathbb{B}})\times\mathbb{B}(1_{\mathbb{B}},\mathbb{C}(\mathbb{F}E,\mathbb{F}B))$$
$$\xrightarrow{\circ_{\mathbb{B}}\times \mathrm{id}} \mathbb{B}(1_{\mathbb{B}},1_{\mathbb{B}})\times\mathbb{B}(1_{\mathbb{B}},\mathbb{C}(\mathbb{F}E,\mathbb{F}B))$$
$$\xrightarrow{\circ_{\mathbb{B}}} \mathbb{B}(1_{\mathbb{B}},\mathbb{C}(\mathbb{F}E,\mathbb{F}B)).$$

Thus, we have

$$g_2 = \left(\mathbb{B}(B,A) \xrightarrow{!} 1 \xrightarrow{\mathbb{F}_{E,B}(i)} \mathbb{B}(1_{\mathbb{B}},\mathbb{C}(\mathbb{F}E,\mathbb{F}B))\right).$$

Finally, the morphism $g_3$ is given by

$$\mathbb{B}(B,A) \xrightarrow{\langle \mathbb{F}_{B,A}, f_3\circ ! \rangle} \mathbb{B}(1_{\mathbb{B}},E)\times\mathbb{B}(E,\mathbb{C}(\mathbb{F}B,\mathbb{F}A)) \xrightarrow{\circ_{\mathbb{B}}} \mathbb{B}(1_{\mathbb{B}},\mathbb{C}(\mathbb{F}B,\mathbb{F}A)).$$

Since $f_3 = e(E) \in_1 \mathbb{B}(E,E)$, it holds that

$$g_3 = \left(\mathbb{B}(B,A) \xrightarrow{\mathbb{F}_{B,A}} \mathbb{B}(1_{\mathbb{B}},\mathbb{C}(\mathbb{F}B,\mathbb{F}A))\right).$$

We return to the morphism $g$. The morphism $g$ is obtained by

$$\mathbb{B}(B,A) \xrightarrow{\langle g_1,g_2,g_3 \rangle} \mathbb{B}(1_{\mathbb{B}},\mathbb{C}(1_{\mathbb{C}},\mathbb{F}E))\times\mathbb{B}(1_{\mathbb{B}},\mathbb{C}(\mathbb{F}E,\mathbb{F}B))\times\mathbb{B}(1_{\mathbb{B}},\mathbb{C}(\mathbb{F}B,\mathbb{F}A))$$
$$\xrightarrow{\circ_{\mathbb{C}}\times \mathrm{id}} \mathbb{B}(1_{\mathbb{B}},\mathbb{C}(1_{\mathbb{C}},\mathbb{F}B))\times\mathbb{B}(1_{\mathbb{B}},\mathbb{C}(\mathbb{F}B,\mathbb{F}A))$$
$$\xrightarrow{\circ_{\mathbb{C}}} \mathbb{B}(1_{\mathbb{B}},\mathbb{C}(1_{\mathbb{C}},\mathbb{F}A)).$$

Since $\mathbb{F}$ preserves the composition morphism, the calculations of $g_1$, $g_2$ and $g_3$ above yield the following commutative diagram:

$$\begin{array}{ccc}
\mathbb{B}(B,A) & \xrightarrow{\langle g_1,g_2,g_3\rangle} & \\
\downarrow{\scriptstyle \langle \mathbb{F}_{B,A}, i\circ !, \mathrm{id}\rangle} & & \\
\mathbb{B}(1_{\mathbb{B}},E)\times\mathbb{B}(E,B)\times\mathbb{B}(B,A) & \xrightarrow{\mathbb{F}_{1_{\mathbb{B}},E}\times\mathbb{F}_{E,B}\times\mathbb{F}_{B,A}} & \mathbb{B}(1_{\mathbb{B}},\mathbb{C}(1_{\mathbb{C}},\mathbb{F}E))\times\mathbb{B}(1_{\mathbb{B}},\mathbb{C}(\mathbb{F}E,\mathbb{F}B)) \times\mathbb{B}(1_{\mathbb{B}},\mathbb{C}(\mathbb{F}B,\mathbb{F}A)) \\
\downarrow{\scriptstyle \circ_{\mathbb{B}}\times\mathrm{id}} & & \downarrow{\scriptstyle \circ_{\mathbb{C}}\times\mathrm{id}} \\
\mathbb{B}(1_{\mathbb{B}},B)\times\mathbb{B}(B,A) & \xrightarrow{\mathbb{F}_{1_{\mathbb{B}},B}\times\mathbb{F}_{B,A}} & \mathbb{B}(1_{\mathbb{B}},\mathbb{C}(1_{\mathbb{C}},\mathbb{F}B))\times\mathbb{B}(1_{\mathbb{B}},\mathbb{C}(\mathbb{F}B,\mathbb{F}A)) \\
\downarrow{\scriptstyle \circ_{\mathbb{B}}} & & \downarrow{\scriptstyle \circ_{\mathbb{C}}} \\
\mathbb{B}(1_{\mathbb{B}},A) & \xrightarrow{\mathbb{F}_{1_{\mathbb{B}},A}} & \mathbb{B}(1_{\mathbb{B}},\mathbb{C}(1_{\mathbb{C}},\mathbb{F}A)).
\end{array}$$

This exactly shows that $g = \mathbb{F}_{1_{\mathbb{B}},A}\circ D_i$ as required. □

Now we show the internal version of Löb's theorem for pre-geminal categories (Theorem 5.4):

**Lemma 6.10 (Löb's theorem for internal pre-geminal categories)** *Let* $(\mathbb{B}, \mathbb{C}, \mathbb{F}, \mathbb{H})$ *be an* $\mathcal{A}$*-internal pre-geminal category. For any global object* $A \in_1 \mathbb{B}$*, there exists a morphism* $Y_A : \mathbb{B}(\mathbb{C}(1_\mathbb{C}, \mathbb{F}A), A) \to \mathbb{B}(1_\mathbb{B}, A)$ *such that the following commutes:*

$$\begin{array}{ccc} \mathbb{B}(\mathbb{C}(1_\mathbb{C}, \mathbb{F}A), A) & \xrightarrow{\langle \mathbb{F}_{1_\mathbb{B},A} \circ Y_A, \mathrm{id}\rangle} & \mathbb{B}(1_\mathbb{B}, \mathbb{C}(1_\mathbb{C}, \mathbb{F}A)) \times \mathbb{B}(\mathbb{C}(1_\mathbb{C}, \mathbb{F}A), A) \\ & \searrow^{Y_A} & \downarrow \circ_\mathbb{B} \\ & & \mathbb{B}(1_\mathbb{B}, A). \end{array}$$

*Proof.* By Proposition 6.3, the quadruple $(\Gamma\mathbb{B}, \mathbb{C}, \Gamma\mathbb{F}, \mathbb{H})$ forms a pre-geminal category. Hence, applying Corollary 5.3, there exists $B \in_1 \mathbb{B}$ such that there is an isomorphism $\mathbb{C}(\mathbb{F}B, \mathbb{F}A) \cong B$ in the category $\Gamma\mathbb{B}$. Letting $E = \mathbb{C}(\mathbb{F}B, \mathbb{F}A) \in_1 \mathbb{B}$, this isomorphism implies that there are global morphisms $i \in_1 \mathbb{B}(E, B)$ and $i^{-1} \in_1 \mathbb{B}(B, E)$ such that the following commute:

$$\begin{array}{ccc} 1 & \xrightarrow{\langle i, i^{-1}\rangle} & \mathbb{B}(E,B) \times \mathbb{B}(B,E) \\ & \searrow^{e(E)} & \downarrow \circ_\mathbb{B} \\ & & \mathbb{B}(E,E), \end{array} \qquad \begin{array}{ccc} 1 & \xrightarrow{\langle i^{-1}, i\rangle} & \mathbb{B}(B,E) \times \mathbb{B}(E,B) \\ & \searrow^{e(B)} & \downarrow \circ_\mathbb{B} \\ & & \mathbb{B}(B,B). \end{array}$$

For this $i \in_1 \mathbb{B}(E, B)$, let $d_i \in_1 \mathbb{B}(E, \mathbb{C}(1_\mathbb{C}, \mathbb{F}A))$ be the global morphism constructed in Lemma 6.9. Further, let $h : \mathbb{B}(\mathbb{C}(1_\mathbb{C}, \mathbb{F}A), A) \to \mathbb{B}(B, A)$ be the following morphism:

$$\begin{aligned} \mathbb{B}(\mathbb{C}(1_\mathbb{C}, \mathbb{F}A), A) &\xrightarrow{\langle i^{-1} \circ !, d_i \circ !, \mathrm{id}\rangle} \mathbb{B}(B, E) \times \mathbb{B}(E, \mathbb{C}(1_\mathbb{C}, \mathbb{F}A)) \times \mathbb{B}(\mathbb{C}(1_\mathbb{C}, \mathbb{F}A), A) \\ &\xrightarrow{\circ_\mathbb{B} \times \mathrm{id}} \mathbb{B}(B, \mathbb{C}(1_\mathbb{C}, \mathbb{F}A)) \times \mathbb{B}(\mathbb{C}(1_\mathbb{C}, \mathbb{F}A), A) \\ &\xrightarrow{\circ_\mathbb{B}} \mathbb{B}(B, A). \end{aligned}$$

We define $Y_A = D_i \circ h$, where $D_i$ is the morphism defined in Lemma 6.9. Explicitly, $Y_A$ is given by the following composite:

$$\begin{aligned} \mathbb{B}(\mathbb{C}(1_\mathbb{C}, \mathbb{F}A), A) &\xrightarrow{\langle \mathbb{F}_{B,A} \circ h, i \circ !, h\rangle} \mathbb{B}(1_\mathbb{B}, E) \times \mathbb{B}(E, B) \times \mathbb{B}(B, A) \\ &\xrightarrow{\circ_\mathbb{B} \times \mathrm{id}} \mathbb{B}(1_\mathbb{B}, B) \times \mathbb{B}(B, A) \\ &\xrightarrow{\circ_\mathbb{B}} \mathbb{B}(1_\mathbb{B}, A). \end{aligned}$$

We now verify that the required diagram commutes. First, by the associativity of $\circ_\mathbb{B}$ and the fact that $i$ and $i^{-1}$ are inverses of each other, the following diagram commutes:

$$\begin{array}{ccc} \mathbb{B}(\mathbb{C}(1_\mathbb{C}, \mathbb{F}A), A) & \xrightarrow{\langle i \circ !, h\rangle} & \mathbb{B}(E, B) \times \mathbb{B}(B, A) \\ {\scriptstyle \langle d_i \circ !, \mathrm{id}\rangle} \downarrow & & \downarrow \circ_\mathbb{B} \\ \mathbb{B}(E, \mathbb{C}(1_\mathbb{C}, \mathbb{F}A)) \times \mathbb{B}(\mathbb{C}(1_\mathbb{C}, \mathbb{F}A), A) & \xrightarrow{\circ_\mathbb{B}} & \mathbb{B}(E, A). \end{array}$$

Using the associativity of $\circ_{\mathbb{B}}$ again, $Y_A$ can be rewritten as

$$\begin{aligned}\mathbb{B}(\mathbb{C}(1_{\mathbb{C}}, \mathbb{F}A), A) &\xrightarrow{\langle \mathbb{F}_{B,A} \circ h,\, d_i \circ !,\, \mathrm{id} \rangle} \mathbb{B}(1_{\mathbb{B}}, E) \times \mathbb{B}(E, \mathbb{C}(1_{\mathbb{C}}, \mathbb{F}A)) \times \mathbb{B}(\mathbb{C}(1_{\mathbb{C}}, \mathbb{F}A), A) \\ &\xrightarrow{\circ_{\mathbb{B}} \times \mathrm{id}} \mathbb{B}(1_{\mathbb{B}}, \mathbb{C}(1_{\mathbb{C}}, \mathbb{F}A)) \times \mathbb{B}(\mathbb{C}(1_{\mathbb{C}}, \mathbb{F}A), A) \\ &\xrightarrow{\circ_{\mathbb{B}}} \mathbb{B}(1_{\mathbb{B}}, A).\end{aligned}$$

On the other hand, by the construction of $d_i$, the following diagram commutes:

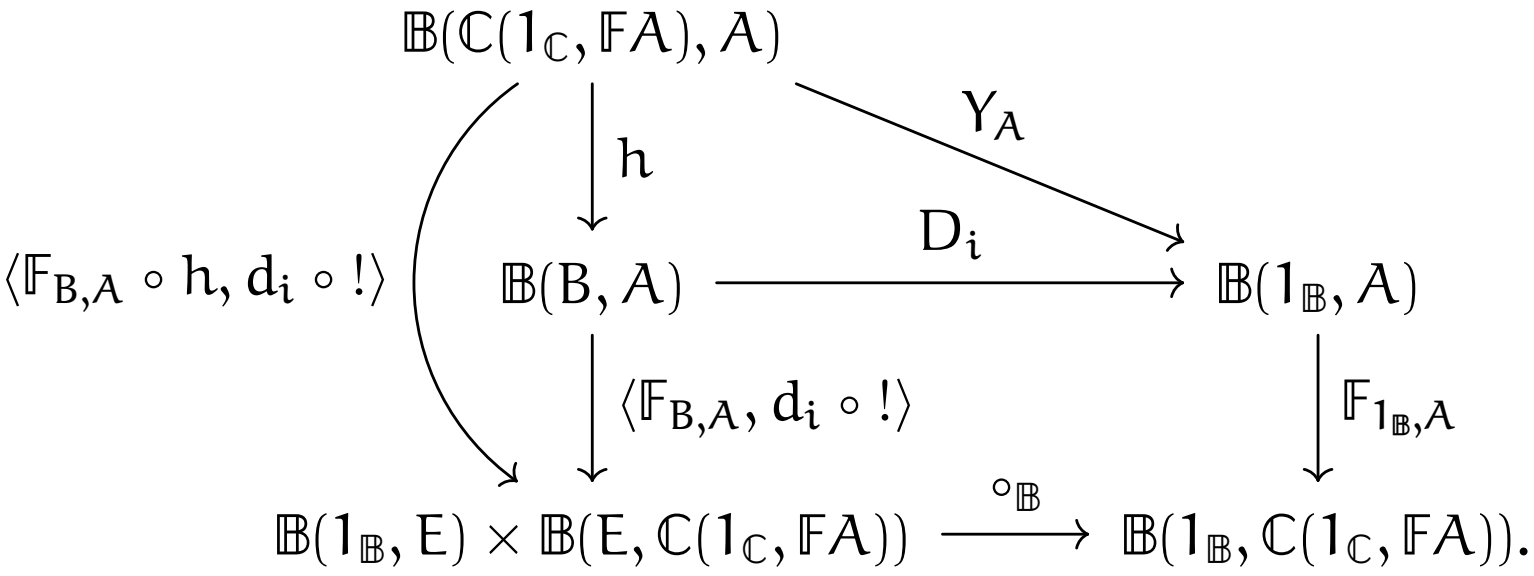

Combining these, it follows that $Y_A$ coincides with the composite

$$\mathbb{B}(\mathbb{C}(1_{\mathbb{C}}, \mathbb{F}A), A) \xrightarrow{\langle \mathbb{F}_{1_{\mathbb{B}},A} \circ Y_A,\, \mathrm{id} \rangle} \mathbb{B}(1_{\mathbb{B}}, \mathbb{C}(1_{\mathbb{C}}, \mathbb{F}A)) \times \mathbb{B}(\mathbb{C}(1_{\mathbb{C}}, \mathbb{F}A), A) \xrightarrow{\circ_{\mathbb{B}}} \mathbb{B}(1_{\mathbb{B}}, A).$$

This shows that the required diagram commutes. □

Applying this lemma to the internal pre-geminal category induced by a geminal category, we obtain the Gödel–Löb axiom for geminal categories.

> **Theorem 6.11 (Gödel–Löb axiom for geminal categories)** *Let $(\mathcal{B}, \mathbb{C}, F, \mathbb{H})$ be a geminal category. For any object $A \in \mathcal{B}$, there exists a morphism $Y_A : \mathbb{C}(F\Box A, FA) \to \Box A$ such that the following commutes:*
>
> $$\begin{array}{ccc} \mathbb{C}(F\Box A, FA) & \xrightarrow{\langle \mathbb{H}_{1_{\mathbb{C}},FA} \circ Y_A,\, \mathrm{id} \rangle} & \mathbb{C}(1_{\mathbb{C}}, F\Box A) \times \mathbb{C}(F\Box A, FA) \\ & \searrow^{Y_A} & \downarrow \circ_{\mathbb{C}} \\ & & \mathbb{C}(1_{\mathbb{C}}, FA). \end{array}$$

*Proof.* First, note that the codomain of $\mathbb{H}_{1_{\mathbb{C}},FA}$ is

$$(\Box\mathbb{C})(\mathbb{H}1_{\mathbb{C}}, \mathbb{H}FA) \cong \Box(\mathbb{C}(1_{\mathbb{C}}, FA)) = \Box\Box A = \mathbb{C}(1_{\mathbb{C}}, F\Box A),$$

where the first isomorphism is obtained by an argument similar to one in the proof of Lemma 4.16. Therefore, the upper path in the given diagram is well-defined.

The morphism $Y_A$ is obtained by applying Lemma 6.10 to the $\mathcal{B}$-internal pre-geminal category $(\mathbb{C}, F\mathbb{C}, \mathbb{H}, F\mathbb{H})$ constructed in Theorem 6.8 and the global object $FA \in_1 \mathbb{C}$. Indeed, the lemma yields the morphism

$$Y_A : \mathbb{C}((F\mathbb{C})(1_{F\mathbb{C}}, \mathbb{H}FA), FA) \to \mathbb{C}(1_{\mathbb{C}}, FA).$$

Since $\ulcorner \mathbb{H} \circ F = F^{\#} \circ F$, we have

$$(F\mathbb{C})(1_{F\mathbb{C}}, \mathbb{H}FA) = (F\mathbb{C})\big(F^{\#}(1_{\mathbb{C}}), F^{\#}FA\big) \cong F(\mathbb{C}(1_{\mathbb{C}}, FA)) = F\Box A.$$

Hence, $Y_A$ can be viewed as a morphism of the form $\mathbb{C}(F\Box A, FA) \to \Box A$. The commutativity of the required diagram directly follows from the one established in Lemma 6.10. □

The morphism $Y_A : \mathbb{C}(F\Box A, FA) \to \Box A$ is a categorical analog of the Gödel–Löb axiom, $\Box(\Box A \to A) \to \Box A$. This interpretation is justified by the fact that the object $\mathbb{C}(F\Box A, FA)$ is the internalization of the hom-set $\mathcal{B}(\Box A, A)$, thereby corresponding to the formula $\Box(\Box A \to A)$. Furthermore, in the context of intensional recursion, $Y_A$ can be regarded as a *modal* or *intensional fixed point combinator* — we refer the reader once again to Kavvos [21] for a detailed discussion.

**Remark 6.12** This interpretation of the axiom suggests that exponentials are not required at all to interpret the Gödel–Löb axiom in categorical structures, even though the axiom appears to rely heavily on implications. This is explained by the fact that the object $\mathbb{C}(FB, FA)$ internalizes the *judgment* $B \vdash A$ rather than the *implication* $B \to A$. From this perspective, the Gödel–Löb axiom may be more appropriately understood as $\Box(\Box A \vdash A) \vdash A$ or $\Box([\Box A]A) \vdash A$ (the latter in the notation of *contextual modal type theory* [33]), both of which involve neither implications nor function types.

**Remark 6.13** There is an alternative strategy for proving the Gödel–Löb axiom that avoids a full internalization of the proof of Löb's theorem. In the modal logic K4, the Gödel–Löb axiom $\Box(\Box A \to A) \to \Box A$ is famously derivable from Löb's rule [3, p. 59]. Specifically, this is achieved by applying Löb's rule to the Gödel–Löb axiom itself. Since the axiom 4 is attainable in geminal categories (see Theorem 6.16), this approach may yield a more elegant proof of Theorem 6.11.

While this strategy is directly applicable to geminal categories with exponentials, a challenge arises in the general case: there is no single object in a geminal category corresponding to the Gödel–Löb axiom $\Box(\Box A \to A) \to \Box A$. Consequently, the standard Löb's theorem (Theorem 5.4) cannot be applied to the axiom itself, suggesting that a "parameterized" version of the theorem would be needed to implement this strategy. We leave the formal development of this parameterized version and its application to the Gödel–Löb axiom as future work, stating it here as a conjecture:

**Conjecture 6.14 (Parameterized Löb's theorem for geminal categories)** *Let* $(\mathcal{B}, \mathbb{C}, F, \mathbb{H})$ *be a geminal category. For any objects* $A, B \in \mathcal{B}$ *and any morphism* $f : \mathbb{C}(FB, FA) \times B \to A$, *there exists a morphism* $a : B \to A$ *such that the following commutes:*

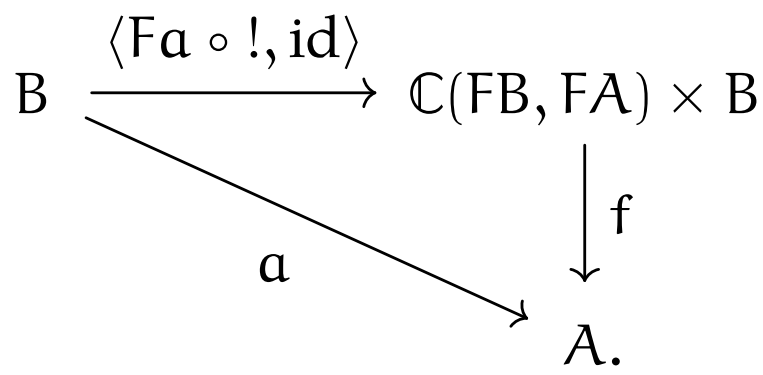

### 6.3 Modal structures

Since we now have the Gödel–Löb axiom for geminal categories, it is natural to ask whether the pair $(\mathcal{B}, \Box)$ forms a categorical model of the provability logic GL, or its corresponding type theory. Indeed, in the context of modal calculi, several categorical models are known to interpret modality via endofunctors preserving finite products [20, 34]. In this subsection, we examine the structure of $(\mathcal{B}, \Box)$ arising from geminal categories in a manner similar to these existing models. It turns out that $(\mathcal{B}, \Box)$ possesses a structure remarkably similar to Kavvos' *Gödel–Löb category* [20], which is known to be a categorical model of the dual-context modal calculus DGL. However, we also observe a subtle difference between the two, at least in their natural formulations.

We start with a simple lemma.

**Lemma 6.15** *Let $\mathcal{B}$ be a category with finite limits, $\mathbb{C}$ be a $\mathcal{B}$-internal category with chosen finite limits, and $F : \mathcal{B} \to \Gamma\mathbb{C}$ be a functor preserving finite limits. For any morphism $f : A \times B \to C$ in $\mathcal{B}$, one can construct $f^\vee : \Box A \to \mathbb{C}(FB, FC)$ such that the following commutes:*

$$\begin{array}{ccc} \Box B \times \Box A & \xrightarrow{\mathrm{id} \times f^\vee} & \mathbb{C}(1_\mathbb{C}, FB) \times \mathbb{C}(FB, FC) \\ \cong \downarrow & & \downarrow \circ_\mathbb{C} \\ \Box(A \times B) & \xrightarrow{\Box f} & \mathbb{C}(1_\mathbb{C}, FC). \end{array}$$

*Proof.* First, let $g$ be the composite

$$\Box A = \mathbb{C}(1_\mathbb{C}, FA) \xrightarrow{\langle !_\mathbb{C}(FB) \circ !, \mathrm{id} \rangle} \mathbb{C}(FB, 1_\mathbb{C}) \times \mathbb{C}(1_\mathbb{C}, FA) \xrightarrow{\circ_\mathbb{C}} \mathbb{C}(FB, FA),$$

and let $h$ be the composite

$$\Box A \xrightarrow{\langle g, e(FB) \circ ! \rangle} \mathbb{C}(FB, FA) \times \mathbb{C}(FB, FB) \xrightarrow{\cong} \mathbb{C}(FB, F(A \times B)).$$

We define $f^\vee$ by the following composite:

$$\Box A \xrightarrow{\langle h, Ff \circ ! \rangle} \mathbb{C}(FB, F(A \times B)) \times \mathbb{C}(F(A \times B), FC) \xrightarrow{\circ_\mathbb{C}} \mathbb{C}(FB, FC).$$

Next, we verify the required commutativity. By the definition of $f^\vee$ and the associativity of $\circ_\mathbb{C}$, the following diagram commutes:

$$\begin{array}{ccc} \Box B \times \Box A & \xrightarrow{\mathrm{id} \times f^\vee} & \mathbb{C}(1_\mathbb{C}, FB) \times \mathbb{C}(FB, FC) \\ \mathrm{id} \times h \downarrow & & \\ \mathbb{C}(1_\mathbb{C}, FB) \times \mathbb{C}(FB, F(A \times B)) & & \downarrow \circ_\mathbb{C} \\ \langle \circ_\mathbb{C}, Ff \circ ! \rangle \downarrow & & \\ \mathbb{C}(1_\mathbb{C}, F(A \times B)) \times \mathbb{C}(F(A \times B), FC) & \xrightarrow{\circ_\mathbb{C}} & \mathbb{C}(1_\mathbb{C}, FC). \end{array}$$

We compute the first component of the left-hand side. By the definition of $h$ and Proposition 2.57, the following commutes:

$$\begin{array}{ccc}
\Box B \times \Box A & \xrightarrow{\mathrm{id} \times h} & \mathbb{C}(1_{\mathbb{C}}, FB) \times \mathbb{C}(FB, F(A \times B)) \\
\downarrow{\scriptstyle \langle \mathrm{id} \times (e(FB) \circ !), \mathrm{id} \times g \rangle} & & \\
(\mathbb{C}(1_{\mathbb{C}}, FB) \times \mathbb{C}(FB, FB)) \times (\mathbb{C}(1_{\mathbb{C}}, FB) \times \mathbb{C}(FB, FA)) & & \downarrow{\scriptstyle \circ_{\mathbb{C}}} \\
\downarrow{\scriptstyle \circ_{\mathbb{C}} \times \circ_{\mathbb{C}}} & & \\
\mathbb{C}(1_{\mathbb{C}}, FB) \times \mathbb{C}(1_{\mathbb{C}}, FA) & \xrightarrow{\cong} & \mathbb{C}(1_{\mathbb{C}}, F(A \times B)).
\end{array}$$

The left-hand side of this diagram coincides with the identity morphism, as shown by the following commutative diagrams:

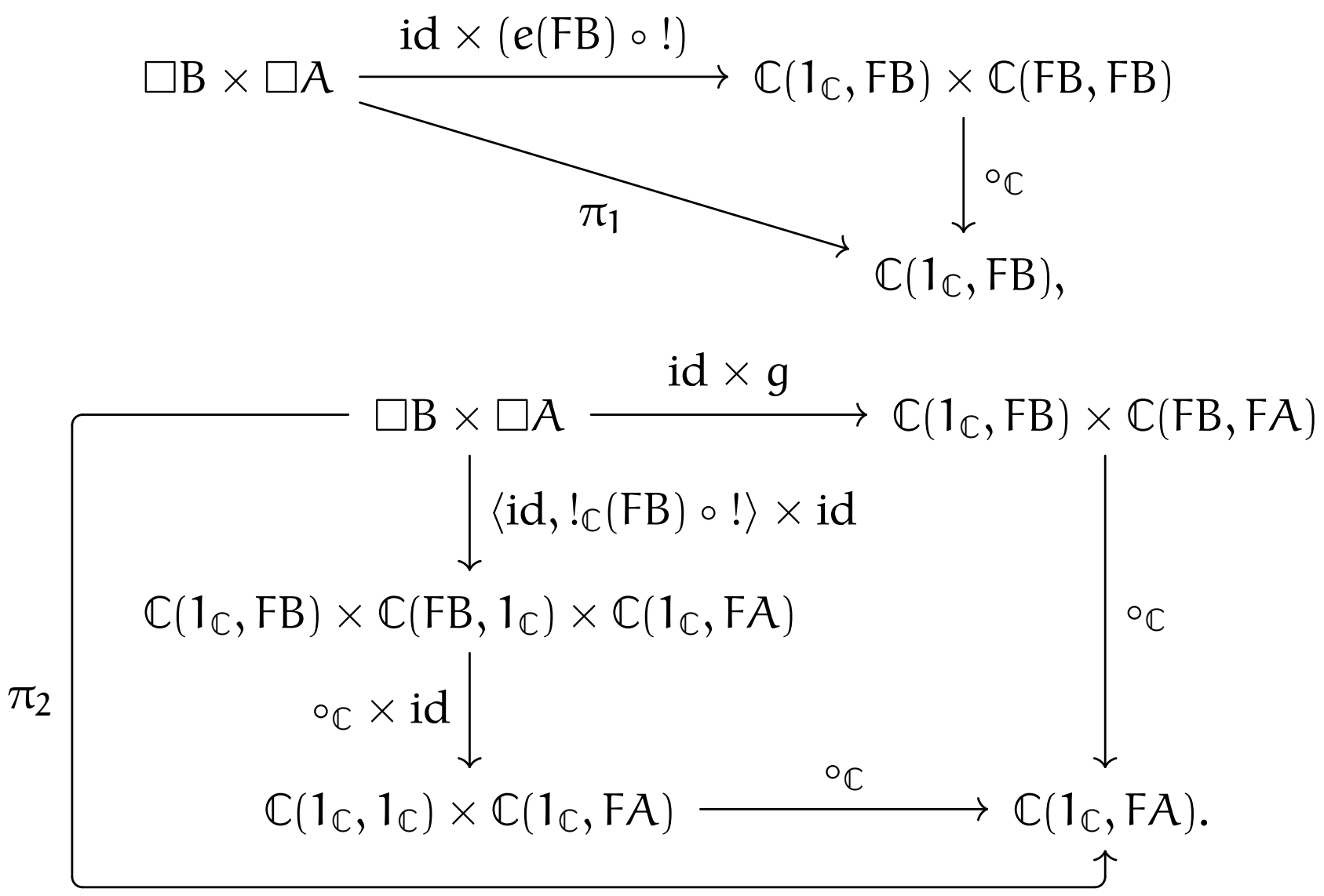

Hence, we have

$$\begin{array}{ccc}
\Box B \times \Box A & \xrightarrow{\mathrm{id} \times h} & \mathbb{C}(1_{\mathbb{C}}, FB) \times \mathbb{C}(FB, F(A \times B)) \\
& \searrow^{\cong} & \downarrow{\scriptstyle \circ_{\mathbb{C}}} \\
& & \Box(A \times B).
\end{array}$$

Therefore, we obtain the following commutative diagram:

$$\begin{array}{ccc}
\Box B \times \Box A & \xrightarrow{\mathrm{id} \times f^{\vee}} & \mathbb{C}(1_{\mathbb{C}}, FB) \times \mathbb{C}(FB, FC) \\
\downarrow{\scriptstyle \cong} & & \\
\Box(A \times B) & & \downarrow{\scriptstyle \circ_{\mathbb{C}}} \\
\downarrow{\scriptstyle \langle \mathrm{id}, Ff \circ ! \rangle} & & \\
\mathbb{C}(1_{\mathbb{C}}, F(A \times B)) \times \mathbb{C}(F(A \times B), FC) & \xrightarrow{\circ_{\mathbb{C}}} & \mathbb{C}(1_{\mathbb{C}}, FC).
\end{array}$$

On the other hand, by the definition of $\Box$, the following diagram commutes:

$$\begin{array}{ccc} \Box(A \times B) & \xrightarrow{\langle \mathrm{id}, Ff \circ ! \rangle} & \mathbb{C}(1_{\mathbb{C}}, F(A \times B)) \times \mathbb{C}(F(A \times B), FC) \\ & \searrow^{\Box f} & \downarrow \circ_{\mathbb{C}} \\ & & \mathbb{C}(1_{\mathbb{C}}, FC). \end{array}$$

Combining these, we obtain the required commutative diagram. $\Box$

Now we describe the modal structure of the pair $(\mathcal{B}, \Box)$.

**Theorem 6.16** *Let $(\mathcal{B}, \mathbb{C}, F, \mathbb{H})$ be a geminal category. Then, the pair $(\mathcal{B}, \Box)$ can be equipped with the following additional structures:*

- *A natural transformation $\delta : \Box \Rightarrow \Box\Box$ satisfying the co-associative law, i.e., the identity $(\delta\Box) \circ \delta = (\Box\delta) \circ \delta$.*

$$\Box \xrightarrow{\delta} \Box\Box \underset{\Box\delta}{\overset{\delta\Box}{\rightrightarrows}} \Box\Box\Box$$

- *A family of functions $(-)^{\dagger}_{A,B} : \mathcal{B}(B \times \Box A, A) \to \mathcal{B}(\Box B, \Box A)$ which is natural in $B \in \mathcal{B}$ and such that, for any $f : B \times \Box A \to A$, the following diagram commutes:*

$$\begin{array}{ccccc} \Box B & \xrightarrow{\langle \mathrm{id}, \delta_A \circ f^{\dagger} \rangle} & \Box B \times \Box\Box A & \xrightarrow{\cong} & \Box(B \times \Box A) \\ & & \searrow^{f^{\dagger}} & & \downarrow \Box f \\ & & & & \Box A. \end{array}$$

*Proof.* For any $A \in \mathcal{B}$, we define $\delta_A$ by

$$\delta_A = \mathbb{H}_{1_{\mathbb{C}}, FA} : \Box A = \mathbb{C}(1_{\mathbb{C}}, FA) \to \Box\Box A = \Box(\mathbb{C}(1_{\mathbb{C}}, FA)).$$

The naturality of $\delta$ follows from the functoriality of $\mathbb{H}$. To verify the required equation, recall the identity $\mathbb{H}^{\#} \circ \mathbb{H} = \Box\mathbb{H} \circ \mathbb{H}$ imposed on geminal categories.

$$\mathbb{C} \xrightarrow{\mathbb{H}} \Box\mathbb{C} \underset{\Box\mathbb{H}}{\overset{\mathbb{H}^{\#}}{\rightrightarrows}} \Box\Box\mathbb{C}$$

Restricting this identity to the hom-object $\mathbb{C}(1_{\mathbb{C}}, FA)$, we obtain the identity expressed in the following diagram:

$$\mathbb{C}(1_{\mathbb{C}}, FA) \xrightarrow{\mathbb{H}_{1_{\mathbb{C}}, FA}} \Box(\mathbb{C}(1_{\mathbb{C}}, FA)) \underset{\Box\mathbb{H}_{1_{\mathbb{C}}, FA}}{\overset{\mathbb{H}_{1_{\mathbb{C}}, \mathbb{C}(1_{\mathbb{C}}, FA)}}{\rightrightarrows}} \Box\Box(\mathbb{C}(1_{\mathbb{C}}, FA)).$$

This precisely shows the co-associativity law, $(\delta\Box)_A \circ \delta_A = (\Box\delta)_A \circ \delta_A$.

Next, we define the operation $(-)^{\dagger}_{A,B}$. For any morphism $f : B \times \Box A \to A$, we define $f^{\dagger} : \Box B \to \Box A$ as the following composite:

$$\Box B \xrightarrow{f^{\vee}} \mathbb{C}(F\Box A, FA) \xrightarrow{Y_A} \Box A,$$

where $f^{\vee}$ and $Y_A$ are the morphisms constructed in Lemma 6.15 and Theorem 6.11, respectively. The required commutativity condition is verified by the following commutative diagram:

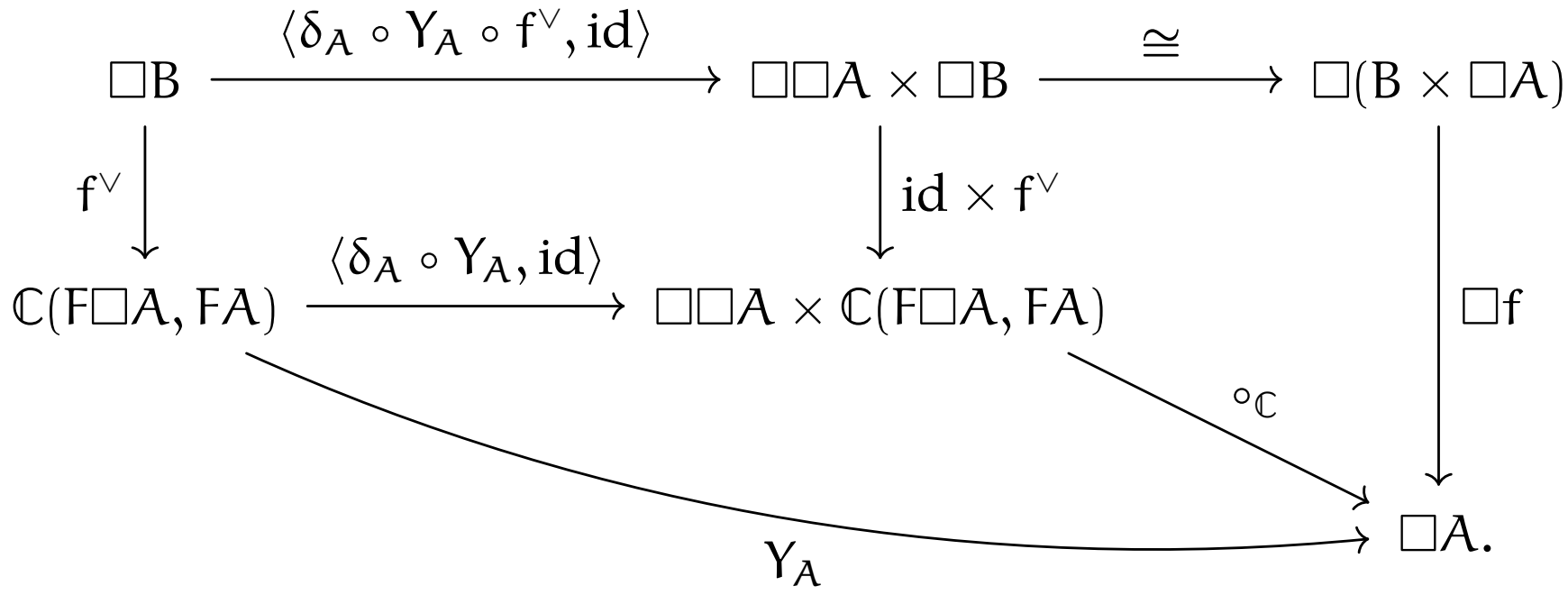

This completes the proof. $\square$

The first three components $(\mathcal{B}, \Box, \delta)$ described above provide a categorical counterpart to the modal logic K4. Indeed, this structure precisely matches Kavvos' *Kripke-4 categories* [20], with the exception that our conditions do not require $\mathcal{B}$ to be Cartesian closed. The final component $(-)^{\dagger}_{A,B}$ corresponds to the following rule admissible in GL:

$$\frac{\Gamma, \Box A \vdash A}{\Box\Gamma \vdash \Box A.}$$

In the Cartesian closed case, this rule precisely corresponds to the Gödel–Löb axiom, as shown in the following proposition.

**Proposition 6.17** *Let $(\mathcal{B}, \Box, \delta)$ be a triple constituting the first three components in Theorem 6.16. Assuming that $\mathcal{B}$ is Cartesian closed, there is a one-to-one correspondence between the following:*

- *Families of functions $(-)^{\dagger}_{A,B} : \mathcal{B}(B \times \Box A, A) \to \mathcal{B}(\Box B, \Box A)$ that are natural in $B$ and satisfy the commutativity condition in Theorem 6.16.*
- *Families of morphisms $Y_A : \Box(A^{\Box A}) \to \Box A$ such that the following diagram commutes:*

$$\begin{array}{ccccc} \Box(A^{\Box A}) & \xrightarrow{\langle \mathrm{id}, \delta_A \circ Y_A \rangle} & \Box(A^{\Box A}) \times \Box\Box A & \xrightarrow{\cong} & \Box(A^{\Box A} \times \Box A) \\ & \searrow^{Y_A} & & & \downarrow{\Box \mathrm{ev}} \\ & & & & \Box A. \end{array}$$

*Proof.* Fix an object $A \in \mathcal{B}$. The correspondence is given by the Yoneda lemma as follows:

$$\dfrac{\dfrac{\mathcal{B}((-) \times \Box A, A) \Rightarrow \mathcal{B}(\Box(-), \Box A)}{\mathcal{B}\big((-), A^{\Box A}\big) \Rightarrow \mathcal{B}(\Box(-), \Box A)}}{\Box\big(A^{\Box A}\big) \to \Box A.}$$

We show that the commutativity conditions for $(-)^{\dagger}_{A,B}$ and $Y_A$ are equivalent under this correspondence. If the family $(-)^{\dagger}_{A,B}$ satisfies the commutativity condition in Theorem 6.16, then the required commutativity for $Y_A = (\mathrm{ev})^{\dagger}$ follows immediately.

Conversely, assume that $Y_A$ satisfies the given condition. For any $f : B \times \Box A \to A$, its corresponding morphism $f^{\dagger} : \Box B \to \Box A$ is defined by the following composite:

$$\Box B \xrightarrow{\Box(\lambda(f))} \Box\big(A^{\Box A}\big) \xrightarrow{Y_A} \Box A.$$

For this morphism, we have the following commutative diagram:

$$\begin{array}{ccccc}
\Box B & \xrightarrow{\langle \mathrm{id}, \delta_A \circ Y_A \circ \Box(\lambda(f)) \rangle} & \Box B \times \Box\Box A & \xrightarrow{\cong} & \Box(B \times \Box A) \\
\downarrow{\scriptstyle \Box(\lambda(f))} & & \downarrow{\scriptstyle \Box(\lambda(f)) \times \mathrm{id}} & & \downarrow{\scriptstyle \Box(\lambda(f) \times \mathrm{id})} \\
\Box\big(A^{\Box A}\big) & \xrightarrow{\langle \mathrm{id}, \delta_A \circ Y_A \rangle} & \Box\big(A^{\Box A}\big) \times \Box\Box A & \xrightarrow{\cong} & \Box\big(A^{\Box A} \times \Box A\big) \\
& & & & \downarrow{\scriptstyle \Box\,\mathrm{ev}} \\
& & & & \Box A
\end{array}$$

This shows that $(-)^{\dagger}$ satisfies the required commutativity condition in Theorem 6.16. $\Box$

Finally, we compare our structure in Theorem 6.16 with existing categorical models for modal calculi. The closest structure is Kavvos' *Gödel–Löb category* [20]. It provides a categorical model for the dual-context modal calculus DGL, a type-theoretic counterpart of the provability logic GL. However, there are three differences to consider.

First, while Gödel–Löb categories are assumed to be Cartesian closed to model the modal $\lambda$-calculus, we do not require this condition. This is a minor distinction, as Kavvos' definition naturally extends to any category with finite products, as seen in *guarded fixpoint categories* [30].

Second, the operation $(-)^{\dagger}_{A,B}$ in a Gödel–Löb category has a slightly different domain, $\mathcal{B}(\Box B \times B \times \Box A, A)$. Our simpler domain $\mathcal{B}(B \times \Box A, A)$ corresponds to focusing only on the structures induced by *modal fixed point combinators* $Y_A$ [20, Definition 5.30].

The third one is the most subtle but also appears to be the most serious: in the definition of a Gödel–Löb category, the morphism $\delta_A \circ f^{\dagger}$ in the commutativity condition in Theorem 6.16 is replaced by $\big(\Box f^{\dagger}\big) \circ \delta_B$. It is currently unclear whether this results in an essentially different structure. This mismatch appears to be deeply rooted in the respective natures of geminal categories and DGL: while the former naturally internalizes the fixed point property of Löb's theorem as $\delta_A \circ f^{\dagger}$, the latter's equational theory naturally leads to $\big(\Box f^{\dagger}\big) \circ \delta_B$. For this reason, we leave to future work the question whether geminal categories precisely model Kavvos' dual-context modal calculus DGL.

The notion of a *guarded fixpoint category* [30], logically corresponding to the *strong Löb logic*, has been proposed as a model for guarded recursion. From any guarded fixpoint category where $\Box$ preserves finite products, one can induce the structure described in Theorem 6.16 in the same manner as Kavvos' construction of Gödel–Löb categories from guarded fixpoint categories [20, Theorem 5.38]. In particular, any *Cartesian traced category* induces the structure in Theorem 6.16 with $\Box$ taken as the identity functor, owing to the established equivalence between traces and fixed point operators [11].

# 7 Conclusion and future work

In this thesis, we have developed the theory of *geminal categories*, originally introduced by Ramesh [38], in a self-contained and reorganized manner. By introducing the notion of a *code structure* on fibrations as an abstraction of Gödel coding, we have successfully re-motivated the framework of geminal categories. Within this setting, we have established the *Gödel–Löb axiom for geminal categories* and a simplified proof of *Löb's theorem for pre-geminal categories*. These results would provide essential steps for establishing the connection between geminal categories and the provability logic GL. Our reorganization makes the theory of geminal categories more accessible to a broader audience, and contributes to the foundational theory of geminal categories as a promising tool for analyzing self-internalizing structures.

Since the theory of geminal categories has recently emerged, there is plenty of room for further development. First, establishing that geminal categories provide categorical models for type-theoretic counterparts of GL, such as Kavvos' DGL, is the most natural way to bridge the gap between our framework and traditional provability logics, as well as their associated computational systems. Conversely, it is also worth considering which categorical models of GL arise from geminal (i.e., "self-internalizing") structures. Notably, Remark 6.12 suggests a full correspondence between geminal categories and systems without function types; to explore this, *contextual modal type theory* [33] that internalizes contexts or judgments would be useful. Furthermore, developing an enriched analog of geminal categories (Remark 4.13) would not only show its generality but may also relate to the *iterated enrichment* approach to modality by Nishiwaki et al. [34].

The notion of *introspective theories*, also introduced by Ramesh [38], remains to be investigated. This concept is deeply interrelated with geminal categories as it captures the "internalizing" process such as the one demonstrated in Section 6. Since performing such internalization is often technically demanding, a formal syntactic treatment would be necessary. It is plausible that appropriate internal languages for introspective theories could be directly related to calculi for metaprogramming.

For geminal categories to play a universal role, they must encompass various examples from both proof theory and computer science. While we treated only initial geminal models as examples, Ramesh's construction based on presheaf categories [38, Section 6.5] deserves further analysis. Furthermore, as guarded fixpoint categories [30] contain rich models such as complete metric spaces and complete partial orders, it should be investigated whether these well-known models of computation can be induced from geminal categories. From a proof-theoretic perspective, it would also be worth considering how more advanced logics such as the polymodal logic GLP [15] can be situated within the context of geminal categories.